# General linear and functor cohomology over finite fields

By Vincent Franjou*, Eric M. Friedlander**, Alexander Scorichenko,
and Andrei Suslin***

### Introduction

In recent years, there has been considerable success in computing Ext-groups of modular representations associated to the general linear group by relating this problem to one of computing Ext-groups in functor categories [F-L-S], [F-S]. In this paper, we extend our ability to make such Ext-group calculations by establishing several fundamental results. Throughout this paper, we work over fields of positive characteristic $p$.

The reader familiar with the representation theory of algebraic objects will recognize the importance of an understanding of Ext-groups. For example, the existence of nonzero Ext-groups of positive degree is equivalent to the existence of objects which are not "direct sums" of simple objects. Indeed, a knowledge of Ext-groups provides considerable knowledge of compound objects. In the study of modular representation theory of finite Chevalley groups such as $\mathrm{GL}_n(F_q)$, Ext-groups play an even more central role: it has been shown in [CPS] that a knowledge of certain $\mathrm{Ext}^1$-groups is sufficient to prove Lusztig's Conjecture concerning the dimension and characters of irreducible representations.

We consider two different categories of functors, the category $\mathcal{F}(\mathbb{F}_q)$ of all functors from finite dimensional $\mathbb{F}_q$-vector spaces to $\mathbb{F}_q$-vector spaces, where $\mathbb{F}_q$ is the finite field of cardinality $q$, and the category $\mathcal{P}(\mathbb{F}_q)$ of strict polynomial functors of finite degree as defined in [F-S]. The category $\mathcal{P}(\mathbb{F}_q)$ presents several advantages over the category $\mathcal{F}(\mathbb{F}_q)$ from the point of view of computing Ext-groups. These are the accessibility of injectives and projectives, the existence of a base change, and an even easier access to Ext-groups of tensor products. This explains the usefulness of our comparison in Theorem 3.10 of Ext-groups in the category $\mathcal{P}(\mathbb{F}_q)$ with Ext-groups in the category $\mathcal{F}(\mathbb{F}_q)$. Weaker forms of this theorem have been known to us since 1995 and to S. Betley independently

*Partially supported by the C.N.R.S., UMR 6629.
**Partially supported by the N.S.F., N.S.A., and the Humboldt Foundation.
***Partially supported by the N.S.F. grant DMS-9510242.

(see its use in [B2, §3]). This early work apparently inspired the paper of N. Kuhn [K2] as well as the present paper.

The calculation (for an arbitrary field $k$ of positive characteristic) of the $\text{Ext}_{\mathcal{P}(k)}$-groups from a Frobenius-twisted divided power functor to a Frobenius-twisted symmetric power functor is presented in Theorem 4.5. This calculation is extended in Section 5 to various other calculations of $\text{Ext}_{\mathcal{P}(k)}$-groups between divided powers, exterior powers, and symmetric powers. This leads to similar $\text{Ext}_{\mathcal{F}(\mathbb{F}_q)}$-group calculations in Theorem 6.3. The results are given with their structure as tri-graded Hopf algebras. A result, stated in this form, has been obtained by N. Kuhn for natural transformations ($\text{Hom}_{\mathcal{F}(\mathbb{F}_q)}$-case) from a divided power functor to a symmetric power functor [K3, §5]. Computing Ext-groups between symmetric and exterior powers is the topic in [F], which contains partial results for the category $\mathcal{F}(\mathbb{F}_p)$.

The final result, proved by the last-named author in the appendix, is the proof of equality of the $\text{Ext}_{\mathcal{F}(\mathbb{F}_q)}$-groups and $\text{Ext}_{\text{GL}_\infty(\mathbb{F}_q)}$-groups associated to finite functors. This has a history of its own that will be briefly recalled at the beginning of the appendix.

These results complement the important work of E. Cline, B. Parshall, L. Scott, and W. van der Kallen in [CPSvdK]. The results in that paper apply to general reductive groups defined and split over the prime field and to general (finite dimensional) rational modules, but lack the computational applicability of the present paper. As observed in [CPSvdK], Ext-groups of rational $G$-modules are isomorphic to Ext-groups of associated Chevalley groups, provided that Frobenius twist is applied sufficiently many times to the rational modules, and provided that the finite field is sufficiently large. One consequence of our work is a strong stability result for the effect of iterating Frobenius twists (Corollary 4.10); it applies to the Ext-groups of rational modules arising from strict polynomial functors of finite degree. A second consequence is an equally precise lower bound for the order of the finite field $\mathbb{F}_q$ required to compare these "stably twisted" rational Ext-groups with the Ext-groups computed for the infinite general linear group $\text{GL}_\infty(\mathbb{F}_q)$. For explicit calculations of Ext-groups for $\text{GL}_n(\mathbb{F}_q)$ for various fundamental $\text{GL}_n(\mathbb{F}_q)$-modules, one then can combine results of this paper with explicit stability results of W. van der Kallen [vdK].

What follows is a brief sketch of the contents of this paper. Section 1 recalls the category $\mathcal{P}(k)$ of strict polynomial functors of finite degree on finite dimensional $k$-vector spaces and further recalls the relationship of $\mathcal{P}$ to the category of rational representations of the general linear group. The investigation of the forgetful functor $\mathcal{P}(k) \to \mathcal{F}(k)$ is begun by observing that $\text{Hom}_{\mathcal{P}(k)}(P, Q) = \text{Hom}_{\mathcal{F}(k)}(P, Q)$ for strict polynomial functors of degree $d$ provided that cardinality of the field $k$ is at least $d$. Theorem 1.7 presents the key property for Ext-groups involving functors of exponential type; this property enables our Ext-group calculations and appears to have no evident analogue for $GL_n$ (with $n$ finite).

In Section 2, we employ earlier work of N. Kuhn and the first named author to provide a first comparison of $\mathrm{Ext}^n_{\mathcal{P}(\mathbb{F}_q)}(P,Q)$ and $\mathrm{Ext}^n_{\mathcal{F}(\mathbb{F}_q)}(P,Q)$. The weakness of this comparison is that the lower bound on $q$ depends upon the Ext-degree $n$ as well as the degrees of $P$ and $Q$. Section 3 remedies this weakness, providing in Theorem 3.10 a comparison of Ext-groups for which the lower bound for $q$ depends only upon the degrees of $P$ and $Q$. Our proof relies heavily upon an analysis of base change, i.e. the effect of an extension $L/k$ of finite fields on $\mathrm{Ext}_{\mathcal{F}}(P,Q)$-groups when $P$ and $Q$ are strict polynomial functors.

For every integer $d$, we compute in Theorem 4.5, $\mathrm{Ext}^*_{\mathcal{P}(k)}(\Gamma^{d(r)}, S^{dp^{r-j}(j)})$ and $\mathrm{Ext}^*_{\mathcal{P}(k)}(\Gamma^{d(r)}, \Lambda^{dp^{r-j}(j)})$, where $\Gamma^d$, $S^d$, $\Lambda^d$ denote the $d$-fold divided power, symmetric power, and exterior power functor respectively. These computations are fundamental, for $\Gamma^d$ (respectively its dual $S^d$) is projective (resp. injective) in $\mathcal{P}(k)$ and tensor products of $\Gamma^i$ (resp. $S^i$) of total degree $d$ constitute a family of projective generators (resp. injective cogenerators) for $\mathcal{P}_d(k)$, the full subcategory of $\mathcal{P}(k)$ consisting of functors homogeneous of degree $d$. For example, Theorem 4.5 leads to the strong stability result (with respect to Frobenius twist) of Corollary 4.10 which is applicable to arbitrary strict polynomial functors of finite degree. The proof of Theorem 4.5 is an intricate nested triple induction argument. Readers of [F-L-S] or [F-S] will recognize here a new ingredient: the computation of the differentials in hypercohomology spectral sequences as Koszul differentials.

In Section 5, we extend the computation of Theorem 4.5 to other fundamental pairs of strict polynomial functors. Finally, in Section 6, we combine these computations of $\mathrm{Ext}_{\mathcal{P}(k)}$-groups with the strong comparison theorem of Section 3 and our understanding of base change for $\mathrm{Ext}_{\mathcal{F}(\mathbb{F}_q)}$-groups. The result is various complete calculations of $\mathrm{Ext}^*_{\mathcal{F}(\mathbb{F}_q)}(-,-)$, as tri-graded Hopf algebras.

The appendix, written by the last-named author, demonstrates a natural isomorphism $\mathrm{Ext}^*_{\mathcal{F}}(A,B) @>\sim>> \mathrm{Ext}^*_{GL(\mathbb{F}_q)}(A(\mathbb{F}_q^\infty), B(\mathbb{F}_q^\infty))$ for finite functors $A$, $B$ in $\mathcal{F}$.

*E. Friedlander gratefully acknowledges the hospitality of the University of Heidelberg.*

## 1. Recollections of functor categories

The purpose of this expository section is to recall the definitions and basic properties of the category $\mathcal{P}$ of strict polynomial functors of finite degree (on $k$-vector spaces) introduced in [F-S] and to be used in subsequent sections. We also contrast the category $\mathcal{P}$ with the category $\mathcal{F}$ of all functors from finite dimensional $k$-vector spaces to $k$-vector spaces. Strict polynomial functors were introduced in order to study rational cohomology of the general linear groups $\mathrm{GL}_n$ over $k$ (i.e., cohomology of comodules for the Hopf algebra $k[\mathrm{GL}_n]$, the coordinate algebra of the algebraic group $\mathrm{GL}_n$). Although one can consider strict polynomial functors over arbitrary commutative rings (as was necessarily done in [S-F-B]), we shall restrict attention throughout this paper to such functors defined on vector spaces over a field $k$ of characteristic $p > 0$. Indeed, in subsequent sections we specialize further to the case in which $k$ is a finite field.

We let $\mathcal{V} = \mathcal{V}_k$ denote the category of $k$-vector spaces and $k$-linear homomorphisms and we denote by $\mathcal{V}^f$ the full subcategory of finite dimensional $k$-vector spaces. A polynomial map $T : V \to W$ between finite dimensional vector spaces is defined to be a morphism of the corresponding affine schemes over $k$: $\mathrm{Spec}(S^*(V^\#)) \to \mathrm{Spec}(S^*(W^\#))$ (where $\mathrm{Spec}(S^*(V^\#))$ is the affine scheme associated to the symmetric algebra over $k$ of the $k$-linear dual of $V$). Equivalently, such a polynomial map is an element of $S^*(V^\#) \otimes W$. The polynomial map $T : V \to W$ is said to be homogeneous of degree $d$ if $T \in S^d(V^\#) \otimes W$. A polynomial map between finite dimensional vector spaces is uniquely determined by its associated set-theoretic function from $V$ to $W$ provided that $k$ is infinite; this is more readily understood as the observation that a polynomial (in any number of variables) is uniquely determined by its values at $k$-rational points provided that the base field is infinite.

## Strict polynomial functors

We recall [F-S, 2.1] that a *strict polynomial functor* $P : \mathcal{V}^f \to \mathcal{V}^f$ is the following collection of data: for any $V \in \mathcal{V}^f$ a vector space $P(V) \in \mathcal{V}^f$; for any $V, W$ in $\mathcal{V}^f$, a polynomial map $P_{V,W} : \mathrm{Hom}_k(V, W) \to \mathrm{Hom}_k(P(V), P(W))$. These polynomial maps should satisfy appropriate compatibility conditions similar to the ones used in the usual definition of a functor. A strict polynomial functor is said to be homogeneous of degree $d$ if $P_{V,W} : \mathrm{Hom}_k(V, W) \to \mathrm{Hom}_k(P(V), P(W))$ is homogeneous of degree $d$ for each pair $V, W$ in $\mathcal{V}^f$. A strict polynomial functor is said to be of finite degree provided that the degrees of the polynomial maps $P_{V,W}$ are bounded independent of $V, W \in \mathcal{V}^f$.

Replacing the polynomial maps $P_{V,W} : \mathrm{Hom}_k(V, W) \to \mathrm{Hom}_k(P(V), P(W))$ by the associated set-theoretic functions, we associate to each strict polynomial functor $P$ a functor in the usual sense $\mathcal{V}^f \to \mathcal{V}^f$ (which we usually denote by the same letter $P$). The above remarks about polynomial maps imply that over infinite fields strict polynomial functors may be viewed as functors (in the usual sense) which satisfy an appropriate additional property. However, over finite fields (and *a fortiori* over more general base rings), this is no longer the case: the concept of a strict polynomial functor incorporates more data than that of a functor in the usual sense.

We denote by $\mathcal{P}$ (or by $\mathcal{P}(k)$ if there is need to specify the base field $k$) the abelian category of strict polynomial functors of finite degree. One easily verifies that $\mathcal{P}$ splits as a direct sum $\oplus_d \mathcal{P}_d$, where $\mathcal{P}_d$ denotes the full subcategory of strict polynomial functors homogeneous of degree $d$.

Typical examples of strict polynomial functors homogeneous of degree $d$ are the $d$-fold tensor power functor $\otimes^d$, the $d$-fold exterior power functor $\Lambda^d$, the $d$-fold symmetric power functor $S^d$ (defined as the $\Sigma_d$-coinvariants of $\otimes^d$), and the $d$-fold divided power functor $\Gamma^d$ (defined as the $\Sigma_d$-invariants of $\otimes^d$). The functor $S^d$ is injective, and $\Gamma^d$ is projective in $\mathcal{P}_d$.

Let $\phi : k \to k$ denote the $p^{\text{th}}$-power map. The *Frobenius twist* functor $I^{(1)}$ sends a vector space $V$ to the base change of $V$ via the map $\phi$ (which we denote by $V^{(1)}$). So defined, $I^{(1)}$ is a strict polynomial functor homogeneous of degree $p$; $I^{(1)}_{V,W}$ is the $p^{\text{th}}$-power map

$$\operatorname{Hom}_k(V^{(1)}, W^{(1)})^\# = (\operatorname{Hom}_k(V, W)^\#)^{(1)} \to S^p(\operatorname{Hom}_k(V, W)^\#)$$

viewed as an element of $S^p(\operatorname{Hom}_k(V, W)^\#) \otimes \operatorname{Hom}_k(V^{(1)}, W^{(1)})$. For any functor $G : \mathcal{V}^f \to \mathcal{V}$, we define $G^{(1)}$ as $G \circ I^{(1)}$. (The reader should consult [F-S] for details.)

As observed in [S-F-B, §2], there is a natural construction of base change of a strict polynomial functor. Namely, if $k \to K$ is a field extension and $V$ is a $k$-vector space, let $V_K$ denote $K \otimes_k V$ with the evident $K$ vector space structure. If $P$ is a strict polynomial functor over $k$, then the base change $P_K$ is defined by setting $P_K(V_K) = P(V)_K$ and setting the polynomial map $(P_K)_{V_K,W_K} : \operatorname{Hom}_K(V_K, W_K) \to \operatorname{Hom}_K(P(V_K), P(W_K))$ to be the base change of $P_{V,W}$ as a morphism of affine schemes. As observed in [S-F-B, 2.6], this base change is both exact and preserves projectives. Consequently, we have the following elementary base change property.

PROPOSITION 1.1 ([S-F-B, 2.7]). *Let $P, Q$ be strict polynomial functors of finite degree over a field $k$. For any field extension $k \to K$, there is a natural isomorphism of $K$-vector spaces*

$$\operatorname{Ext}^*_{\mathcal{P}(K)}(P_K, Q_K) \cong \operatorname{Ext}^*_{\mathcal{P}(k)}(P, Q) \otimes_k K .$$

If $P$ is a strict polynomial functor, then for any $n > 0$ the vector space $P(k^n)$ inherits a natural structure of a rational $\operatorname{GL}_n$-module. (In fact, functoriality of $P$ implies that $P(k^n)$ inherits a natural structure of a rational $M_n$-module whose restriction provides the rational $\operatorname{GL}_n$-structure). We recall the following relationship between Ext-groups in the category $\mathcal{P}$ and in the category of rational $\operatorname{GL}_n$-modules.

THEOREM 1.2 ([F-S, 3.13]). *Let $P, Q$ be strict polynomial functors homogeneous of degree $d$ and let $n \geq d$. Then there is a natural isomorphism*

$$\operatorname{Ext}^*_{\mathcal{P}}(P, Q) \xrightarrow{\simeq} \operatorname{Ext}^*_{\operatorname{GL}_n}(P(k^n), Q(k^n))$$

*induced by the exact functor sending a strict polynomial functor to its value on $k^n$.*

The Ext-groups of the previous theorem can also be computed as the Ext-groups of the the classical Schur algebra $S(n,d)$ (with $n \geq d$ as above). An explicit determination of the cohomological dimension of $S(n,d)$ is given in [T].

By a theorem of H. Andersen (cf. [J, II, 10.14]), the Frobenius twist induces an injection on rational Ext-groups: for any two finite dimensional rational $\operatorname{GL}_n$-modules $M$ and $N$, the natural map induced by the Frobenius twist (which we view as an exact functor on the category of rational $\operatorname{GL}_n$-modules)

$$\operatorname{Ext}^*_{\operatorname{GL}_n}(M, N) \to \operatorname{Ext}^*_{\operatorname{GL}_n}(M^{(1)}, N^{(1)})$$

is injective. Thus, Theorem 1.2 gives us the following useful corollary.

COROLLARY 1.3. *Let $P, Q$ be strict polynomial functors homogeneous of degree $d$. The Frobenius twist is an exact functor on $\mathcal{P}$ which induces an injective map on* Ext-*groups*:

$$\operatorname{Ext}^*_{\mathcal{P}}(P, Q) \to \operatorname{Ext}^*_{\mathcal{P}}(P^{(1)}, Q^{(1)}).$$

Consider now the abelian category $\mathcal{F}$ of functors $\mathcal{V}^f \to \mathcal{V}$. If we need to indicate the base field $k$ explicitly, we shall denote this category by $\mathcal{F}(k)$. The forgetful functor $\mathcal{P} \to \mathcal{F}$ is clearly exact, thereby inducing a natural map on Ext-groups

$$\operatorname{Ext}^*_{\mathcal{P}}(P, Q) \to \operatorname{Ext}^*_{\mathcal{F}}(P, Q)$$

where we have abused notation by using $P$, $Q$ to denote strict polynomial functors and their images in $\mathcal{F}$.

The following elementary proposition provides a key to understanding the forgetful functor $\mathcal{P} \to \mathcal{F}$.

PROPOSITION 1.4. *Assume that $k$ has at least $d$ elements. Then for any $P$, $Q$ in $\mathcal{P}_d$, the natural inclusion $\operatorname{Hom}_{\mathcal{P}}(P, Q) \hookrightarrow \operatorname{Hom}_{\mathcal{F}}(P, Q)$ is an isomorphism.*

*Proof.* Let $f \in \operatorname{Hom}_{\mathcal{F}}(P, Q)$ be a homomorphism of functors. To check that $f$ is a homomorphism of strict polynomial functors, we have to verify that for any $V, W$ in $\mathcal{V}^f$ the following diagram of polynomial maps commutes:

$$\begin{CD}
\mathrm{Hom}_k(V,W) @>Q_{V,W}>> \mathrm{Hom}_k(Q(V), Q(W)) \\
@VP_{V,W}VV @VV\circ f_V V \\
\mathrm{Hom}_k(P(V), P(W)) @>f_W \circ>> \mathrm{Hom}_k(P(V), Q(W))
\end{CD}$$

To do so, we observe that both compositions are homogeneous polynomial maps of degree $d$ from $\mathrm{Hom}_k(V, W)$ to $\mathrm{Hom}_k(P(V), Q(W))$ whose values at all rational points coincide. Finally we observe that if a field $k$ contains at least $d$ elements then a homogeneous polynomial of degree $d$ (in any number of variables) which takes zero values at all rational points is necessarily zero (as a polynomial). $\square$

## Polynomial functors

We recall polynomial functors in the category $\mathcal{F}$. For a functor $F$ in $\mathcal{F}$, define its difference functor $\Delta(F)$ by

$$\Delta(F)(V) = \mathrm{Ker}\{F(V \oplus k) \to F(V)\}.$$

The functor $F$ is said to be polynomial if the $r^{\mathrm{th}}$ difference functor $\Delta^r(F)$ vanishes for $r$ sufficiently large. The Eilenberg-MacLane degree of a polynomial functor $F$ is the least integer $d$ such that $\Delta^{d+1}(F) = 0$. The same functors (or rather their images under the forgetful functor $\mathcal{P} \to \mathcal{F}$) used as examples of strict polynomial functors homogeneous of degree $d$ provide examples of polynomial functors of Eilenberg-MacLane degree $d$: the $d$-fold tensor power functor $\otimes^d$, the $d$-fold divided power functor $\Gamma^d$, etc. More generally one checks immediately that for any strict polynomial functor $P$ in $\mathcal{P}_d$ its image under the forgetful functor $\mathcal{P} \to \mathcal{F}$ is a polynomial of Eilenberg-MacLane degree less than or equal to $d$.

Following N. Kuhn [K1] we say that a functor $F$ in $\mathcal{F}$ is finite if it is of finite Eilenberg-MacLane degree and takes values in $\mathcal{V}^f$. The previous remarks imply that the image in $\mathcal{F}$ of any strict polynomial functor $P \in \mathcal{P}$ is finite.

For any vector space $V \in \mathcal{V}^f$ define a functor $P_V \in \mathcal{F}$ by the formula $P_V(W) = k[\mathrm{Hom}_k(V, W)]$. The Yoneda Lemma shows that for any $Q$ in $\mathcal{F}$ we have a natural isomorphism $\mathrm{Hom}_{\mathcal{F}}(P_V, Q) = Q(V)$. This implies immediately that the functor $P_V$ is projective in $\mathcal{F}$. A functor $P \in \mathcal{F}$ is said to be of finite type if it admits an epimorphism from a finite direct sum of functors of the form $P_V$. To say that a functor $Q \in \mathcal{F}$ admits a projective resolution of finite type is clearly equivalent to saying that $Q$ admits a resolution each term of which is isomorphic to a finite direct sum of functors of the form $P_V$. On several occasions we shall need the following useful fact.

PROPOSITION 1.5 ([S], [F-L-S, 10.1]). *Assume that the field $k$ is finite. Then every finite functor $Q$ in $\mathcal{F}$ admits a projective resolution of finite type.*

A clear difference between the two categories of functors $\mathcal{F}$ and $\mathcal{P}$ appears when we consider the Frobenius twist. If $k$ is perfect, then the Frobenius map $\phi : k \to k$ is an isomorphism so that the Frobenius twist $(-)^{(1)}$ becomes invertible when viewed in $\mathcal{F}$. Indeed, if $k$ is the prime field $\mathbb{F}_p$, then $I^{(1)} = I$ in $\mathcal{F}$; note that even for $k = \mathbb{F}_p$, $(-)^{(1)}$ is not invertible in $\mathcal{P}$.

Observe that if $k$ is perfect, then for any $P$, $Q$ in $\mathcal{F}$ we have a natural isomorphism $\mathrm{Ext}^*_{\mathcal{F}}(P,Q) @>\sim>> \mathrm{Ext}^*_{\mathcal{F}}(P^{(1)}, Q^{(1)})$. Hence for any strict polynomial functors $P$, $Q$ in $\mathcal{P}$ we get a natural map

$$\varinjlim_r \mathrm{Ext}^*_{\mathcal{P}}(P^{(r)}, Q^{(r)}) \to \varinjlim_r \mathrm{Ext}^*_{\mathcal{F}}(P^{(r)}, Q^{(r)}) \cong \mathrm{Ext}^*_{\mathcal{F}}(P,Q) \ .$$

Theorem 3.10 gives conditions for this map to be an isomorphism.

## Exponential functors

An *exponential functor* is a graded functor $A^* = (A^0, A^1, \ldots, A^n, \ldots)$ from $\mathcal{V}^f$ to $\mathcal{V}^f$ together with natural isomorphisms

$$A^0(V) \cong k \ , \qquad A^n(V \oplus W) \cong \bigoplus_{m=0}^n A^m(V) \otimes A^{n-m}(W) \ , \quad n > 0 \ .$$

LEMMA 1.6. *Let $A^*$ be an exponential functor.*

(1) *The functors $A^1$, $A^2$, ... are without constant term, i.e. $A^i(0) = 0$ for $i > 0$.*

(2) *The natural maps*

$$A^n(V) = A^n(V) \otimes A^0(W) \hookrightarrow \bigoplus_{m=0}^n A^m(V) \otimes A^{n-m}(W) = A^n(V \oplus W) \ ,$$

$$A^n(V \oplus W) = \bigoplus_{m=0}^n A^m(V) \otimes A^{n-m}(W) \twoheadrightarrow A^n(V) \otimes A^0(W) = A^n(V)$$

*coincide (up to an automorphism of the functor $A^n$) with the map induced by the inclusion of the first factor $A^n(i_1)$ and the map induced by the projection onto the first factor $A^n(p_1)$ respectively.*

(3) *The Eilenberg-MacLane degree of the functor $A^n$ is at most $n$ and is equal to $n$ provided that $A^1 \neq 0$. In particular the functors $A^n$ are finite.*

*Proof.* (1) The exponential condition shows that for $n > 0$ we have an isomorphism

$$A^n(0) = A^n(0 \oplus 0) \cong A^0(0) \otimes A^n(0) \oplus \cdots \oplus A^n(0) \otimes A^0(0) = A^n(0) \oplus \cdots \oplus A^n(0) \ .$$

Thus $\dim A^n(0) \geq 2 \dim A^n(0)$ and hence $\dim A^n(0) = 0$.

(2) Denote the homomorphisms in question by $i^n_{V,W}$ and $p^n_{V,W}$ respectively. The functoriality of all the maps involved implies the commutativity of the

following diagram:

$$A^n(V) @>i_{V,0}^n>> A^n(V \oplus 0) = A^n(V)$$
$$@V=VV @VV A^n(i_1) V$$
$$A^n(V) @>i_{V,W}^n>> A^n(V \oplus W) .$$

Now it suffices to note that (according to (1)) $i_{-,0}^n$ is an automorphism of the functor $A^n$. The same reasoning applies to $p_{V,W}^n$.

(3) The exponential condition and (2) show that the functor $\Delta(A^n)$ is isomorphic to the direct sum $A^1(k) \otimes A^{n-1} \oplus \cdots \oplus A^n(k)$. Immediate induction on $n$ now concludes the proof. □

Typical examples of exponential functors are given by the symmetric algebra $S^* = (S^0, S^1, \ldots S^n, \ldots)$, the exterior algebra $\Lambda^* = (\Lambda^0, \Lambda^1, \ldots, \Lambda^n, \ldots)$ and the divided power algebra $\Gamma^* = (\Gamma^0, \Gamma^1, \ldots, \Gamma^n, \ldots)$. Another exponential functor $L^* = (L^0, L^1, \ldots, L^n, \ldots)$ is obtained as the quotient of the symmetric power algebra by the ideal of $p^{\text{th}}$ powers; it coincides with the exterior algebra when $p = 2$.

Note also that if we define the tensor product of graded functors via the usual formula

$$(A^* \otimes B^*)^n = \bigoplus_{m=0}^{n} A^m \otimes B^{n-m},$$

then the tensor product of two exponential functors is again exponential. Clearly, the Frobenius twist of an exponential functor is again exponential.

The previous definition generalizes immediately to the case of strict polynomial functors. We skip the obvious details.

The following theorem, in the case of the category $\mathcal{F}$, was used in [F]. It generalizes a result due to Pirashvili [P] (much used in [F-S] and [F-L-S]) which asserts that if $A$ is an additive functor and $B$, $C$ are functors without constant term, then all Ext-groups from $A$ to $B \otimes C$ are 0. The isomorphism of Theorem 1.7 provides an important tool enabling computations of functor cohomology. From now on we assume (if not specified otherwise) that the base field $k$ is finite.

THEOREM 1.7. *Let $A^*$ be an exponential functor. For any $B$, $C$ in $\mathcal{F}$, there exist natural isomorphisms*

$$\operatorname{Ext}^*_{\mathcal{F}}(A^n, B \otimes C) = \bigoplus_{m=0}^{n} \operatorname{Ext}^*_{\mathcal{F}}(A^m, B) \otimes \operatorname{Ext}^*_{\mathcal{F}}(A^{n-m}, C).$$

*Furthermore, if B and C take values in the category $\mathcal{V}^f$ then there also exist natural isomorphisms*

$$\operatorname{Ext}^*_{\mathcal{F}}(B \otimes C, A^n) = \bigoplus_{m=0}^{n} \operatorname{Ext}^*_{\mathcal{F}}(B, A^m) \otimes \operatorname{Ext}^*_{\mathcal{F}}(C, A^{n-m}).$$

*Similarly, if $A^*$ is an exponential strict polynomial functor (of finite degree), then for any $B$, $C$ in $\mathcal{P}$ there are natural isomorphisms*

$$\operatorname{Ext}^*_{\mathcal{P}}(A^n, B \otimes C) = \bigoplus_{m=0}^{n} \operatorname{Ext}^*_{\mathcal{P}}(A^m, B) \otimes \operatorname{Ext}^*_{\mathcal{P}}(A^{n-m}, C)$$

$$\operatorname{Ext}^*_{\mathcal{P}}(B \otimes C, A^n) = \bigoplus_{m=0}^{n} \operatorname{Ext}^*_{\mathcal{P}}(B, A^m) \otimes \operatorname{Ext}^*_{\mathcal{P}}(C, A^{n-m}).$$

*Proof.* In case of the category $\mathcal{F}$ the first part is proved in [F, 1.4.2] (we recall the proof below). The same argument gives the second part. Alternatively, the second part follows from the first by the duality isomorphism 1.12. The proof for the category $\mathcal{P}$ is identical, using the theory of strict polynomial bifunctors of finite degree as developed in [S-F-B]. It should be noted also that in case of the category $\mathcal{P}$ the above theorem holds over arbitrary (not necessarily finite) fields.

Let bi$-\mathcal{F}$ denote the abelian category of bifunctors $\mathcal{V}^f \times \mathcal{V}^f \to \mathcal{V}$. Consider a pair of adjoint (on both sides) functors

$$\mathcal{V}^f \underset{@<\Pi<<}{\overset{@>D>>}{\rightleftarrows}} \mathcal{V}^f \times \mathcal{V}^f.$$

Here $\Pi$ is the direct sum functor $\Pi(V,W) = V \oplus W$ and $D$ is the diagonal functor $D(V) = (V,V)$. Taking compositions on the right with these functors we get a pair of adjoint (on both sides) functors between $\mathcal{F}$ and bi$-\mathcal{F}$.

$$\mathcal{F} \underset{@<Q \circ D \leftarrow Q<<}{\overset{@>P \mapsto P \circ \Pi>>}{\rightleftarrows}} \text{bi}-\mathcal{F}.$$

Since both functors are exact they also preserve projectives and injectives and we get the usual adjunction isomorphisms.

(1.7.1) *For any $P$ in $\mathcal{F}$ and any $Q$ in* bi$-\mathcal{F}$ *we have natural isomorphisms*

$$\operatorname{Ext}^*_{\mathcal{F}}(P, Q \circ D) = \operatorname{Ext}^*_{\text{bi}-\mathcal{F}}(P \circ \Pi, Q) ,$$
$$\operatorname{Ext}^*_{\mathcal{F}}(Q \circ D, P) = \operatorname{Ext}^*_{\text{bi}-\mathcal{F}}(Q, P \circ \Pi) .$$

For any functors $B$, $C$ in $\mathcal{F}$ we define their external tensor product $B \boxtimes C \in$ bi$-\mathcal{F}$ via the formula $B \boxtimes C(V,W) = B(V) \otimes C(W)$. Exactness of tensor products and an obvious formula for the external tensor product of projective generators: $P_V \boxtimes P_W = P_{(V,W)}$, give us the following Künneth-type formula:

(1.7.2) *Assume that $A_1$, $A_2$ are functors in $\mathcal{F}$ which admit projective resolutions of finite type; then for any $B$, $C$ in $\mathcal{F}$ there are natural isomorphisms*

$$\mathrm{Ext}^*_{\mathrm{bi}-\mathcal{F}}(A_1 \boxtimes A_2, B \boxtimes C) = \mathrm{Ext}^*_{\mathcal{F}}(A_1, B) \otimes \mathrm{Ext}^*_{\mathcal{F}}(A_2, C).$$

Using (1.7.1) and the exponential property of $A^*$ we get isomorphisms:

$$\mathrm{Ext}^*_{\mathcal{F}}(A^n, B \otimes C) = \mathrm{Ext}^*_{\mathcal{F}}(A^n, (B \boxtimes C) \circ D) = \mathrm{Ext}^*_{\mathrm{bi}-\mathcal{F}}(A^n \circ \Pi, B \boxtimes C)$$
$$= \mathrm{Ext}^*_{\mathrm{bi}-\mathcal{F}}(\bigoplus_{m=0}^{n} A^m \boxtimes A^{n-m}, B \boxtimes C) = \bigoplus_{m=0}^{n} \mathrm{Ext}^*_{\mathrm{bi}-\mathcal{F}}(A^m \boxtimes A^{n-m}, B \boxtimes C).$$

Finally we observe that all functors $A^i$ are finite according to Lemma 1.6 and hence admit projective resolutions of finite type (see Proposition 1.5), so that we may use (1.7.2) to conclude the proof. □

Observe that for homogeneous strict polynomial functors $A^n$, $B$ and $C$, the sum in Theorem 1.7 cannot have more than one nonzero term. Thus, in the case of the category $\mathcal{P}$, Theorem 1.7 generalizes [F-S, Prop. 5.2]. This, and the injectivity of the symmetric powers in $\mathcal{P}$, make computing Ext-groups in the category $\mathcal{P}$ easier than in the category $\mathcal{F}$.

Applying Theorem 1.7 to the exponential functors $\underbrace{A^* \otimes \cdots \otimes A^*}_{n}$ and $\underbrace{B^* \otimes \cdots \otimes B^*}_{m}$ we obtain the following corollary.

COROLLARY 1.8. *Let $A^*$, $B^*$ be exponential functors. For any nonnegative integers $k_1, \ldots, k_n$; $l_1, \ldots, l_m$ there are natural isomorphisms*

$$\mathrm{Ext}^*_{\mathcal{F}}(A^{k_1} \otimes \cdots \otimes A^{k_n}, B^{l_1} \otimes \cdots \otimes B^{l_m})$$
$$= \bigoplus_{\substack{k_{1,1}+\cdots+k_{1,m}=k_1 \\ l_{1,1}+\cdots+l_{1,n}=l_1}} \cdots \bigoplus_{\substack{k_{n,1}+\cdots+k_{n,m}=k_n \\ l_{m,1}+\cdots+l_{m,n}=l_m}} \bigotimes_{s=1,\ t=1}^{n\ \ m} \mathrm{Ext}^*_{\mathcal{F}}(A^{k_{s,t}}, B^{l_{t,s}}).$$

*Similarly, if $A^*$ and $B^*$ are exponential strict polynomial functors of finite degree, then there are corresponding natural isomorphisms obtained by replacing $\mathrm{Ext}^*_{\mathcal{F}}$ by $\mathrm{Ext}^*_{\mathcal{P}}$.*

All our examples of exponential functors are bi-algebras. Indeed, for any exponential functor $A^*$, the natural maps

$$A^i(V) \otimes A^j(V) \hookrightarrow A^{i+j}(V \oplus V) @>A^{i+j}(\Sigma)>> A^{i+j}(V),$$
$$A^{i+j}(V) @>A^{i+j}(\Delta)>> A^{i+j}(V \oplus V) \twoheadrightarrow A^i(V) \otimes A^j(V)$$

define natural product and coproduct operations $A^i \otimes A^j \to A^{i+j}$, $A^{i+j} \to A^i \otimes A^j$. Here, $\Sigma$ is the sum map $(x,y) \mapsto x+y$ and $\Delta$ is the diagonal map $x \mapsto (x,x)$.

*Definition* 1.9. $A^*$ is a *Hopf exponential functor* provided that for any $V$ in $\mathcal{V}^f$ the above product operation makes $A^*(V)$ into a (graded) associative $k$-algebra with unit $1 \in A^0(V) = k$.

The name *Hopf functor* is justified by the following lemma.

LEMMA 1.10. *Assume that $A^*$ is an Hopf exponential functor. Then for any $V$ in $\mathcal{V}^f$ the above operations make $A^*(V)$ into a (graded) Hopf algebra with co-unit $\varepsilon : A^*(V) \to A^0(V) = k$.*

*Proof.* Note first that the (right) multiplication by $1 \in A^0(V)$ coincides with the natural isomorphism

$$i_{V,0}^n : A^n(V) = A^n(V) \otimes A^0(0) @>\sim>> A^n(V \oplus 0) = A^n(V)$$

(cf. the proof of Lemma 1.6 (2)). In the same way the right comultiplication by $\varepsilon$ coincides with the natural isomorphism

$$p_{V,0}^n : A^n(V) = A^n(V \oplus 0) @>\sim>> A^n(V) \otimes A^0(0) = A^n(V).$$

Since $p_{V,0}^n = (i_{V,0}^n)^{-1}$ we conclude that $1 \in A^0(V)$ is a unit if and only if $\varepsilon$ is a co-unit of $A^*(V)$.

Observe next that the natural homomorphism $A^m(V) \otimes A^{n-m}(W) \to A^n(V \oplus W)$ in the definition of the exponential functor may be expressed in terms of the product operation as the composition

$$A^m(V) \otimes A^{n-m}(W) @> A^m(i_1) \otimes A^{n-m}(i_2) >> A^m(V \oplus W) \otimes A^{n-m}(V \oplus W)$$
$$@> \text{mult} >> A^n(V \oplus W).$$

This remark implies immediately that the associativity of $A^*(V)$ (for all $V$) is equivalent to the fact that the two possible identifications of graded tri-functors $A^*(U) \otimes A^*(V) \otimes A^*(W)$ and $A^*(U \oplus V \oplus W)$ coincide. Now the verification of the fact that $A^*$ is a Hopf algebra becomes a straightforward computation. For example to check the coassociativity we have to verify the commutativity of the following diagram

$$A^*(V) @> \text{comult} >> A^*(V) \otimes A^*(V)$$
$$@V \text{comult} VV @V \text{comult} \otimes 1 VV$$
$$A^*(V) \otimes A^*(V) @> 1 \otimes \text{comult} >> A^*(V) \otimes A^*(V) \otimes A^*(V).$$

However one checks easily that both compositions coincide with the homomorphism $A^*(V) @> A^*(\Delta) >> A^*(V \oplus V \oplus V) = A^*(V) \otimes A^*(V) \otimes A^*(V)$, where $\Delta : V \to V \oplus V \oplus V$ is the diagonal map. □

We say that the exponential functor $A^*$ is commutative (respectively, skew-commutative) if for all nonnegative integers $i$, $j$ and every $V$ in $\mathcal{V}^f$ the map $\tau : V \oplus V \to V \oplus V$, $(x,y) \mapsto (y,x)$ gives rise to a commutative diagram (resp. a diagram commutative up to a sign $(-1)^{ij}$):

$$\begin{CD}
A^i(V) \otimes A^j(V) @>>> A^{i+j}(V \oplus V) \\
@VTVV @VVA(\tau)V \\
A^j(V) \otimes A^i(V) @>>> A^{i+j}(V \oplus V)
\end{CD}$$

where $T$ is the twist map. For a commutative exponential functor $A^*$ we set $\varepsilon(A^*) = 1$ and for a skew commutative one we set $\varepsilon(A^*) = -1$. We say in these cases that $A^*$ is $\varepsilon(A^*)$-commutative. When an exponential functor is commutative (resp. skew-commutative), the product in the algebra $A^*(V)$ is (skew-) commutative and the coproduct is (skew-) cocommutative. We readily verify that the functors $\Gamma^*$, $\Lambda^*$, $S^*$, $L^*$ and their Frobenius twists are Hopf exponential functors and that $\varepsilon(\Gamma^*) = \varepsilon(S^*) = \varepsilon(L^*) = +1$, $\varepsilon(\Lambda^*) = -1$.

Assume that $A^*$ and $B^*$ are exponential functors (respectively, exponential strict polynomial functors of finite degree). In this case the tri-graded vector space $\operatorname{Ext}^*(A^*, B^*)$ (a notation to denote $\operatorname{Ext}^*_{\mathcal{F}}(A^*, B^*)$ as well as $\operatorname{Ext}^*_{\mathcal{P}}(A^*, B^*)$) acquires a natural structure of a (tri-graded) bi-algebra. The product operation is defined as the composition

$$\begin{aligned}
\operatorname{Ext}^*(A^*, B^*) \otimes \operatorname{Ext}^*(A^*, B^*) &\xrightarrow{\sim} \operatorname{Ext}^*(A^* \boxtimes A^*, B^* \boxtimes B^*) \\
&= \operatorname{Ext}^*(A^* \circ \Pi, B^* \circ \Pi) \to \operatorname{Ext}^*(A^* \circ \Pi \circ D, B^* \circ \Pi \circ D) \to \operatorname{Ext}^*(A^*, B^*).
\end{aligned}$$

Here the second arrow is induced by the exact functor $\text{bi-}\mathcal{F} \xrightarrow{Q \mapsto Q \circ D} \mathcal{F}$ and the last arrow is defined by the adjunction homomorphisms $\Delta : I \to \Pi \circ D$ and $\Sigma : \Pi \circ D \to I$. Equivalently the product of the classes $e \in \operatorname{Ext}^i(A^n, B^m)$, $e' \in \operatorname{Ext}^{i'}(A^{n'}, B^{m'})$ may be described as the image of the tensor product $e \otimes e' \in \operatorname{Ext}^{i+i'}(A^n \otimes A^{n'}, B^m \otimes B^{m'})$ under the homomorphism $\operatorname{Ext}^*(A^n \otimes A^{n'}, B^m \otimes B^{m'}) \to \operatorname{Ext}^*(A^{n+n'}, B^{m+m'})$ induced by the product operation $B^m \otimes B^{m'} \to B^{m+m'}$ and the coproduct operation $A^{n+n'} \to A^n \otimes A^{n'}$. Finally the product operation may be also described (using the isomorphisms of Theorem 1.7) as eigher of the following two compositions

$$\operatorname{Ext}^*(A^*, B^*) \otimes \operatorname{Ext}^*(A^*, B^*) \xrightarrow{\sim} \operatorname{Ext}^*(A^* \otimes A^*, B^*) \to \operatorname{Ext}^*(A^*, B^*)$$
$$\operatorname{Ext}^*(A^*, B^*) \otimes \operatorname{Ext}^*(A^*, B^*) \xrightarrow{\sim} \operatorname{Ext}^*(A^*, B^* \otimes B^*) \to \operatorname{Ext}^*(A^*, B^*).$$

In a similar fashion the coproduct is defined as the composition

$$\begin{aligned}
\operatorname{Ext}^*(A^*, B^*) &\to \operatorname{Ext}^*(A^* \circ \Pi, B^* \circ \Pi) = \operatorname{Ext}^*(A^* \boxtimes A^*, B^* \boxtimes B^*) \\
&= \operatorname{Ext}^*(A^*, B^*) \otimes \operatorname{Ext}^*(A^*, B^*).
\end{aligned}$$

The coproduct operation may be also described (using the isomorphisms of Theorem 1.7) as either one of the following compositions

$$\operatorname{Ext}^*(A^*, B^*) \to \operatorname{Ext}^*(A^* \otimes A^*, B^*) \overset{\sim}{\twoheadrightarrow} \operatorname{Ext}^*(A^*, B^*) \otimes \operatorname{Ext}^*(A^*, B^*)$$

$$\operatorname{Ext}^*(A^*, B^*) \to \operatorname{Ext}^*(A^*, B^* \otimes B^*) \overset{\sim}{\twoheadrightarrow} \operatorname{Ext}^*(A^*, B^*) \otimes \operatorname{Ext}^*(A^*, B^*).$$

The proof of the following lemma is straightforward (but tiresome).

LEMMA 1.11. *Let $A^*$ and $B^*$ be Hopf exponential functors (resp. Hopf exponential strict polynomial functors of finite degree). Then the tri-graded vector space* $\operatorname{Ext}^*(A^*, B^*)$ *has a natural structure of a (tri-graded) Hopf algebra. Moreover if $A^*$ is $\varepsilon(A^*)$-commutative and $B^*$ is $\varepsilon(B^*)$-commutative, then the following diagrams commute up to a sign* $(-1)^{st + \frac{\varepsilon(A)-1}{2} \cdot ik + \frac{\varepsilon(B)-1}{2} \cdot jl}$:

$$\begin{CD}
\operatorname{Ext}^s(A^i, B^j) \otimes \operatorname{Ext}^t(A^k, B^l) @>\operatorname{mult}>> \operatorname{Ext}^{s+t}(A^{i+k}, B^{j+l}) \\
@V\cong VV @V=VV \\
\operatorname{Ext}^t(A^k, B^l) \otimes \operatorname{Ext}^s(A^i, B^j) @>\operatorname{mult}>> \operatorname{Ext}^{s+t}(A^{i+k}, B^{j+l})
\end{CD}$$

$$\begin{CD}
\operatorname{Ext}^{s+t}(A^{i+k}, B^{j+l}) @>\operatorname{comult}>> \operatorname{Ext}^s(A^i, B^j) \otimes \operatorname{Ext}^t(A^k, B^l) \\
@V=VV @V\cong VV \\
\operatorname{Ext}^{s+t}(A^{i+k}, B^{j+l}) @>\operatorname{comult}>> \operatorname{Ext}^t(A^k, B^l) \otimes \operatorname{Ext}^s(A^i, B^j).
\end{CD}$$

## Dual functors

For a functor $P : \mathcal{V}^f \to \mathcal{V}$, define its dual $P^\# : \mathcal{V}^f \to \mathcal{V}$ by $P^\#(V) = P(V^\#)^\#$. Two examples of this duality which will be important in what follows are

$$(S^d)^\# = \Gamma^d, \quad (\Lambda^d)^\# = \Lambda^d.$$

The contravariant functor $\# : \mathcal{F} \to \mathcal{F}$ is clearly exact and hence for any $P, Q$ in $\mathcal{F}$ we get a natural duality homomorphism $\# : \operatorname{Ext}^*_{\mathcal{F}}(P, Q) \to \operatorname{Ext}^*_{\mathcal{F}}(Q^\#, P^\#)$. The same construction obviously applies to the category $\mathcal{P}$ as well.

LEMMA 1.12. *Assume that $P, Q \in \mathcal{F}$ take values in the category $\mathcal{V}^f$. Then the duality homomorphism*

$$\# : \operatorname{Ext}^*_{\mathcal{F}}(P, Q) \to \operatorname{Ext}^*_{\mathcal{F}}(Q^\#, P^\#)$$

*is an isomorphism. Similarly for any $P, Q \in \mathcal{P}$ the duality homomorphism*

$$\# : \operatorname{Ext}^*_{\mathcal{P}}(P, Q) \to \operatorname{Ext}^*_{\mathcal{P}}(Q^\#, P^\#)$$

*is an isomorphism.*

*Proof.* For any $F$ in $\mathcal{F}$ we have a natural homomorphism $F \to F^{\#\#}$, which is an isomorphism provided that $F$ takes values in $\mathcal{V}^f$. In particular, $P^{\#\#} = P$, $Q^{\#\#} = Q$. Now it suffices to show that the composition

$$\operatorname{Ext}^*_{\mathcal{F}}(P, Q) \xrightarrow{\#} \operatorname{Ext}^*_{\mathcal{F}}(Q^{\#}, P^{\#}) \xrightarrow{\#} \operatorname{Ext}^*_{\mathcal{F}}(P^{\#\#}, Q^{\#\#}) = \operatorname{Ext}^*_{\mathcal{F}}(P, Q)$$

is the identity map. Let $e \in \operatorname{Ext}^n_{\mathcal{F}}(P, Q)$ be represented by an extension

$$0 \to Q \to P_n \to \cdots \to P_1 \to P \to 0.$$

Then $e^{\#\#}$ is represented by the extension

$$0 \to Q^{\#\#} \to P_n^{\#\#} \to \cdots \to P_1^{\#\#} \to P^{\#\#} \to 0$$

and our statement follows from the commutativity of the diagram

$$\begin{CD}
0 @>>> Q @>>> P_n @>>> \ldots @>>> P_1 @>>> P @>>> 0 \\
@. @VV{\cong}V @VVV @. @VVV @VV{\cong}V \\
0 @>>> Q^{\#\#} @>>> P_n^{\#\#} @>>> \ldots @>>> P_1^{\#\#} @>>> P^{\#\#} @>>> 0.
\end{CD}$$

The case of the category $\mathcal{P}$ is trivial since in this case $\# : \mathcal{P} \to \mathcal{P}$ is an anti-equivalence. $\square$

Clearly, the dual of an exponential functor is again exponential. The following result (to be used in §5) is straightforward from the definitions.

LEMMA 1.13. *Let $A^*, B^*$ be Hopf exponential functors (resp. Hopf exponential strict polynomial functors of finite degree). The natural isomorphism of Lemma* 1.12

$$\operatorname{Ext}^*(A^*, B^*) \xrightarrow{\sim} \operatorname{Ext}^*(B^{*\#}, A^{*\#})$$

*is an anti-isomorphism of graded Hopf algebras (i.e. $(xy)^{\#} = y^{\#} x^{\#}$ etc.).*

## 2. The weak comparison theorem

In this section, we investigate the map on Ext-groups induced by the forgetful functor $\mathcal{P} \to \mathcal{F}$. Throughout this section, and in much of the next, we shall restrict our attention to the case in which the base field $k$ is a finite field $k = \mathbb{F}_q$ of characteristic $p$, with $q = p^N$. We show for strict polynomial functors $A, B$ of degree $d$ that the natural map

$$\operatorname{Ext}^i_{\mathcal{P}}(A^{(m)}, B^{(m)}) \to \operatorname{Ext}^i_{\mathcal{F}}(A^{(m)}, B^{(m)}) = \operatorname{Ext}^i_{\mathcal{F}}(A, B)$$

is an isomorphism provided that $m$ and $q$ are sufficiently large compared to $d$ and $i$. In the following section, we show that the condition on $q$ can be weakened to the simple condition that $q \geq d$.

We begin with the following vanishing theorem which complements the fact that symmetric power functors $S^d$ are injective in $\mathcal{P}_d$. The proof of this theorem uses the finiteness of the base field $k$ in an essential way.

THEOREM 2.1 ([F, §4]). *Let $A : \mathcal{V}^f \to \mathcal{V}^f$ be a finite functor of Eilenberg-MacLane degree $\leq d$. For any $i > 0$,*

$$\mathrm{Ext}^i_{\mathcal{F}}(A, S^{p^h}) = 0$$

*provided that $h \geq \log_p \frac{i+2}{2} + [\frac{d-1}{p-1}]$.*

*Proof.* By [F, 4.1.2], the natural embedding $S^{p^h(l+1)} \hookrightarrow S^{p^{h+1}(l)}$ induces an isomorphism

$$\mathrm{Ext}^i_{\mathcal{F}}(A, S^{p^h(l+1)}) @>\sim>> \mathrm{Ext}^i_{\mathcal{F}}(A, S^{p^{h+1}(l)})$$

provided that $h \geq \log_p \frac{i+2}{2} + [\frac{d-1}{p-1}]$. Applying this remark $N$ times we see that (for the same values of $h$) the natural embedding

$$S^{p^h} = S^{p^h(N)} \hookrightarrow S^{p^{h+1}(N-1)} \hookrightarrow \cdots \hookrightarrow S^{p^{h+N}}$$

induces an isomorphism

$$\mathrm{Ext}^i_{\mathcal{F}}(A, S^{p^h}) @>\sim>> \mathrm{Ext}^i_{\mathcal{F}}(A, S^{p^{h+N}}).$$

On the other hand, the vanishing theorem of N. Kuhn [K3, Th. 1.1] shows that

$$\varinjlim_{k \geq 0} \mathrm{Ext}^i_{\mathcal{F}}(A, S^{p^{h+kN}}) = 0$$

for any functor $A$ admitting a projective resolution of finite type. Since each finite functor $A$ admits such a resolution (see Proposition 1.5), this completes the proof. $\square$

The following classical lemma was used in [F] and thus implicitly in the proof of Theorem 2.1. In Proposition 2.3, it enables us to extend the vanishing statement of Theorem 2.1 to tensor products of symmetric power functors.

LEMMA 2.2. *Let $n$ be an integer with $p$-adic expansion*

$$n = n_0 + n_1 p + \cdots + n_k p^k \quad (0 \leq n_i < p) .$$

*Then the functor $S^n$ is canonically a direct summand in $(S^1)^{\otimes n_0} \otimes \cdots \otimes (S^{p^k})^{\otimes n_k}$.*

*Proof.* For any integer $m$ with $0 \leq m \leq n$, the multiplication and comultiplication of the exponential functor $S^*$ yield natural homomorphisms

$$S^n @>\mathrm{comult}>> S^m \otimes S^{n-m} @>\mathrm{mult}>> S^n .$$

Explicitly, $S^n @>\text{comult}>> S^m \otimes S^{n-m}$ sends an $n$-fold product $x_1 \cdots x_n \in S^n(V)$ to the sum indexed by $\Sigma_n/\Sigma_m \times \Sigma_{n-m}$ of tensors of the form $x_{i_1} \cdots x_{i_m} \otimes x_{i_{m+1}} \cdots x_{i_n}$. Composition of these homomorphisms is multiplication by $\binom{n}{m}$, so $S^n$ is a direct summand in $S^m \otimes S^{n-m}$ provided that $\binom{n}{m}$ is prime to $p$. One easily verifies that $\binom{n}{m}$ is prime to $p$ if and only if all $p$-adic digits of $m$ are less or equal to the corresponding digits of $n$. The lemma now follows by an induction on the sum of $p$-adic digits of $n$. □

PROPOSITION 2.3. *Let $A^* \in \mathcal{F}$ be an exponential functor. Assume that*

$$h \geq \log_p \frac{s+2}{2} + [\frac{d-1}{p-1}], \quad \text{and} \quad n_1 + \cdots + n_k \equiv 0 \mod p^h.$$

*Then for all $0 < i \leq s$*

$$\mathrm{Ext}^i_{\mathcal{F}}(A^d, S^{n_1} \otimes \cdots \otimes S^{n_k}) = 0.$$

*Proof.* Using Lemma 2.2, we are easily reduced to the case where $n_i = p^{m_i}$ is a power of $p$ for all $i$. By Theorem 1.7,

$$\mathrm{Ext}^*_{\mathcal{F}}(A^d, S^{p^{m_1}} \otimes \cdots \otimes S^{p^{m_k}}) = \bigoplus_{d_1+\cdots+d_k=d} \bigotimes_{j=1}^{k} \mathrm{Ext}^*_{\mathcal{F}}(A^{d_j}, S^{p^{m_j}}).$$

Since the functor $A^0$ is constant, all its Ext-groups to functors without constant term are trivial; thus the summand corresponding to the $k$-tuple $(d_1, \ldots, d_k)$ is trivial provided that $d_j = 0$ for some $j$. Consider now a summand corresponding to a $k$-tuple $(d_1, \ldots, d_k)$ with all $d_j > 0$. In this case we have inequalities $d_j \leq d - (k-1)$. On the other hand, an easy induction on $k$ establishes that if $p^{m_1} + \cdots + p^{m_k} \equiv 0 \mod p^h$, then $m_j \geq h - [\frac{k-1}{p-1}]$ for all $j$. Consequently,

$$m_j \geq h - \left[\frac{k-1}{p-1}\right] \geq \log_p \frac{s+2}{2} + \left[\frac{d-1}{p-1}\right] - \left[\frac{k-1}{p-1}\right] \geq \log_p \frac{s+2}{2} + \left[\frac{d_j-1}{p-1}\right].$$

Theorem 2.1 shows that $\mathrm{Ext}^i_{\mathcal{F}}(A^{d_j}, S^{p^{m_j}}) = 0$ for $0 < i \leq s$, thus concluding the proof. □

COROLLARY 2.4. *Assume that $h \geq \log_p \frac{s+2}{2} + [\frac{d-1}{p-1}]$. Let $d_1, \ldots, d_k$; $n_1, \ldots, n_l$ be integers such that $d_1 + \cdots + d_k = d$ and $n_1 + \cdots + n_l \equiv 0 \mod p^h$. Then*

$$\mathrm{Ext}^i_{\mathcal{F}}((\Gamma^{d_1} \otimes \cdots \otimes \Gamma^{d_k})^{(m)}, S^{n_1} \otimes \cdots \otimes S^{n_l}) = 0$$

*for all $m$ and all $0 < i \leq s$.*

*Proof.* This follows immediately from Proposition 2.3 since $(\Gamma^{d_1} \otimes \cdots \otimes \Gamma^{d_k})^{(m)}$ is a direct summand in $A^d$, where $A^*$ is the exponential functor $(\underbrace{\Gamma^* \otimes \cdots \otimes \Gamma^*}_{k})^{(m)}$. □

Corollary 2.4 easily implies our first, and weakest, comparison of the Ext-groups $\mathrm{Ext}^*_{\mathcal{P}}$ and $\mathrm{Ext}^*_{\mathcal{F}}$. A first form of this proposition was proved by E. Friedlander (1995, unpublished) and by S. Betley independently [B2], and discussion of T. Pirashvili with the first author.

PROPOSITION 2.5. *Let $A$ and $B$ be strict polynomial functors homogeneous of degree $d$, over a field with $q$ elements. Assume that*

$$m \geq \log_p \frac{s+2}{2} + [\frac{d-1}{p-1}] \,, \quad q \geq dp^m.$$

*Then the canonical homomorphisms*

$$\mathrm{Ext}^i_{\mathcal{P}}(A^{(m)}, B^{(m)}) \to \mathrm{Ext}^i_{\mathcal{F}}(A^{(m)}, B^{(m)}) \cong \mathrm{Ext}^i_{\mathcal{F}}(A, B)$$

*are isomorphisms for all $0 \leq i \leq s$.*

*Proof.* Consider first the special case $A = \Gamma^{d_1} \otimes \cdots \otimes \Gamma^{d_k}$, $d_1 + \cdots + d_k = d$. The functor $B^{(m)}$ admits an injective resolution in $\mathcal{P}_{dp^m}$

$$0 \to B^{(m)} \to I^\bullet$$

in which every $I^j$ is a direct sum of functors of the form $S^{n_1} \otimes \cdots \otimes S^{n_l}$ with $n_1 + \cdots + n_l = dp^m$ [F-S, §2]. We employ the hypercohomology spectral sequence

$$E_1^{ij} = \mathrm{Ext}^i_{\mathcal{F}}(A^{(m)}, I^j) \Longrightarrow \mathrm{Ext}^{i+j}_{\mathcal{F}}(A^{(m)}, B^{(m)}).$$

Corollary 2.4 shows that $E_1^{ij} = 0$ for all $0 < i \leq s$. This shows that for $0 \leq i \leq s$ the group $\mathrm{Ext}^i_{\mathcal{F}}(A^{(m)}, B^{(m)})$ coincides with the $i^{\mathrm{th}}$ homology group of the complex $\mathrm{Hom}_{\mathcal{F}}(A^{(m)}, I^\bullet)$, which equals $\mathrm{Hom}_{\mathcal{P}}(A^{(m)}, I^\bullet)$ by Proposition 1.4. The homology of this latter complex equals $\mathrm{Ext}^*_{\mathcal{P}}(A^{(m)}, B^{(m)})$.

In the general case we consider a projective resolution

$$0 \leftarrow A \leftarrow P_\bullet$$

in which each $P_i$ is a direct sum of functors of the form $\Gamma^{d_1} \otimes \cdots \otimes \Gamma^{d_k}$, $d_1 + \cdots + d_k = d$. This gives us two hypercohomology spectral sequences

$$_{\mathcal{P}}E_1^{ij} = \mathrm{Ext}^i_{\mathcal{P}}(P_j^{(m)}, B^{(m)}) \Longrightarrow \mathrm{Ext}^{i+j}_{\mathcal{P}}(A^{(m)}, B^{(m)}) \,,$$
$$_{\mathcal{F}}E_1^{ij} = \mathrm{Ext}^i_{\mathcal{F}}(P_j^{(m)}, B^{(m)}) \Longrightarrow \mathrm{Ext}^{i+j}_{\mathcal{F}}(A^{(m)}, B^{(m)})$$

and a homomorphism of spectral sequences $_{\mathcal{P}}E \to {_{\mathcal{F}}E}$. According to the special case considered above, the homomorphisms $_{\mathcal{P}}E_1^{ij} \to {_{\mathcal{F}}E_1^{ij}}$ on $E_1^{ij}$-terms are isomorphisms for all $j$ and all $0 \leq i \leq s$. The standard comparison theorem for spectral sequences thus implies that $\mathrm{Ext}^i_{\mathcal{P}}(A^{(m)}, B^{(m)}) \to \mathrm{Ext}^i_{\mathcal{F}}(A^{(m)}, B^{(m)}) = \mathrm{Ext}^i_{\mathcal{F}}(A, B)$ is an isomorphism for $0 \leq i \leq s$. □

Combining Proposition 2.5 with the injectivity of Frobenius twists established in Corollary 1.3, we obtain the following theorem asserting that $\mathrm{Ext}_{\mathcal{P}}$-groups stabilize with respect to repeated Frobenius twists. The lower bound of the number of twists required to stabilize $\mathrm{Ext}^s_{\mathcal{P}}(A, B)$ depends upon both the cohomological degree $s$ and the degree $d$ of the functors $A$ and $B$; in Corollary 4.10, we shall eliminate the dependence of this bound on $d$.

THEOREM 2.6 (weak twist stability). *Let $A$ and $B$ be homogeneous strict polynomial functors of degree $d$. The Frobenius twist homomorphism*

$$\mathrm{Ext}^s_{\mathcal{P}}(A^{(m)}, B^{(m)}) \to \mathrm{Ext}^s_{\mathcal{P}}(A^{(m+1)}, B^{(m+1)})$$

*is injective for all $m$ and is an isomorphism for $m \geq \log_p \frac{s+2}{2} + [\frac{d-1}{p-1}]$.*

*Proof.* Injectivity is given by Corollary 1.3. To prove the asserted surjectivity, we pick a finite extension $E/k$ with more than $dp^{m+1}$ elements. Consider the following commutative diagram

$$\begin{CD}
\mathrm{Ext}^s_{\mathcal{P}(E)}(A^{(m)}_E, B^{(m)}_E) @>>> \mathrm{Ext}^s_{\mathcal{P}(E)}(A^{(m+1)}_E, B^{(m+1)}_E) \\
@V \cong VV @V \cong VV \\
\mathrm{Ext}^s_{\mathcal{F}(E)}(A^{(m)}_E, B^{(m)}_E) @>\cong>> \mathrm{Ext}^s_{\mathcal{F}(E)}(A^{(m+1)}_E, B^{(m+1)}_E).
\end{CD}$$

By Proposition 2.5, the vertical arrows in this diagram are isomorphisms, whereas the lower horizontal arrow is an isomorphism because Frobenius twist is an invertible functor in $\mathcal{F}(E)$. Thus, the top horizontal arrow is also an isomorphism. The proof is completed by applying Proposition 1.1: the upper horizontal arrow of this commutative square is the extension of scalars of $\mathrm{Ext}^s_{\mathcal{P}(k)}(A^{(m)}, B^{(m)}) \to \mathrm{Ext}^s_{\mathcal{P}(k)}(A^{(m+1)}, B^{(m+1)})$. □

Using Theorem 2.6, we significantly strengthen Proposition 2.5 by having the bound on $q$ (the cardinality of our base field) depend only upon the Ext-degree and the degree of the homogeneous polynomial and not also upon the the number of Frobenius twists. Theorem 2.7 follows immediately from Proposition 2.5 and Theorem 2.6.

THEOREM 2.7 (the weak comparison theorem). *Let $A$ and $B$ be strict polynomial functors homogeneous of degree $d$, over a field with $q$ elements. Let $m_0$ be the least integer $\geq \log_p \frac{s+2}{2} + [\frac{d-1}{p-1}]$. Assume that $q \geq dp^{m_0}$. Then for any $m \geq m_0$ the canonical homomorphisms*

$$\mathrm{Ext}^i_{\mathcal{P}}(A^{(m)}, B^{(m)}) \to \mathrm{Ext}^i_{\mathcal{F}}(A^{(m)}, B^{(m)}) = \mathrm{Ext}^i_{\mathcal{F}}(A, B)$$

*are isomorphisms for all $0 \leq i \leq s$.*

The following is an immediate corollary of Theorem 2.7.

COROLLARY 2.8. *In conditions and notation of Theorem 2.7, the canonical homomorphism*

$$\varinjlim_m \mathrm{Ext}^i_{\mathcal{P}}(A^{(m)}, B^{(m)}) \to \mathrm{Ext}^i_{\mathcal{F}}(A, B)$$

*is an isomorphism for $0 \leq i \leq s$, provided that $q \geq dp^{m_0}$.*

## 3. Extension of scalars and Ext-groups

The aim of this section is to understand the effect on $\mathrm{Ext}_{\mathcal{F}}$-groups of base change from $\mathcal{P}(k)$ to $\mathcal{P}(K)$, where $k \to K$ is an extension of finite fields. In Theorem 3.9, we see that this effect is simply one of base change provided that the cardinality of the base field $k$ is greater than or equal to the degree of the strict polynomial functors involved. This will enable us to provide in Theorem 3.10 a considerably stronger version of the comparison of Theorem 2.7: the condition on cardinality of the (finite) base field is weakened to the condition that the cardinality be greater or equal to the degree of the strict polynomial functors involved.

For any finite extension $k \to K$ of fields, we have pairs of functors

$$\mathcal{V}^f_K \overset{t}{\underset{\tau}{\leftrightarrows}} \mathcal{V}^f_k \qquad \mathcal{V}_K \overset{t}{\underset{\tau}{\leftrightarrows}} \mathcal{V}_k$$

where $\tau$ is the forgetful (or restriction) functor and $t$ is the extension of scalars (or induction) functor. Clearly, both $t$ and $\tau$ are exact and $t$ is left adjoint to $\tau$.

The following proposition is less evident and presumably less well known. One can informally describe this result as saying that induction equals coinduction. We include a sketch of proof for the sake of completeness.

PROPOSITION 3.1. *Let $k \to K$ be a finite extension of fields. A choice of a nonzero $k$-linear homomorphism $T: K \to k$ determines adjunction transformations*

$$\theta: \tau \circ t \to I_{\mathcal{V}_k}, \qquad \eta: I_{\mathcal{V}_K} \to t \circ \tau$$

*which establish that $t$ is right adjoint to $\tau$ (as functors on either $\mathcal{V}$ or $\mathcal{V}^f$).*

*Proof.* The choice of $T$ determines a nondegenerate symmetric $k$-bilinear form

$$\begin{aligned} K \times K &\to k \\ (a, b) &\mapsto T(ab) \,. \end{aligned}$$

We define
$$x = \sum_{i=1}^{n} e_i \otimes e_i^{\#} \in K \otimes_k K ,$$

where $e_1, \ldots, e_n$ is a basis for $K$ over $k$ and $e_1^{\#}, \ldots, e_n^{\#}$ is the dual basis with respect to the above bilinear form. For any $V$ in $\mathcal{V}_k$, we define
$$\theta_V = T \otimes_k 1_V : \tau \circ t(V) = K \otimes_k V \to V ;$$

for any $W$ in $\mathcal{V}_K$, we define
$$\eta_W = x \cdot (-) : W \to t \circ \tau(W) = K \otimes_k W$$

(i.e., multiplication by $x \in K \otimes_k K$).

One easily checks that $x$ is independent of the choice of the basis $e_1, \ldots, e_n$. This, in turn, easily implies that

(3.1.1) $\qquad x \cdot (a \otimes 1 - 1 \otimes a) = 0 \in K \otimes_k K \quad \forall a \in K,$

so that $\eta_W$ is $K$-linear. Similarly, one checks that

(3.1.2) $\qquad (T \otimes 1_K)(x) = (1_K \otimes T)(x) = 1 \in K$

which easily implies that

$$\begin{array}{cc} \mathrm{Hom}_K(W, K \otimes_k V) \to \mathrm{Hom}_k(W, V) & \mathrm{Hom}_k(W, V) \to \mathrm{Hom}_K(W, K \otimes_k V) \\ \phi \mapsto \theta_V \phi & \psi \mapsto (1_K \otimes_k \psi)\eta_W \end{array}$$

are mutually inverse natural isomorphisms. $\square$

*Remark* 3.1.3. If $K/k$ is separable, then we shall without special mention always choose $T : K \to k$ to be the trace map.

The functors $\tau$ and $t$ determine adjoint (on both sides) pairs

$$\mathcal{F}(k) = \mathrm{Funct}(\mathcal{V}_k^f, \mathcal{V}_k) \underset{\tau_*}{\overset{t_*}{\rightleftarrows}} \mathcal{F}(K/k)$$
$$= \mathrm{Funct}(\mathcal{V}_k^f, \mathcal{V}_K) \underset{t^*}{\overset{\tau^*}{\rightleftarrows}} \mathcal{F}(K) = \mathrm{Funct}(\mathcal{V}_K^f, \mathcal{V}_K) .$$

Here $t_*$ and $\tau_*$ are given by the composition with $t$ and $\tau$ on the left, whereas $t^*$ and $\tau^*$ are given by composition with $t$ and $\tau$ on the right. Note that all functors in the above diagram are exact. Since these functors have exact adjoints, we conclude further that they take injectives to injectives and projectives to projectives.

These observations immediately imply the following proposition.

PROPOSITION 3.2. *Let $k \to K$ be a finite extension of fields. For any functors $A$ in $\mathcal{F}(K)$, $B$ in $\mathcal{F}(K/k)$, $C$ in $\mathcal{F}(k)$, there are natural isomorphisms* *of graded $K$ vector-spaces*:

$$\operatorname{Ext}^*_{\mathcal{F}(K/k)}(B, t^*A) = \operatorname{Ext}^*_{\mathcal{F}(K)}(\tau^*B, A), \ \operatorname{Ext}^*_{\mathcal{F}(K/k)}(B, t_*C) = \operatorname{Ext}^*_{\mathcal{F}(k)}(\tau_*B, C),$$
$$\operatorname{Ext}^*_{\mathcal{F}(K/k)}(t^*A, B) = \operatorname{Ext}^*_{\mathcal{F}(K)}(A, \tau^*B), \ \operatorname{Ext}^*_{\mathcal{F}(K/k)}(t_*C, B) = \operatorname{Ext}^*_{\mathcal{F}(k)}(C, \tau_*B).$$

COROLLARY 3.3. *Let $k \to K$ be a finite extension of fields. Let $A \in \mathcal{F}(K)$, $A_0 \in \mathcal{F}(k)$ be functors such that $t^*A = t_*A_0 \in \mathcal{F}(K/k)$. For any $C$ in $\mathcal{F}(k)$, we have natural isomorphisms of graded $K$-vector spaces*

$$K \otimes_k \operatorname{Ext}^*_{\mathcal{F}(k)}(A_0, C) = \operatorname{Ext}^*_{\mathcal{F}(K)}(A, \tau^*t_*C),$$
$$K \otimes_k \operatorname{Ext}^*_{\mathcal{F}(k)}(C, A_0) = \operatorname{Ext}^*_{\mathcal{F}(K)}(\tau^*t_*C, A).$$

*Proof.* According to Proposition 3.2 we have natural isomorphisms

$$\operatorname{Ext}^*_{\mathcal{F}(K)}(A, \tau^*t_*C) = \operatorname{Ext}^*_{\mathcal{F}(K/k)}(t^*A, t_*C) = \operatorname{Ext}^*_{\mathcal{F}(K/k)}(t_*A_0, t_*C)$$
$$= \operatorname{Ext}^*_{\mathcal{F}(k)}(A_0, \tau_*t_*C).$$

It now suffices to observe that the functor $\tau_*t_*C$ is given by the formula $V \mapsto K \otimes_k C(V)$. If we forget about the action of $K$, then this is just a direct sum of finitely many copies of $C$, so that the natural homomorphism

$$K \otimes_k \operatorname{Ext}^*_{\mathcal{F}(k)}(A_0, C) \to \operatorname{Ext}^*_{\mathcal{F}(k)}(A_0, \tau_*t_*C)$$

is an isomorphism of graded $K$-vector spaces. □

Using Corollary 3.3, we obtain the following understanding of the effect of base change on $\operatorname{Ext}_{\mathcal{F}}$-groups, an understanding which we later compare to the corresponding property of $\operatorname{Ext}_{\mathcal{P}}$-groups, as given by Proposition 1.1.

THEOREM 3.4. *For any strict polynomial functors $P, Q$ in $\mathcal{P}(k)$ we have natural isomorphisms of graded $K$-vector spaces*

$$K \otimes_k \operatorname{Ext}^*_{\mathcal{F}(k)}(P, Q) = \operatorname{Ext}^*_{\mathcal{F}(K)}(P_K, \tau^*t^*Q_K) = \operatorname{Ext}^*_{\mathcal{F}(K)}(\tau^*t^*P_K, Q_K).$$

*Proof.* This follows immediately from Corollary 3.3, when we take into account the obvious formula

$$t^*P_K = t_*P \in \mathcal{F}(K/k).$$

□

*Remark* 3.4.1. It is easy to see from the discussion above that the isomorphisms

$$K \otimes_k \operatorname{Ext}^*_{\mathcal{F}(k)}(P, Q) = \operatorname{Ext}^*_{\mathcal{F}(K)}(P_K, \tau^*t^*Q_K) = \operatorname{Ext}^*_{\mathcal{F}(K)}(P_K, Q_K \circ (t \circ \tau)),$$
$$K \otimes_k \operatorname{Ext}^*_{\mathcal{F}(k)}(P, Q) = \operatorname{Ext}^*_{\mathcal{F}(K)}(\tau^*t^*P_K, Q_K) = \operatorname{Ext}^*_{\mathcal{F}(K)}(P_K \circ (t \circ \tau), Q_K)$$

are obtained as $K$-linear extensions of the canonical $k$-linear maps

$$\operatorname{Ext}^*_{\mathcal{F}(k)}(P,Q) \to \operatorname{Ext}^*_{\mathcal{F}(K)}(P_K, Q_K \circ (t \circ \tau)),$$
$$\operatorname{Ext}^*_{\mathcal{F}(k)}(P,Q) \to \operatorname{Ext}^*_{\mathcal{F}(K)}(P_K \circ (t \circ \tau), Q_K),$$

which may be described as follows. The exact functor $\tau^* t_* : \mathcal{F}(k) \to \mathcal{F}(K)$ defines a canonical homomorphism $\tau^* t_* : \operatorname{Ext}^*_{\mathcal{F}(k)}(P,Q) \to \operatorname{Ext}^*_{\mathcal{F}(K)}(\tau^* t_* P, \tau^* t_* Q)$
$= \operatorname{Ext}^*_{\mathcal{F}(K)}(\tau^* t^* P_K, \tau^* t^* Q_K) = \operatorname{Ext}^*_{\mathcal{F}(K)}(P_K \circ (t \circ \tau), Q_K \circ (t \circ \tau))$. Furthermore, left adjointness of $t$ to $\tau$ defines a functorial homomorphism $\gamma : t \circ \tau \to I$, whereas right adjointness of $t$ to $\tau$ defines a functorial homomorphism $\eta : I \to t \circ \tau$. The homomorphisms in question are the compositions of the homomorphism induced by $\tau^* t_*$ with homomorphisms

$$\operatorname{Ext}^*_{\mathcal{F}(K)}(P_K \circ (t \circ \tau), Q_K \circ (t \circ \tau)) @>X \mapsto X \cdot P_K(\eta)>> \operatorname{Ext}^*_{\mathcal{F}(K)}(P_K, Q_K \circ (t \circ \tau)),$$
$$\operatorname{Ext}^*_{\mathcal{F}(K)}(P_K \circ (t \circ \tau), Q_K \circ (t \circ \tau)) @>X \mapsto Q_K(\gamma) \cdot X>> \operatorname{Ext}^*_{\mathcal{F}(K)}(P_K \circ (t \circ \tau), Q_K)$$

defined by the functor homomorphisms $P_K(\eta) : P_K \to P_K \circ (t \circ \tau)$ and $Q_K(\gamma) : Q_K \circ (t \circ \tau) \to Q_K$ respectively.

The following theorem relates the effects of the base change on Ext-groups in the categories $\mathcal{P}$ and $\mathcal{F}$.

THEOREM 3.5. *Let $k \to K$ be a finite extension of fields and let $P, Q$ be in $\mathcal{P}(k)$. The following diagram commutes*

(3.51) $\qquad K \otimes_k \operatorname{Ext}^*_{\mathcal{P}(k)}(P,Q) @>\sim>> \operatorname{Ext}^*_{\mathcal{P}(K)}(P_K, Q_K)$

$\qquad\qquad\qquad\qquad @VVV @VVV$

$\qquad K \otimes_k \operatorname{Ext}^*_{\mathcal{F}(k)}(P,Q) @>\sim>> \operatorname{Ext}^*_{\mathcal{F}(K)}(P_K, Q_K \circ (t \circ \tau)).$

*Here the top horizontal arrow is the isomorphism of Proposition* 1.1, *the bottom horizontal arrow is the isomorphism of Theorem* 3.4, *the left vertical arrow is the standard homomorphism, and the right vertical arrow is the composition of the standard homomorphism* $\operatorname{Ext}^*_{\mathcal{P}(K)}(P_K, Q_K) \to \operatorname{Ext}^*_{\mathcal{F}(K)}(P_K, Q_K)$ *and the homomorphism on* Ext-*groups induced by the functor homomorphism* $Q_K(\eta) : Q_K \mapsto Q_K \circ (t \circ \tau)$.

*Similarly, there is a commutative diagram*

(3.5.2) $\qquad K \otimes_k \operatorname{Ext}^*_{\mathcal{P}(k)}(P,Q) @>\sim>> \operatorname{Ext}^*_{\mathcal{P}(K)}(P_K, Q_K)$

$\qquad\qquad\qquad\qquad @VVV @VVV$

$\qquad K \otimes_k \operatorname{Ext}^*_{\mathcal{F}(k)}(P,Q) @>\sim>> \operatorname{Ext}^*_{\mathcal{F}(K)}(P_K \circ (t \circ \tau), Q_K).$

*Proof.* Let $x$ in $\operatorname{Ext}^m_{\mathcal{P}(k)}(P,Q)$ be represented by an extension of strict polynomial functors

$$x : \quad 0 \to Q \to X_m \to \cdots \to X_1 \to P \to 0.$$

The extension $\tau^*t_*(x)$ coincides with $\tau^*t^*(x_K)$. Note further that we have a canonical homomorphism of extensions

$$\begin{array}{ccccccccc}
\tau^*t^*(x_K) & : & 0 \to & Q_K \circ (t \circ \tau) & \to & \cdots & \to & P_K \circ (t \circ \tau) & \to & 0 \\
& & & @AQ_K(\eta)AA & & @AP_K(\eta)AA & & & \\
x_K & : & 0 \to & Q_K & \to & \cdots & \to & P_K & \to & 0.
\end{array}$$

The existence of such a homomorphism gives us the desired equality

$$\tau^*t_*(x) \cdot P_K(\eta) = \tau^*t^*(x_K) \cdot P_K(\eta) = Q_K(\eta) \cdot x_K.$$

Essentially the same proof gives (3.5.2); alternatively, (3.5.2) follows by applying the duality isomorphism of (1.12) to (3.5.1). □

The following lemma will enable us to reinterpret in the special case of an extension of finite fields the Ext-groups $\mathrm{Ext}^*_{\mathcal{F}(K)}(P_K, Q_K \circ (t \circ \tau))$ occurring in the lower right corner of the commutative square (3.5.1).

LEMMA 3.6. *Let $k$ be a finite field with $q = p^N$ elements and let $k \to K$ be a finite field extension of degree $n$. The corresponding functor $t \circ \tau : \mathcal{V}_K \to \mathcal{V}_K$ coincides with $I \oplus I^{(N)} \oplus \cdots \oplus I^{((n-1)N)}$. Under this identification, the adjunction homomorphism $\eta : I \to t \circ \tau$ (respectively, $\gamma : t \circ \tau \to I$) corresponds to the canonical embedding of $I$ into the above direct sum (resp. the canonical projection of the above direct sum onto $I$). In particular, the composition $I \xrightarrow{\eta} t \circ \tau \xrightarrow{\gamma} I$ is the identity endomorphism.*

*Proof.* The functor in question is given by the formula $W \mapsto K \otimes_k W = (K \otimes_k K) \otimes_K W$. The $K-K$ bimodule $K \otimes_k K$ is canonically isomorphic to the direct sum $\oplus_{\sigma \in \mathrm{Gal}(K/k)} {}^\sigma K$, where ${}^\sigma K$ is $K$ with the standard right $K$-module structure and the left $K$-module structure is defined via $\sigma$. In case $\sigma$ is the $lN^{\mathrm{th}}$ power of the Frobenius automorphism, the $K$-vector space ${}^\sigma K \otimes_K W$ coincides with $W^{(lN)}$. These remarks allow us to identify $t \circ \tau$ with $I \oplus I^{(N)} \oplus \cdots \oplus I^{((n-1)N)}$. The homomorphism $\gamma_W : K \otimes_k W \to W$ is given by the formula $a \otimes w \mapsto aw$ and hence corresponds to the multiplication homomorphism $K \otimes_k K \to K$, i.e. to the projection of $K \otimes_k K$ onto the summand corresponding to $\sigma = \mathrm{Id}_K$. Finally the homomorphism $\eta_W : W \to K \otimes_k W$ is given by multiplication by $x \in K \otimes_k K$. The property (3.1.1) shows that $x$ is killed by the kernel of the projection $K \otimes_k K \to K$ and hence all $\sigma$-components of $x$ with $\sigma \neq \mathrm{Id}_K$ are trivial. Finally the property (3.1.2) implies easily that the $\mathrm{Id}_K$-component of $x$ is equal to 1. □

As we now show, injectivity of $\mathrm{Ext}^*_{\mathcal{P}} \to \mathrm{Ext}^*_{\mathcal{F}}$ is an easy consequence of Theorems 2.7 and 3.5 provided that we work over a finite base field. (The hy-pothesis of finiteness of the base field is required in the proof of Theorem 2.7.)

COROLLARY 3.7. *Let $k$ be a finite field and consider strict polynomial functors $P$ and $Q$ of finite degree. The canonical homomorphism*

$$\mathrm{Ext}^s_{\mathcal{P}(k)}(P,Q) \to \mathrm{Ext}^s_{\mathcal{F}(k)}(P,Q)$$

*is injective for all $s$.*

*Proof.* According to Theorem 2.7, we may find an extension $K/k$ and an integer $m \geq 0$ such that the homomorphism $\mathrm{Ext}^s_{\mathcal{P}(K)}(P_K^{(m)}, Q_K^{(m)}) \to \mathrm{Ext}^s_{\mathcal{F}(K)}(P_K, Q_K)$ is an isomorphism. Applying Theorem 3.5 to the functors $P^{(m)}$ and $Q^{(m)}$ and noting that the right vertical arrow of (3.5.1) is the composition of an isomorphism $\mathrm{Ext}^s_{\mathcal{P}(K)}(P_K^{(m)}, Q_K^{(m)}) \to \mathrm{Ext}^s_{\mathcal{F}(K)}(P_K^{(m)}, Q_K^{(m)})$ and a split monomorphism $\mathrm{Ext}^s_{\mathcal{F}(K)}(P_K^{(m)}, Q_K^{(m)}) @>X \mapsto X \cdot Q_K^{(m)}(\eta)>> \mathrm{Ext}^s_{\mathcal{F}(K)}(P_K^{(m)}, Q_K^{(m)} \circ (t \circ \tau))$, we conclude that the map $\mathrm{Ext}^s_{\mathcal{P}(k)}(P^{(m)}, Q^{(m)}) \to \mathrm{Ext}^s_{\mathcal{F}(k)}(P^{(m)}, Q^{(m)}) = \mathrm{Ext}^s_{\mathcal{F}(k)}(P,Q)$ is injective. Injectivity of the map in question now follows from the injectivity of the twist map (given in Corollary 1.3). □

We construct now an extension-of-scalars homomorphism on $\mathrm{Ext}^*_{\mathcal{F}}$-groups applied to strict polynomial functors.

PROPOSITION 3.8. *Let $k \to K$ be a finite extension of fields. For $P$ and $Q$ in $\mathcal{P}(k)$, there is a natural "extension-of-scalars" $k$-homomorphism*

$$\mathrm{Ext}^*_{\mathcal{F}(k)}(P,Q) \to \mathrm{Ext}^*_{\mathcal{F}(K)}(P_K, Q_K)$$

*which fits in a commutative square of $K$-linear maps*

(3.8.1)
$$\begin{CD} K \otimes_k \mathrm{Ext}^*_{\mathcal{P}(k)}(P,Q) @>\sim>> \mathrm{Ext}^*_{\mathcal{P}(K)}(P_K, Q_K) \\ @VVV @VVV \\ K \otimes_k \mathrm{Ext}^*_{\mathcal{F}(k)}(P,Q) @>>> \mathrm{Ext}^*_{\mathcal{F}(K)}(P_K, Q_K) \end{CD}$$

*whose upper horizontal arrow is the isomorphism of Proposition 1.1 and whose vertical maps are the standard homomorphisms.*

*Proof.* Let $\gamma : t \circ \tau \to I$ be the adjunction homomorphism establishing the left adjointness of $t$ to $\tau$. The extension-of-scalars homomorphism $\mathrm{Ext}^*_{\mathcal{F}(k)}(P,Q) \to \mathrm{Ext}^*_{\mathcal{F}(K)}(P_K, Q_K)$ is defined to be the composition of the homomorphism $\mathrm{Ext}^*_{\mathcal{F}(k)}(P,Q) \to \mathrm{Ext}^*_{\mathcal{F}(K)}(P_K, Q_K \circ (t \circ \tau))$ of Theorem 3.4 with the map

$$\mathrm{Ext}^*_{\mathcal{F}(K)}(P_K, Q_K \circ (t \circ \tau)) @>X \mapsto Q_K(\gamma) \cdot X>> \mathrm{Ext}^*_{\mathcal{F}(K)}(P_K, Q_K).$$

Theorem 3.5 and the fact that $Q_K(\gamma) \circ Q_K(\eta) = Q_K(\gamma \circ \eta) = \mathrm{Id}$ imply the commutativity of (3.8.1). $\square$

Using a weight argument together with Lemma 3.6, we now verify that the extension-of-scalars homomorphism induces an isomorphism

$$K \otimes_k \mathrm{Ext}^*_{\mathcal{F}(k)}(P,Q) \overset{\sim}{\longrightarrow} \mathrm{Ext}^*_{\mathcal{F}(K)}(P_K, Q_K)$$

provided that $P$ and $Q$ are homogeneous strict polynomial functors of degree $d$ and $q = |k| \geq d$.

THEOREM 3.9. *Let $k$ be a finite field with $q = p^N$ elements and let $k \to K$ be a degree $n$ extension of finite fields. Consider strict polynomial functors $P$, $Q$ homogeneous of degree $d$. If $q \geq d$, then the extension-of-scalars homomorphism*

$$K \otimes_k \mathrm{Ext}^*_{\mathcal{F}(k)}(P,Q) \to \mathrm{Ext}^*_{\mathcal{F}(K)}(P_K, Q_K)$$

*is an isomorphism.*

*Proof.* Theorem 3.4 shows that the left-hand side coincides with $\mathrm{Ext}^*_{\mathcal{F}(K)}(P_K, Q_K \circ (t \circ \tau))$. After this identification our homomorphism coincides with the map induced by the functor homomorphism $Q_K \circ (t \circ \tau) \overset{Q_K(\gamma)}{\longrightarrow} Q_K$. This functor homomorphism is split by $Q_K(\eta) : Q_K \to Q_K \circ (t \circ \tau)$, which gives us a direct sum decomposition $Q_K \circ (t \circ \tau) = Q_K \oplus Q_K'$. To prove the theorem, we must show that $\mathrm{Ext}^*_{\mathcal{F}(K)}(P_K, Q_K') = 0$.

Assume first that there exists an exponential strict polynomial functor $Q^*$ such that each $Q^i$ is homogeneous of degree $i$ and $Q = Q^d$. By Lemma 3.6 and the exponential property of $Q^*$, we have a natural direct sum decomposition

$$Q_K \circ (t \circ \tau) = Q_K \circ (I \oplus \cdots \oplus I^{((n-1)N)}) = \bigoplus_{d_0 + \cdots + d_{n-1} = d} Q_K^{d_0} \otimes \cdots \otimes (Q_K^{d_{n-1}})^{((n-1)N)}.$$

The summand corresponding to the $n$-tuple $(d, 0, \ldots, 0)$ is $Q_K$ and hence $Q_K'$ coincides with the direct sum over $n$-tuples $(d_0, \ldots, d_{n-1}) \neq (d, 0, \ldots, 0)$. The summand in the above direct sum decomposition indexed by $(d_0, \ldots, d_{n-1})$ is homogeneous of degree $d_0 + d_1 q + \cdots + d_{n-1} q^{n-1} = d + d_1(q-1) + \cdots + d_{n-1}(q^{n-1} - 1)$. On the other hand $P_K$ is homogeneous of degree $d$. Since all Ext-groups in the category $\mathcal{F}_K$ between two homogeneous functors are zero unless their degrees are congruent modulo $q^n - 1$ (see e.g. [K1, 3.3]), it suffices to note now that for any $n$-tuple $\neq (d, \ldots, 0)$ we have inequalities

$$0 < d_1(q-1) + \cdots + d_{n-1}(q^{n-1} - 1) \leq d(q^{n-1} - 1) \leq q(q^{n-1} - 1) < q^n - 1.$$

In the general case, we have an injective resolution $0 \to Q \to Q_\bullet$ in which each functor $Q_i$ is a direct sum of functors of the form $S^{l_1} \otimes \cdots \otimes S^{l_d}$, $l_1 + \ldots + l_d = d$. Since the functor $\mathcal{P}(k) \to \mathcal{F}(K)$, $Q \mapsto Q_K{}'$ is obviously exact, we get a resolution
$$0 \to Q'_K \to (Q_\bullet)_K{}'.$$

The above special case shows that $\operatorname{Ext}^*_{\mathcal{F}(K)}(P_K, (Q_i)_K{}') = 0$, so that considering the hypercohomology spectral sequence we easily conclude the proof. □

We now have all the ingredients in place to easily conclude our strong comparison theorem.

THEOREM 3.10 (strong comparison theorem).  *Let $k$ be a finite field with $q$ elements and further let $P$ and $Q$ be homogeneous strict polynomial functors of degree $d$. If $q \geq d$, then the canonical homomorphism*
$$\varinjlim_m \operatorname{Ext}^*_{\mathcal{P}(k)}(P^{(m)}, Q^{(m)}) \to \operatorname{Ext}^*_{\mathcal{F}(k)}(P, Q)$$
*is an isomorphism (in all degrees).*

*Proof.* For any $s \geq 0$ we can find, according to Corollary 2.8, a finite extension $K/k$ such that the standard homomorphism $\varinjlim_m \operatorname{Ext}^*_{\mathcal{P}(K)}(P_K^{(m)}, Q_K^{(m)}) \to \operatorname{Ext}^*_{\mathcal{F}(K)}(P_K, Q_K)$ is an isomorphism for all $0 \leq * \leq s$. According to Proposition 3.8 and Theorem 3.9 we have a commutative diagram

$$\begin{array}{ccc} K \otimes_k \varinjlim_m \operatorname{Ext}^*_{\mathcal{P}(k)}(P^{(m)}, Q^{(m)}) & \xrightarrow{\sim} & \varinjlim_m \operatorname{Ext}^*_{\mathcal{P}(K)}(P_K^{(m)}, Q_K^{(m)}) \\ \downarrow & & \downarrow \cong \\ K \otimes_k \operatorname{Ext}^*_{\mathcal{F}(k)}(P, Q) & \xrightarrow{\sim} & \operatorname{Ext}^*_{\mathcal{F}(K)}(P_K, Q_K) \end{array}$$

in which all arrows except possibly the left vertical one are isomorphisms (for $* \leq s$). This implies immediately that the left vertical arrow is an isomorphism as well. □

Theorem 3.10 is complemented by the following elementary observation.

LEMMA 3.11.  *Let $k$ be a finite field with $q$ elements. Assume that $P, Q \in \mathcal{P}(k)$ are homogeneous strict polynomial functors with $0 \leq \deg P \neq \deg Q < q$. Then*
$$\varinjlim_m \operatorname{Ext}^*_{\mathcal{P}(k)}(P^{(m)}, Q^{(m)}) = \operatorname{Ext}^*_{\mathcal{F}(k)}(P, Q) = 0.$$

*Proof.* Vanishing of the left-hand side is obvious and vanishing of the right hand side follows from the fact that $\deg P \not\equiv \deg Q \mod (q-1)$; cf. [K1, 3.3]. □

## 4. Computation of $\operatorname{Ext}^*_{\mathcal{P}}(\Gamma^{d(r)}, S^{dp^{r-j}(j)})$ and $\operatorname{Ext}^*_{\mathcal{P}}(\Gamma^{d(r)}, \Lambda^{dp^{r-j}(j)})$

In Theorem 4.5, we compute the Ext-groups of the title (for all $d \geq 0$ and all $0 \leq j \leq r$). Throughout this section, $\mathcal{P}$ will denote the category of strict polynomial functors of finite degree on the category of finite dimensional vector spaces over an arbitrary (but fixed) field $k$. The starting point is the determination of the Ext-groups

$$V_j = V_{r,j} = \operatorname{Ext}_{\mathcal{P}}^*(I^{(r)}, S^{p^{r-j}(j)}) \text{ and } W_j = W_{r,j} = \operatorname{Ext}_{\mathcal{P}}^*(I^{(r)}, \Lambda^{p^{r-j}(j)}).$$

This was achieved in [F-S, 4.5, 4.5.1] by use of the hypercohomology spectral sequences for the Koszul and de Rham complexes, as initiated by [F-L-S, §3, 5]. Again, our computation relies heavily on the same technique, using the complexes:

$$\Omega_{dp^{r-j+1}}^{\bullet(j-1)} : 0 \to S^{dp^{r-j+1}(j-1)} \to S^{dp^{r-j}-1(j-1)}$$
$$\otimes \Lambda^{1(j-1)} \cdots \to \Lambda^{dp^{r-j+1}(j-1)} \to 0$$
$$Kz_{dp^{r-j}}^{\bullet(j)} : 0 \to \Lambda^{dp^{r-j}(j)} \to \Lambda^{dp^{r-j}-1(j)} \otimes S^{1(j)} \cdots \to S^{dp^{r-j}(j)} \to 0.$$

Since a theorem of Cartier tells us that the cohomology of the de Rham complex has a form similar to the de Rham itself but with an additional Frobenius twist, the de Rham complex provides the structure for an inductive argument on the number $j$ of Frobenius twists. This requires, however, that we consider $S^{dp^{r-j}(j)}$ and $\Lambda^{dp^{r-j}(j)}$ simultaneously, and the Koszul complex serves to link tightly Ext-groups for symmetric powers with those for exterior powers.

Since all but the outer terms of the complexes involve the tensor product of a symmetric and exterior functor, Theorem 1.7 suggests an inductive argument on the integer $d$.

Finally, we compute by an induction on the cohomological degree the hypercohomology spectral sequences for the two complexes. When studying the second hypercohomology spectral sequences for the de Rham complex, a new feature appears, the use of the generalized Koszul complex associated to a map. Since the cohomology of such a complex is readily computed, it is useful to view the (only) nontrivial differential in our hypercohomology spectral sequence as fitting into a generalized Koszul complex. The map involved is nothing but the nontrivial differential for the case $d = 1$, which relates the graded vector spaces $V_j$ and $W_j$. The same feature reappears when one computes the first hypercohomology spectral sequence for the Koszul complex; this has only one nonzero differential, which we identify as the differential for the generalized Koszul complex associated once again to the graded map from $W_j$ to $V_j$ arising in the case $d = 1$.

**Hypercohomology spectral sequences**

In a familiar manner, there are two spectral sequences associated to applying $\mathrm{RHom}_{\mathcal{P}}(\Gamma^{d(r)}, -)$ to a complex of functors $F^\bullet$ in $\mathcal{P}$. The first has an $E_1^{n,m}$-term $\mathrm{Ext}_{\mathcal{P}}^m(\Gamma^{d(r)}, F^n)$, whereas the second has an $E_2^{n,m}$-term $\mathrm{Ext}_{\mathcal{P}}^n(\Gamma^{d(r)}, \mathrm{H}^m(F^\bullet))$. Since the Koszul complex is exact, the second spectral sequence associated to $Kz_{dp^{r-j}}^{\bullet(j)}$ is trivial (and hence the abutment of the first one is zero). We shall use the notation $\widetilde{E}$ for the first spectral sequence associated to the Koszul complex $Kz_{dp^{r-j}}^{\bullet(j)}$:

$$\widetilde{E}_1^{n,m} = \mathrm{Ext}_{\mathcal{P}}^m(\Gamma^{d(r)}, S^{n(j)} \otimes \Lambda^{dp^{r-j}-n(j)}) \Longrightarrow 0.$$

Furthermore, we shall use the notation I (resp. II) for the first (resp. second) spectral sequence associated to the de Rham complex $\Omega_{dp^{r-j+1}}^{\bullet(j-1)}$:

$$\mathrm{I}_1^{n,m} = \mathrm{Ext}_{\mathcal{P}}^m(\Gamma^{d(r)}, S^{dp^{r-j+1}-n(j-1)} \otimes \Lambda^{n(j-1)}) \Longrightarrow \mathrm{Ext}_{\mathcal{P}}^{n+m}(\Gamma^{d(r)}, \Omega_{dp^{r-j+1}}^{\bullet(j-1)}),$$

$$\mathrm{II}_2^{n,m} = \mathrm{Ext}_{\mathcal{P}}^n(\Gamma^{d(r)}, S^{dp^{r-j}-m(j)} \otimes \Lambda^{m(j)}) \Longrightarrow \mathrm{Ext}_{\mathcal{P}}^{n+m}(\Gamma^{d(r)}, \Omega_{dp^{r-j+1}}^{\bullet(j-1)}).$$

To compute the differentials in the spectral sequences $\widetilde{E}$, I, II, we shall use the following proposition.

PROPOSITION 4.1. *Let $T_1^\bullet, T_2^\bullet$ be bounded below complexes of homogeneous strict polynomial functors of degrees $d_1 p^r$ and $d_2 p^r$ respectively. If $d = d_1 + d_2$, then the first (respectively, the second) spectral sequence associated to applying the functor $\mathrm{RHom}_{\mathcal{P}}(\Gamma^{d(r)}, -)$ to the complex $T_1^\bullet \otimes T_2^\bullet$ is naturally isomorphic to the tensor product of the first (resp., second) spectral sequences associated to applying the functors $\mathrm{RHom}(\Gamma^{d_i(r)}, -)$ to $T_i^\bullet$.*

*Proof.* Let $T_1^\bullet \to I_1^{\bullet\bullet}$, $T_2^\bullet \to I_2^{\bullet\bullet}$ be the Cartan-Eilenberg resolutions of $T_1^\bullet$ and $T_2^\bullet$ respectively. An immediate verification shows that $\mathrm{Tot}(T_1^\bullet \otimes T_2^\bullet) \to I_1^{\bullet\bullet} \otimes I_2^{\bullet\bullet}$ is the Cartan-Eilenberg resolution of $\mathrm{Tot}(T_1^\bullet \otimes T_2^\bullet)$ (here $I_1^{\bullet\bullet} \otimes I_2^{\bullet\bullet}$ should be considered as a bicomplex with $(I_1 \otimes I_2)^{n,m} = \bigoplus_{n_1+n_2=n, m_1+m_2=m} I_1^{n_1,m_1} \otimes I_2^{n_2,m_2}$). The two spectral sequences in question are the two spectral sequences of the bicomplex

$$\mathrm{Hom}_{\mathcal{P}}(\Gamma^{d(r)}, I_1^{\bullet\bullet} \otimes I_2^{\bullet\bullet})$$

which can be identified by use of Theorem 1.7 with

$$\mathrm{Hom}_{\mathcal{P}}(\Gamma^{d_1(r)}, I_1^{\bullet\bullet}) \otimes \mathrm{Hom}_{\mathcal{P}}(\Gamma^{d_2(r)}, I_2^{\bullet\bullet}).$$

The result now follows immediately [M, 11, §3, Ex. 6]. □

The following corollary makes explicit the determination of the differentials in the tensor product of spectral sequences as considered in Proposition 4.1.

COROLLARY 4.2. *All differentials of the first hypercohomology spectral sequence associated to applying the functor* $\mathrm{RHom}_\mathcal{P}(\Gamma^{d(r)}, -)$ *to the complex* $(\Omega^{\bullet(j-1)}_{p^{r-j}+1})^{\otimes d}$ *are trivial.*

*In the second hypercohomology spectral sequence associated to applying the functor* $\mathrm{RHom}_\mathcal{P}(\Gamma^{d(r)}, -)$ *to* $(\Omega^{\bullet(j-1)}_{p^{r-j}+1})^{\otimes d}$ (*respectively, the first hypercohomology spectral sequence associated to applying the functor* $\mathrm{RHom}_\mathcal{P}(\Gamma^{d(r)}, -)$ *to* $(Kz^{\bullet(j)}_{p^{r-j}})^{\otimes d}$), *all differentials except for* $d_{p^{r-j}+1}$ (*resp. except for* $d_{p^{r-j}}$) *are trivial. The only nontrivial differential sends a decomposable, homogeneous element* $v_1 \otimes \cdots \otimes v_d$ *to:*

$$\sum_{i=1}^d (-1)^{\sigma(i)} v_1 \otimes \cdots \otimes \partial(v_i) \otimes \cdots \otimes v_d$$

*where* $\partial$ *is the only nonzero differential in the second hypercohomology spectral sequence associated to applying the functor* $\mathrm{RHom}_\mathcal{P}(I^{(r)}, -)$ *to* $\Omega^{\bullet(j-1)}_{p^{r-j}+1}$ (*resp. the first hypercohomology spectral sequence associated to applying* $\mathrm{RHom}_\mathcal{P}(I^{(r)}, -)$ *to* $Kz^{\bullet(j)}_{p^{r-j}}$) *and* $\sigma(i)$ *denotes the number of terms of odd total degree among* $v_1, \ldots, v_{i-1}$.

*Proof.* This follows immediately from Proposition 4.1 and the known information about differentials in the spectral sequences associated to applying the functor $\mathrm{RHom}_\mathcal{P}(I^{(r)}, -)$ to $\Omega^{\bullet(j-1)}_{p^{r-j}+1}$ and $Kz^{\bullet(j)}_{p^{r-j}}$; see [F-S, §4]. □

### Generalized Koszul complexes

The graded vector space $V_j$ is concentrated in even degrees. It is one dimensional in degrees $s \equiv 0 \bmod 2p^{r-j}$, $0 \leq s < 2p^r$, and zero otherwise [F-S, 4.5]. The graded vector space $W_j$ is related to $V_j$ by two homomorphisms. The first one, which we denote by $\theta$, is the (only nonzero) differential

$$(4.2.1) \qquad \theta := d_{p^{r-j}} : \mathrm{Ext}^*_\mathcal{P}(I^{(r)}, \Lambda^{p^{r-j}(j)}) \to \mathrm{Ext}^{*-(p^{r-j}-1)}_\mathcal{P}(I^{(r)}, S^{p^{r-j}(j)})$$

in the first hypercohomology spectral sequence associated to applying the functor $\mathrm{RHom}(I^{(r)}, -)$ to the Koszul complex $Kz^{\bullet(j)}_{p^{r-j}}$.

This homomorphism is an isomorphism, shifting degrees by $p^{r-j} - 1$, and gives an inverse to the Yoneda multiplication by the Koszul complex. This implies, in particular, that $W_j$ is concentrated in even degrees for $p \neq 2$ and in odd degrees for $p = 2$, $j < r$.

The second homomorphism from $W_j$ to $V_j$, which we denote by $\partial$, is the (only nontrivial) differential

$$(4.2.2) \qquad \partial := d_{p^{r-j}+1} : \mathrm{Ext}^*_\mathcal{P}(I^{(r)}, \Lambda^{p^{r-j}(j)}) \to \mathrm{Ext}^{*+p^{r-j}+1}_\mathcal{P}(I^{(r)}, S^{p^{r-j}(j)})$$

in the second hypercohomology spectral sequence associated to applying the functor $\mathrm{RHom}(I^{(r)}, -)$ to the de Rham complex $\Omega^{\bullet(j-1)}_{p^{r-j+1}}$. This homomorphism fits into an exact sequence

$$0 \to W_{j-1} \xrightarrow{\alpha} W_j \xrightarrow{\partial} V_j \to V_{j-1} \to 0$$

which can be interpreted in terms of edge homomorphisms of the corresponding spectral sequence (see [F-S, §4]). Here the homomorphism $V_j \to V_{j-1}$ is the natural map of degree zero induced by the embedding $S^{p^{r-j}(j)} \hookrightarrow S^{p^{r-j+1}(j-1)}$, whereas the homomorphism $\alpha$ is the Yoneda multiplication by a class in $\mathrm{Ext}_{\mathcal{P}}^{p^{r-j+1}-p^{r-j}}(\Lambda^{p^{r-j+1}(j-1)}, \Lambda^{p^{r-j}(j)})$ obtained by truncating the de Rham complex $\Omega^{\bullet(j-1)}_{p^{j-j+1}}$ at its last non-vanishing cohomology degree (namely, $p^{r-j}$) and using Cartier's Theorem. In particular, $\alpha$ increases degrees by $p^{r-j+1} - p^{r-j}$.

We shall use the notation $K_j$ (resp. $C_j$) for the kernel (resp. cokernel) of $\partial$, so that the previous exact sequence gives natural isomorphisms of graded vector spaces

$$W_{j-1}[p^{r-j} - p^{r-j+1}] \xrightarrow{\sim} K_j = \mathrm{Ker}\ \partial, \qquad \mathrm{Coker}\ \partial = C_j \xrightarrow{\sim} V_{j-1}.$$

In order to analyze the spectral sequences $\tilde{E}$ and II, we shall observe that iterations of $\theta$ and $\partial$ form "generalized Koszul complexes" in the following sense. Let $f : W \to V$ be a homomorphism of finite dimensional $k$-vector spaces. Set

$$Q^i_d(f) = S^{d-i}(V) \otimes \Lambda^i(W), \qquad Q^*_*(f) = \bigoplus_{0 \leq i \leq d} Q^i_d(f).$$

We refer to $d$ as the total degree and to $i$ as the cohomological degree.

Note that $Q^*_*(f)$ has a natural structure of a bigraded, graded commutative (with respect to the cohomological degree) algebra. We endow the algebra $Q(f)$ with a Koszul differential $\kappa$ (of cohomological degree $-1$), defined in terms of $f$ and product and coproduct operations as follows:

$$\begin{aligned}\kappa : Q^i_d(f) &= S^{d-i}(V) \otimes \Lambda^i(W) \xrightarrow{1_{S^{d-i}(V)} \otimes \mathrm{comult}} S^{d-i}(V) \otimes W \otimes \Lambda^{i-1}(W) \\ &\xrightarrow{1 \otimes f \otimes 1} S^{d-i}(V) \otimes V \otimes \Lambda^{i-1}(W) \\ &\xrightarrow{\mathrm{mult} \otimes 1_{\Lambda^{i-1}(W)}} S^{d-i+1}(V) \otimes \Lambda^{i-1}(W) = Q^{i-1}_d(f).\end{aligned}$$

We refer to the resulting complex $Q^\bullet_d(f)$ as the generalized Koszul complex for $f$.

One checks easily that $\kappa$ is a (graded) derivation and hence $\mathrm{H}^*(Q(f))$ has a natural structure of a bigraded graded commutative (with respect to the cohomological degree) algebra. Moreover, we have canonical identifications

$$\mathrm{H}^1_1(Q(f)) = \mathrm{Ker}\ f \qquad \mathrm{H}^0_1(Q(f)) = \mathrm{Coker}\ f.$$

The usefulness of this formalism lies in the following elementary lemma.

LEMMA 4.3. *The induced homomorphism of bigraded algebras*

$$S^*(\operatorname{Coker} f) \otimes \Lambda^*(\operatorname{Ker} f) \to \operatorname{H}^*(Q(f))$$

*is an isomorphism.*

*Proof.* Choose (noncanonical) splittings

$$W = \operatorname{Ker} f \oplus U, \quad V = U \oplus \operatorname{Coker} f$$

in such a way that $f$ identifies the two copies of $U$. The differential graded algebra $Q(f)$ is isomorphic to the tensor product of a trivial differential graded algebra $S^*(\operatorname{Coker} f) \otimes \Lambda^*(\operatorname{Ker} f)$ and the usual Koszul differential graded algebra corresponding to the vector space $U$. Our statement follows now from the acyclicity of the usual Koszul algebra. □

*Remark* 4.3.1. In applications, the vector spaces $V$ and $W$ are usually graded and the homomorphism $f$ preserves the grading (or shifts it by a certain integer $\deg f$). In this case $Q_*^*(f)$, and $\operatorname{H}^*(Q(f))$, acquire an additional grading which we call the *internal* grading.

The product operation in the tri-graded Hopf algebra $\operatorname{Ext}_{\mathcal{P}}^*(\Gamma^{*(r)}, S^{*(j)})$ defines a canonical homomorphism of graded vector spaces

$$V_j^{\otimes d} \to \operatorname{Ext}_{\mathcal{P}}^*(\Gamma^{d(r)}, S^{dp^{r-j}(j)}).$$

Lemma 1.11 (and the fact that $V_j$ is concentrated in even degrees) implies that this homomorphism commutes with the action of the symmetric group $\Sigma_d$, and thus factors to give a map:

(4.3.2) $$S^d(V_j) \to \operatorname{Ext}_{\mathcal{P}}^*(\Gamma^{d(r)}, S^{dp^{r-j}(j)}).$$

We shall also consider the analogous map

$$W_j^{\otimes d} \to \operatorname{Ext}_{\mathcal{P}}^*(\Gamma^{d(r)}, \Lambda^{dp^{r-j}(j)}).$$

Lemma 1.11 implies that this map commutes with the action of the symmetric group $\Sigma_d$ up to a sign, and thus factors through $\Lambda^d(W_j)$, provided that $p \neq 2$. In the case $p = 2$, the fact that this map factors through $\Lambda^d(W_j)$ is established below in Lemma 4.4. We shall employ the notation $v_1 \cdot \ldots \cdot v_i$ (resp. $w_1 \wedge \cdots \wedge w_i$) for the image of $v_1 \otimes \cdots \otimes v_i \in V_j^{\otimes i}$ (resp. $w_1 \otimes \cdots \otimes w_i \in W_j^{\otimes i}$) in $\operatorname{Ext}_{\mathcal{P}}^*(\Gamma^{i(r)}, S^{ip^{r-j}(j)})$ (resp. in $\operatorname{Ext}_{\mathcal{P}}^*(\Gamma^{i(r)}, \Lambda^{ip^{r-j}(j)})$).

LEMMA 4.4. *The natural map*

$$W_j^{\otimes d} \to \operatorname{Ext}_{\mathcal{P}}^*(\Gamma^{d(r)}, \Lambda^{dp^{r-j}(j)})$$

*factors through* $\Lambda^d(W_j)$.

*Proof.* This result can be obtained as a byproduct of our proof of Theorem 4.5 below. Such a presentation however would make the central part of the argument even harder to follow than it is right now. To avoid this, we give a direct proof of Lemma 4.4 based on the use of the canonical injective resolution of the functor $S^{*(r)}$ constructed in [F-S, §8].

As was explained above it suffices to consider the case $p = 2$ and in this case it suffices to show that the square of each homogeneous element $w$ in $W_j$ dies in $\mathrm{Ext}^*_\mathcal{P}(\Gamma^{2(r)}, \Lambda^{2^{r-j+1}(j)})$. We start with the special case $j = r$. First:

(4.4.1)  $\mathrm{Ext}^{\mathrm{odd}}_\mathcal{P}(\Gamma^{2(r)}, S^{2(r)}) = 0.$

To prove this we use the injective resolution of $S^{2(r)}$ constructed in [F-S, §8] $0 \to S^{2(r)} \to C^0 \to C^1 \to \ldots$ , where

$$C^n = \bigoplus_{\substack{m_0 + \cdots + m_{2^r-1} = 2^{r+1} \\ \sum i \cdot m_i = n \cdot 2^{r-1}}} S^{m_0} \otimes \cdots \otimes S^{m_{2^r-1}}.$$

Theorem 1.7 implies immediately that $\mathrm{Hom}_\mathcal{P}(\Gamma^{2(r)}, C^n) = 0$ if $n$ is odd.

The cohomology long exact sequence corresponding to the short exact sequence of strict polynomial functors: $0 \to \Lambda^{2(r)} \to I^{(r)} \otimes I^{(r)} \to S^{2(r)} \to 0$ now implies:

(4.4.2) *The standard embedding $\Lambda^{2(r)} \hookrightarrow I^{(r)} \otimes I^{(r)}$ induces injective maps*

$$\mathrm{Ext}^{\mathrm{even}}_\mathcal{P}(\Gamma^{2(r)}, \Lambda^{2(r)}) \hookrightarrow \mathrm{Ext}^{\mathrm{even}}_\mathcal{P}(\Gamma^{2(r)}, I^{(r)} \otimes I^{(r)}).$$

To finish the argument in the case $j = r$ it suffices to note now that the composition

$$(W_r \otimes W_r)^{\mathrm{even}} \to \mathrm{Ext}^{\mathrm{even}}_\mathcal{P}(\Gamma^{2(r)}, \Lambda^{2(r)})$$
$$\hookrightarrow \mathrm{Ext}^{\mathrm{even}}_\mathcal{P}(\Gamma^{2(r)}, I^{(r)} \otimes I^{(r)}) = (W_r \otimes W_r)^{\mathrm{even}}$$

coincides with $1 + \sigma$ (where $\sigma$ permutes the two factors) and hence is zero on elements of the form $w \otimes w$.

Assume now that $j < r$. Let $kz \in \mathrm{Ext}^{2^{r-j}-1}_\mathcal{P}(S^{2^{r-j}}, \Lambda^{2^{r-j}})$ be the extension class represented by the Koszul complex and let $\phi_{r-j} : I^{(r-j)} \hookrightarrow S^{2^{r-j}}$ denote the standard embedding. The standard arguments (see [F-S, §4]) show that the left Yoneda multiplication with $kz^{(j)}$ defines an isomorphism (inverse to $\theta$)

$$\mathrm{Ext}^*_\mathcal{P}(I^{(r)}, S^{2^{r-j}(j)}) @>kz^{(j)}>> \mathrm{Ext}^{*+2^{r-j}-1}_\mathcal{P}(I^{(r)}, \Lambda^{2^{r-j}(j)}).$$

On the other hand $\phi^{(j)}_{r-j} : I^{(r)} \hookrightarrow S^{2^{r-j}(j)}$ induces epimorphisms on $\mathrm{Ext}^*_\mathcal{P}(I^{(r)}, -)$ [F-S, §4]. This shows that each homogeneous element in $W_j = \mathrm{Ext}^*_\mathcal{P}(I^{(r)}, \Lambda^{2^{r-j}(j)})$ may be written as a Yoneda product $w = kz^{(j)} \cdot \phi^{(j)}_{r-j} \cdot x$, for a certain homogeneous element $x$ in $\operatorname{Ext}_{\mathcal{P}}^*(I^{(r)}, I^{(r)})$. With this notation we have the following formula:

$$w \otimes w = ((kz \otimes kz) \cdot (\phi_{r-j} \otimes \phi_{r-j}))^{(j)} \cdot (x \otimes x) \in \operatorname{Ext}_{\mathcal{P}}^*(I^{(r)} \otimes I^{(r)}, \Lambda^{2^{r-j}(j)} \otimes \Lambda^{2^{r-j}(j)}).$$

The image of $w \otimes w$ in $\operatorname{Ext}_{\mathcal{P}}^*(\Gamma^{2(r)}, \Lambda^{2^{r-j+1}(j)})$ is the Yoneda product $m_\Lambda^{(j)} \cdot (w \otimes w) \cdot c_\Gamma^{(r)}$, where $m_\Lambda : \Lambda^{2^{r-j}} \otimes \Lambda^{2^{r-j}} \to \Lambda^{2^{r-j+1}}$ is the multiplication homomorphism and $c_\Gamma : \Gamma^2 \hookrightarrow I \otimes I$ is the standard embedding. The long exact sequence of Ext-groups associated to the short exact sequence of strict polynomial functors $0 \to \Gamma^{2(r)} @>c^{(r)}>> I^{(r)} \otimes I^{(r)} @>m^{(r)}>> \Lambda^{2(r)} \to 0$ and the vanishing of $m_\Lambda^{(r)} \cdot (x \otimes x) \cdot c_\Gamma^{(r)}$ established above show that $(x \otimes x) \cdot c_\Gamma^{(r)}$ may be written in the form $c_\Gamma^{(r)} \cdot y$ for some $y \in \operatorname{Ext}_{\mathcal{P}}^*(\Gamma^{2(r)}, \Gamma^{2(r)})$. Now it suffices to note that the element under consideration is a right multiple of $(m_\Lambda \cdot (kz \otimes kz) \cdot (\phi_{r-j} \otimes \phi_{r-j}) \cdot c_\Gamma^{(r-j)})^{(j)}$ and $m_\Lambda \cdot (kz \otimes kz) \cdot (\phi_{r-j} \otimes \phi_{r-j}) \cdot c_\Gamma^{(r-j)} \in \operatorname{Ext}_{\mathcal{P}}^*(\Gamma^{2(r-j)}, \Lambda^{2^{r-j+1}}) = 0$; see [F-S, Prop. 5.4]. □

## Basic computation of Ext-groups

The following theorem, whose proof occupies much of the remainder of this section, formulates the result of our computation.

THEOREM 4.5.  *For all $d$ and all $0 \leq j \leq r$, the natural homomorphisms*

$$V_j^{\otimes d} \to \operatorname{Ext}_{\mathcal{P}}^*(\Gamma^{d(r)}, S^{dp^{r-j}(j)}) \quad \text{and} \quad W_j^{\otimes d} \to \operatorname{Ext}_{\mathcal{P}}^*(\Gamma^{d(r)}, \Lambda^{dp^{r-j}(j)})$$

*factor to induce isomorphisms of graded vector spaces*

$$S^d(V_j) @>\sim>> \operatorname{Ext}_{\mathcal{P}}^*(\Gamma^{d(r)}, S^{dp^{r-j}(j)}), \qquad \Lambda^d(W_j) @>\sim>> \operatorname{Ext}_{\mathcal{P}}^*(\Gamma^{d(r)}, \Lambda^{dp^{r-j}(j)}).$$

*Proof.* We proceed by a triple induction. The main induction is on the number $j$ of twists, $j \leq r$. The case $j = 0$ follows from [F-S, 2.10, 5.4] and duality.

In the sequel we shall assume that $0 < j \leq r$ and the result holds for $j-1$.

The second induction is on $d$. The cases $d = 0, 1$ of our theorem are tautological, so that we assume in the sequel that $d > 1$ and the result holds for all $d' < d$. For the induction step, we use the hypercohomology spectral sequences I, II obtained by applying the functor $\operatorname{RHom}_{\mathcal{P}}(\Gamma^{d(r)}, -)$ to the de Rham complex $\Omega_{dp^{r-j+1}}^{\bullet(j-1)}$. □

LEMMA 4.6. *All the differentials in the spectral sequence* I *are trivial and hence* $\mathrm{I}_1 = \mathrm{I}_\infty$.

*Proof.* The induction hypothesis on $j$ together with Theorem 1.7 give the following computation of the $\mathrm{I}_1$-terms:

$$\mathrm{I}_1^{n,m} = \begin{cases} 0 & \text{if } n \not\equiv 0 \mod p^{r-j+1} \\ [S^{d-s}(V_{j-1}) \otimes \Lambda^s(W_{j-1})]^m & \text{if } n = sp^{r-j+1}. \end{cases}$$

Here $[-]^m$ denotes the $m^{\text{th}}$ homogeneous component with respect to the internal grading. This formula shows that the map of spectral sequences defined by the homomorphism of complexes (given by the products in symmetric and exterior algebras)

$$(\Omega^{\bullet(j-1)}_{p^{r-j}+1})^{\otimes d} \to \Omega^{\bullet(j-1)}_{dp^{r-j}+1}$$

is surjective on the $E_1$-terms; thus the lemma follows from the first part of Corollary 4.2. □

COROLLARY 4.7.

$$\dim_k \operatorname{Ext}^n_{\mathcal{P}}(\Gamma^{d(r)}, \Omega^{\bullet(j-1)}_{p^{r-j}+1}) = \sum_{s=0}^{d} \dim_k \operatorname{I}_1^{sp^{r-j+1}, n-sp^{r-j+1}}$$

$$= \sum_{s=0}^{d} \dim_k [S^{d-s}(V_{j-1}) \otimes \Lambda^s(W_{j-1})]^{n-sp^{r-j+1}}$$

$$= \sum_{s=0}^{d} \dim_k [S^{d-s}(C_j) \otimes \Lambda^s(K_j)]^{n-sp^{r-j}}.$$

Here $[-]^{n-sp^{r-j+1}}$ stands for the homogeneous component with respect to the internal grading. The shift in degrees in the last identification comes from natural isomorphisms of graded vector spaces $W_{j-1} = K_j[p^{r-j+1} - p^{r-j}]$, $V_{j-1} = C_j$.

Consider now the spectral sequence II:

$$\operatorname{II}_2^{n,m} = \operatorname{Ext}^n_{\mathcal{P}}(\Gamma^{d(r)}, S^{dp^{r-j}-m(j)} \otimes \Lambda^{m(j)}) \Longrightarrow \operatorname{Ext}^{n+m}_{\mathcal{P}}(\Gamma^{d(r)}, \Omega^{\bullet(j-1)}_{dp^{r-j}+1}).$$

The induction hypothesis on $d$ and Theorem 1.7 allow us to identify all $\operatorname{II}_2^{n,m}$-terms except for the two edge rows $m = 0, dp^{r-j}$:

$$\operatorname{II}_2^{n,m} = \begin{cases} 0 & \text{if } m \not\equiv 0 \mod p^{r-j} \\ [S^{d-s}(V_j) \otimes \Lambda^s(W_j)]^n & \text{if } m = sp^{r-j}, \ 0 < s < d. \end{cases}$$

By dimension considerations all differentials $d_k$ with $2 \leq k \leq p^{r-j}$ are trivial. Furthermore, considering the homomorphism of the second hypercohomology spectral sequences induced by the homomorphism of complexes $(\Omega^{\bullet(j-1)}_{p^{r-j}+1})^{\otimes d} \to \Omega^{\bullet(j-1)}_{dp^{r-j}+1}$ and using Corollary 4.2 we identify the action of the differential $d_{p^{r-j}+1}$ on decomposable elements as

(4.7.1) $\quad d_{p^{r-j}+1}(v_1 \cdot \ldots \cdot v_{d-s} \otimes w_1 \wedge \cdots \wedge w_s)$

$$= \sum_{i=1}^{s} (-1)^{i-1} v_1 \cdot \ldots \cdot v_{d-s} \cdot \partial(w_i) \otimes w_1 \wedge \cdots \wedge \hat{w}_i \wedge \cdots \wedge w_s.$$

Here $v_1, \ldots, v_{d-s}$ are homogeneous elements in $V_j$, $w_1, \ldots, w_s$ are homogeneous elements in $W_j$ and $\partial : W_j \to V_j$ is the differential $d_{p^{r-j}+1}$ in the second hypercohomology spectral sequence corresponding to $\Omega^{\bullet(j)}_{p^{r-j}+1}$, discussed previously. The sign in the formula (4.7.1) is explained by the fact that all $v$'s are of even total degree, while all $w$'s are of odd total degree. Thus, $d_{p^{r-j}+1}$ may be viewed as the generalized Koszul differential in $Q^\bullet_d(\partial)$.

Lemma 4.3 together with (4.7.1) implies that we get a canonical homomorphism

$$S^d(C_j) = \mathrm{Coker}(S^{d-1}(V_j) \otimes W_j \to S^d(V_j))$$
$$\to \mathrm{II}^{*,0}_{p^{r-j}+2} = \mathrm{Coker}(S^{d-1}(V_j) \otimes W_j \to \mathrm{Ext}^*_{\mathcal{P}}(\Gamma^{d(r)}, S^{dp^{r-j}(j)})) \ .$$

LEMMA 4.8. *The resulting homomorphism*

$$S^d(C_j) \to \mathrm{II}^{*,0}_{p^{r-j}+2} \to \mathrm{II}^{*,0}_\infty \to \mathrm{Ext}^*_{\mathcal{P}}(\Gamma^{d(r)}, \Omega^{\bullet(j-1)}_{dp^{r-j}+1})$$

*is injective.*

*Proof.* It is clear from the construction that the homomorphism in question is induced by the composition

$$S^d(V_j) \to \mathrm{Ext}^*_{\mathcal{P}}(\Gamma^{d(r)}, S^{dp^{r-j}(j)}) \to \mathrm{Ext}^*_{\mathcal{P}}(\Gamma^{d(r)}, \Omega^{\bullet(j-1)}_{dp^{r-j}+1})$$

where the second arrow is defined by the homomorphism of complexes $S^{dp^{r-j}(j)} \to \Omega^{\bullet(j-1)}_{dp^{r-j}+1}$ (identifying $S^{dp^{r-j}(j)}$ with the zero-dimensional homology group of $\Omega^{\bullet(j-1)}_{dp^{r-j}+1}$). Consider also the natural homomorphism $\Omega^{\bullet(j-1)}_{dp^{r-j}+1} \to S^{dp^{r-j+1}(j-1)}$ (projection onto the zero-dimensional component). The composition of these homomorphisms coincides with the standard embedding $S^{dp^{r-j}(j)} \hookrightarrow S^{dp^{r-j+1}(j-1)}$. The induction hypothesis on $j$ shows that $\mathrm{Ext}^*_{\mathcal{P}}(\Gamma^{d(r)}, S^{dp^{r-j+1}(j-1)}) = S^d(V_{j-1})$. The composite

$$S^d(C_j) \to \mathrm{Ext}^*_{\mathcal{P}}(\Gamma^{d(r)}, \Omega^{\bullet(j-1)}_{dp^{r-j}+1}) \to \mathrm{Ext}^*_{\mathcal{P}}(\Gamma^{d(r)}, S^{dp^{r-j+1}(j-1)}) = S^d(V_{j-1})$$

is by construction the isomorphism $C_j @>\sim>> V_{j-1}$ for $d=1$ and thus (by its multiplicative nature) is an isomorphism for all $d$. □

Now we can start the final induction – that on the internal index. More precisely we are going to prove, using induction on $t$, the following statement. As in Corollary 4.7, $[-]^m$ denotes the $m^{\mathrm{th}}$ homogeneous component with respect to the internal grading.

PROPOSITION 4.9. *With the above notation and assumptions, the following assertions hold for any $t \geq -1$.*

$1_t$. *The natural map $[S^d(V_j)]^t \to \mathrm{Ext}^t_{\mathcal{P}}(\Gamma^{d(r)}, S^{dp^{r-j}(j)})$ is an isomorphism.*

$2_t$. *The natural map $[\Lambda^d(W_j)]^t \to \mathrm{Ext}^t_{\mathcal{P}}(\Gamma^{d(r)}, \Lambda^{dp^{r-j}(j)})$ is an isomorphism.*

$3_t$. $\mathrm{II}^{n,m}_{p^{r-j}+2} = \mathrm{II}^{n,m}_\infty$ *for all $n \leq t - p^{r-j} - 1$.*

*Proof.* All the statements are trivial for $t = -1$, which gives us the base for the inductive argument. Assume now that $t \geq 0$ and the statement holds for all $t' < t$.

Consider once again the spectral sequence II. Using (4.7.1), we get the following commutative diagram

(4.9.1)
$$[Q_d^{s+1}(\partial)]^{n-p^{r-j}-1} @>\kappa>> [Q_d^s(\partial)]^n @>\kappa>> [Q_d^{s-1}(\partial)]^{n+p^{r-j}+1}$$
$$@VVV @VVV @VVV$$
$$\mathrm{II}_{p^{r-j}+1}^{n-p^{r-j}-1,(s+1)p^{r-j}} @>d_{p^{r-j}+1}>> \mathrm{II}_{p^{r-j}+1}^{n,sp^{r-j}} @>d_{p^{r-j}+1}>> \mathrm{II}_{p^{r-j}+1}^{n+p^{r-j}+1,(s-1)p^{r-j}}.$$

Here the top row is a piece of the generalized Koszul complex, corresponding to the homomorphism (shifting degrees by $p^{r-j}+1$) of graded vector spaces $\partial : W_j \to V_j$. Our induction hypothesis on $d$ immediately implies that the three vertical maps in the above diagram are isomorphisms unless $s = 0, 1, d-1, d$. Examining each of these four cases separately and using our induction hypothesis on $t$, we easily conclude the following:

(4.9.2) *The induced map* $[S^{d-s}(C_j) \otimes \Lambda^s(K_j)]^n \to \mathrm{II}_{p^{r-j}+2}^{n,sp^{r-j}}$, *given as the map on cohomology of* (4.9.1),

a) *is an isomorphism provided that* $s \neq 0, 1, d-1, d$ *or* $n < t - p^{r-j} - 1$,

b) *is a monomorphism provided that* $n < t$.

Consider the terms of total degree $t$ in $\mathrm{II}_{p^{r-j}+2}$. The differentials that could possibly hit $\mathrm{II}_{p^{r-j}+2}^{t,0}$ should come from $\mathrm{II}^{n,m}$-terms with $n < t - p^{r-j} - 1$. However all such differentials are trivial in view of the inductive assumption $3_{t-1}$. Thus $\mathrm{II}_{p^{r-j}+2}^{t,0} = \mathrm{II}_\infty^{t,0}$. The same reasoning shows that $\mathrm{II}_{p^{r-j}+2}^{t-p^{r-j},p^{r-j}} = \mathrm{II}_\infty^{t-p^{r-j},p^{r-j}}$. Moreover for $s \geq 2$
$$\mathrm{II}_\infty^{t-sp^{r-j},sp^{r-j}} = \mathrm{II}_{p^{r-j}+2}^{t-sp^{r-j},sp^{r-j}} = [S^{d-s}(C_j) \otimes \Lambda^s(K_j)]^{t-sp^{r-j}}$$

according to $3_{t-1}$ and (4.9.2 a). Thus

(4.9.3)
$$\dim_k \mathrm{Ext}_{\mathcal{P}}^t(\Gamma^{d(r)}, \Omega_{p^{r-j}+1}^{\bullet(j-1)}) = \dim_k \mathrm{II}_{p^{r-j}+2}^{t,0} + \dim_k \mathrm{II}_{p^{r-j}+2}^{t-p^{r-j},p^{r-j}}$$
$$+ \sum_{s=2}^d \dim_k [S^{d-s}(C_j) \otimes \Lambda^s(K_j)]^{t-sp^{r-j}}.$$

Finally we have the following inequalities:
$$\dim_k \mathrm{II}_{p^{r-j}+2}^{t,0} \geq \dim_k [S^d(C_j)]^t, \text{ according to (4.8)},$$
$$\dim_k \mathrm{II}_{p^{r-j}+2}^{t-p^{r-j},p^{r-j}} \geq \dim_k [S^{d-1}(C_j) \otimes K_j]^{t-p^{r-j}}.$$

Comparing (4.9.3) with the result given by Corollary 4.7, we conclude that both inequalities above are actually equalities. This shows that the natural homomorphisms $[S^d(C_j)]^t \to \mathrm{II}^{t,0}_{p^{r-j}+2}$ and $[S^{d-1}(C_j) \otimes K_j]^{t-p^{r-j}} \to \mathrm{II}^{t-p^{r-j},p^{r-j}}_{p^{r-j}+2}$ are isomorphisms. Inspecting the diagram (4.9.1) for $s = 1$, $n = t - p^{r-j} - 1$ (extended by zeros to the right) and using the Five-Lemma we conclude that $[S^d(V_j)]^t \to \mathrm{Ext}^t_{\mathcal{P}}(\Gamma^{d(r)}, S^{dp^{r-j}(j)})$ is an isomorphism; in other words, $1_t$ holds. Using this fact we extend easily the conclusion made in (4.9.2a) to degrees $n \leq t - p^{r-j} - 1$:

(4.9.4) $[S^{d-s}(C_j) \otimes \Lambda^s(K_j)]^n \to \mathrm{II}^{n,sp^{r-j}}_{p^{r-j}+2}$ *is an isomorphism for all $s$ provided that $n \leq t - p^{r-j} - 1$.*

We now prove $3_t$. Let $E$ denote the second hypercohomology spectral sequence obtained by applying $\mathrm{RHom}(I^{(r)}, -)$ to the de Rham complex $\Omega^{\bullet(j-1)}_{p^{r-j+1}}$. The natural homomorphism of complexes $(\Omega^{\bullet(j-1)}_{p^{r-j+1}})^{\otimes d} \to \Omega^{\bullet(j-1)}_{dp^{r-j+1}}$ defines (in view of Proposition 4.1) a homomorphism of spectral sequences $E^{\otimes d} \to \mathrm{II}$. We conclude from (4.9.4) that the homomorphism

$$(E^{\otimes d})^{n,m}_{p^{r-j}+2} \to \mathrm{II}^{n,m}_{p^{r-j}+2}$$

is surjective for $n \leq t - p^{r-j} - 1$. Since the differentials $d_k$ with $k > p^{r-j} + 1$ are trivial in $E$, they are trivial in $E^{\otimes d}$. We conclude that all the differentials $d_k$ in II with $k > p^{r-j} + 1$ starting at $\mathrm{II}^{n,m}$-terms with $n \leq t - p^{r-j} - 1$ are trivial. So $\mathrm{II}^{n,m}_{p^{r-j}+2} = \mathrm{II}^{n,m}_{\infty}$ for these values of $n$, i.e. $3_t$ holds.

To prove that $2_t$ holds, we consider the spectral sequence $\widetilde{E}$ (the first spectral sequence corresponding to $Kz^{\bullet(j)}_{dp^{r-j}}$). The $\widetilde{E}_1$-term of this spectral sequence is given (using Theorem 1.7) by the formula

$$\widetilde{E}^{m,n}_1 = \mathrm{Ext}^n_{\mathcal{P}}(\Gamma^{d(r)}, S^{m(j)} \otimes \Lambda^{dp^{r-j}-m(j)})$$
$$= \begin{cases} 0 & \text{if } m \not\equiv 0 \mod p^{r-j} \\ [S^s(V_j) \otimes \Lambda^{d-s}(W_j)]^n & \text{if } m = sp^{r-j} \ 0 < s < d. \end{cases}$$

The differentials $d_1, \ldots, d_{p^{r-j}-1}$ are trivial by dimension considerations. Corollary 4.2 allows us to determine the action of $d_{p^{r-j}}$ on decomposable elements, so that as for (4.9.1) we obtain a commutative diagram
(4.9.5)
$$[Q^{d-s+1}_d(\theta)]^{n+p^{r-j}-1} @>\kappa>> [Q^{d-s}_d(\theta)]^n @>\kappa>> [Q^{d-s-1}_d(\theta)]^{n-p^{r-j}+1}$$
$$@VVV @VVV @VVV$$
$$\widetilde{E}^{(s-1)p^{r-j},n+p^{r-j}-1}_{p^{r-j}} @>d_{p^{r-j}}>> \widetilde{E}^{sp^{r-j},n}_{p^{r-j}} @>d_{p^{r-j}}>> \widetilde{E}^{(s+1)p^{r-j},n-p^{r-j}+1}_{p^{r-j}}.$$

Here the top row is a piece of the (acyclic) generalized Koszul complex, corresponding to the isomorphism (shifting degrees by $p^{r-j} - 1$) of graded vector spaces $\theta : W_j \to V_j$ (cf. (4.2.1)). To conclude that the bottom row of (4.9.5) is acyclic it suffices by the acyclicity of the top row to know that the central and the right-hand side vertical arrows are isomorphisms. This remark plus our induction hypotheses on $d$ and $t$ imply:

(4.9.6) $\widetilde{E}^{sp^{r-j},n}_{p^{r-j}+1} = 0$ provided that $s \neq 0$, $d-1$, $d$, or $n < t$.

For $k \geq p^{r-j} + 1$, the differential $d_k : \widetilde{E}^{0,t}_k \to \widetilde{E}^{k,t-k+1}_k$ is trivial by (4.9.6). Thus $\widetilde{E}^{0,t}_{p^{r-j}+1} = \widetilde{E}^{0,t}_\infty = 0$. Inspecting finally the diagram (4.9.5) for $s = 1, n = t-p^{r-j}+1$ (extended by zeros to the left) and using the Five-Lemma we conclude that

$$[\Lambda^d(W_j)]^t \to \mathrm{Ext}^t_{\mathcal{P}}(\Gamma^{d(r)}, \Lambda^{dp^{r-j}}(j)) = \widetilde{E}^{0,t}_{p^{r-j}}$$

is an isomorphism. This shows that $2_t$ also holds, thus completing the induction step, the proof of Proposition 4.9 and thus of Theorem 4.5. □

COROLLARY 4.10 (strong twist stability). *Let $A$ and $B$ be homogeneous strict polynomial functors of degree $d$. The Frobenius twist map*

$$\mathrm{Ext}^s_{\mathcal{P}}(A^{(m)}, B^{(m)}) \to \mathrm{Ext}^s_{\mathcal{P}}(A^{(m+1)}, B^{(m+1)})$$

*is an isomorphism provided that $m \geq \log_p \frac{s+1}{2}$.*

*Proof.* Consider first the special case $A = \Gamma^d$, $B = S^d$. Theorem 4.5 shows that

$$\mathrm{Ext}^*_{\mathcal{P}}(\Gamma^{d(m)}, S^{d(m)}) = S^d(\mathrm{Ext}^*_{\mathcal{P}}(I^{(m)}, I^{(m)})).$$

Since the Frobenius twist map $\mathrm{Ext}^*_{\mathcal{P}}(I^{(m)}, I^{(m)}) \to \mathrm{Ext}^*_{\mathcal{P}}(I^{(m+1)}, I^{(m+1)})$ is an isomorphism in degrees $\leq 2p^m - 1$ [F-S;4.9], we conclude that the induced map $S^d(\mathrm{Ext}^*_{\mathcal{P}}(I^{(m)}, I^{(m)})) \to S^d(\mathrm{Ext}^*_{\mathcal{P}}(I^{(m+1)}, I^{(m+1)}))$ is an isomorphism in the same range.

Consider next the special case $A = \Gamma^{d_1} \otimes \cdots \otimes \Gamma^{d_a}$, $B = S^{d'_1} \otimes \cdots \otimes S^{d'_b}$, $d_1 + \cdots + d_a = d'_1 + \cdots + d'_b = d$. In this case, Corollary 1.8 gives us the following formula

$$\mathrm{Ext}^*_{\mathcal{P}}(A^{(m)}, B^{(m)}) = \bigoplus_{\substack{d_{11}+\cdots+d_{1b}=d_1 \\ \cdots \\ d_{a1}+\cdots+d_{ab}=d_a}} \bigoplus_{\substack{d_{11}+\cdots+d_{a1}=d'_1 \\ \cdots \\ d_{1b}+\cdots+d_{ab}=d'_b}} \bigotimes_{i=1,j=1}^{a\ \ b} \mathrm{Ext}^*_{\mathcal{P}}(\Gamma^{d_{ij}(m)}, S^{d_{ij}(m)}).$$

Thus the result in this case follows from the previous one.

In the general case, we can find resolutions

$$0 \leftarrow A \leftarrow P_\bullet \qquad 0 \to B \to I^\bullet$$

in which each $P_i$ (resp. each $I^j$) is a direct sum of functors of the form $\Gamma^{d_1} \otimes \cdots \otimes \Gamma^{d_a}$ (resp. of the form $S^{d'_1} \otimes \cdots \otimes S^{d'_b}$). The result now follows by applying the standard comparison theorem to the map of the hypercohomology spectral sequences converging to $\operatorname{Ext}^*_{\mathcal{P}}(A^{(m)}, B^{(m)})$ and $\operatorname{Ext}^*_{\mathcal{P}}(A^{(m+1)}, B^{(m+1)})$ respectively. □

## 5. Ext-groups between the classical functors in the category $\mathcal{P}$

In this section we continue the computation of Ext-groups in the category $\mathcal{P}$ of strict polynomial functors (over an arbitrary fixed field $k$) between various classical functors. As the reader will see, the results obtained in this section are (relatively) easy applications of Theorem 4.5. Again, we begin by fixing a positive integer $r$.

The general philosophy behind these computations is as follows. There is a sort of hierarchy between the functors $\Gamma^{()}$, $\Lambda^{()}$ and $S^{()}$: $\Gamma^{()}$ has a tendency to be projective (it really is if there are no twists), $S^{()}$ has a tendency to be injective and $\Lambda^{()}$ is in between. We can provide a complete computation for the Ext-groups from a more projective functor to a less projective one but not the other way round. More specifically, we compute all Ext-groups from $\Gamma^{()}$ (with any number of twists) to $\Gamma^{()}$, $\Lambda^{()}$ and $S^{()}$, we compute all Ext-groups from $\Lambda^{()}$ to $\Lambda^{()}$ and $S^{()}$ and finally from $S^{()}$ to itself; we suspect that there is no comparable computation of the Ext-groups from $S^{()}$ to $\Lambda^{()}$ or $\Gamma^{()}$ or Ext-groups from $\Lambda^{()}$ to $\Gamma^{()}$ (actually it seems that there is just no easy answer for these Ext-groups).

For nonnegative integers $0 \leq j \leq r$, set $U_j = U_{r,j} = \operatorname{Ext}_{\mathcal{P}}(\Gamma^{p^{r-j}(j)}, I^{(r)})$. Observe that $U_r = V_r = W_r = \operatorname{Ext}^*_{\mathcal{P}}(I^{(r)}, I^{(r)})$; moreover, the duality isomorphism (1.12) defines canonical isomorphisms of graded vector spaces

$$U_j = \operatorname{Ext}^*_{\mathcal{P}}(\Gamma^{p^{r-j}(j)}, I^{(r)}) \xrightarrow{\sim} \operatorname{Ext}^*_{\mathcal{P}}(I^{(r)}, S^{p^{r-j}(j)}) = V_j$$

for any $j \leq r$.

For any nonnegative integer $d$ the product operations define canonical homomorphisms

$$U_j^{\otimes d} = \operatorname{Ext}^*_{\mathcal{P}}(\Gamma^{dp^{r-j}(j)}, I^{(r)} \otimes \cdots \otimes I^{(r)}) \to \operatorname{Ext}^*_{\mathcal{P}}(\Gamma^{dp^{r-j}(j)}, S^{d(r)}),$$
$$U_j^{\otimes d} = \operatorname{Ext}^*_{\mathcal{P}}(\Gamma^{dp^{r-j}(j)}, I^{(r)} \otimes \cdots \otimes I^{(r)}) \to \operatorname{Ext}^*_{\mathcal{P}}(\Gamma^{dp^{r-j}(j)}, \Lambda^{d(r)}).$$

Using these homomorphisms, we formulate our first result, thereby completing the computation of Theorem 4.5 by considering the situation in which the number of twists on the source functor is less (or equal) to the number of twists on the target functor.

THEOREM 5.1. *For any nonnegative integer $d$ the above maps define natural isomorphisms*

$$S^d(U_j) \xrightarrow{\sim} \operatorname{Ext}^*_{\mathcal{P}}(\Gamma^{dp^{r-j}(j)}, S^{d(r)}),$$
$$\Lambda^d(U_j) \xrightarrow{\sim} \operatorname{Ext}^*_{\mathcal{P}}(\Gamma^{dp^{r-j}(j)}, \Lambda^{d(r)}).$$

*Proof.* The first isomorphism follows immediately from Theorem 4.5, in view of the following commutative diagram:

$$\begin{array}{ccc} S^d(U_j) & \longrightarrow & \operatorname{Ext}^*_{\mathcal{P}}(\Gamma^{dp^{r-j}(j)}, S^{d(r)}) \\ \cong \downarrow & & \cong \downarrow \\ S^d(V_j) & \xrightarrow{\sim} & \operatorname{Ext}^*_{\mathcal{P}}(\Gamma^{d(r)}, S^{dp^{r-j}(j)}). \end{array}$$

Here the vertical arrows are the duality isomorphisms of (1.12), the bottom horizontal arrow is the isomorphism of Theorem 4.5 and the top horizontal arrow is the map under consideration.

To prove the second isomorphism, we proceed by induction on $d$. Consider the first hypercohomology spectral sequence corresponding to the Koszul complex

$$Kz_d^{\bullet(r)}: \quad 0 \to \Lambda^{d(r)} \to \Lambda^{d-1(r)} \otimes S^{1(r)} \to \cdots \to \Lambda^{1(r)} \otimes S^{d-1(r)} \to S^{d(r)} \to 0,$$

$$E_1^{m,n} = \operatorname{Ext}^n_{\mathcal{P}}(\Gamma^{dp^{r-j}(j)}, \Lambda^{d-m(r)} \otimes S^{m(r)}) \Longrightarrow 0.$$

The inductive assumption on $d$ and the already-proved part of the theorem give the computation of all the $E_1$-terms except those in the first column. Moreover, we have the following commutative diagram of complexes

$$\begin{array}{ccccccccc} 0 & \longrightarrow & \Lambda^d(U_j) & \xrightarrow{\kappa} & \Lambda^{d-1}(U_j) \otimes U_j & \xrightarrow{\kappa} & \cdots & \xrightarrow{\kappa} & S^d(U_j) & \longrightarrow & 0 \\ & & \downarrow \cong & & \downarrow \cong & & & & \downarrow \cong \\ 0 & \longrightarrow & E_1^{0,*} & \xrightarrow{d_1} & E_1^{1,*} & \xrightarrow{d_1} & \cdots & \xrightarrow{d_1} & E_1^{d,*} & \longrightarrow & 0. \end{array}$$

Here the top row is the Koszul complex corresponding to the vector space $U_j$, and the bottom row is the $E_1$-term of our spectral sequence. The acyclicity of the Koszul complex implies immediately that $E_2^{m,*} = 0$ for $m > 1$. This implies further vanishing of all higher differentials and hence the identification $E_2 = E_\infty = 0$. Thus the bottom complex in the above diagram is also acyclic, which gives a natural identification $\operatorname{Ext}^*_{\mathcal{P}}(\Gamma^{dp^{r-j}(j)}, \Lambda^{d(r)}) = E_1^{0,*} = \Lambda^d(U_j)$. $\square$

The duality isomorphism (1.12) together with Theorems 4.5 and 5.1 immediately leads to the following statements.

COROLLARY 5.2. *For all $0 \leq j \leq r$ and all $d \geq 0$ there are natural isomorphisms*

$$\Lambda^d(\mathrm{Ext}^*_\mathcal{P}(\Lambda^{p^{r-j}(j)}, I^{(r)})) @>\sim>> \mathrm{Ext}^*_\mathcal{P}(\Lambda^{dp^{r-j}(j)}, S^{d(r)}) ,$$
$$\Lambda^d(\mathrm{Ext}^*_\mathcal{P}(I^{(r)}, S^{p^{r-j}(j)})) @>\sim>> \mathrm{Ext}^*_\mathcal{P}(\Lambda^{d(r)}, S^{dp^{r-j}(j)}).$$

The computation of $\mathrm{Ext}^*_\mathcal{P}(\Gamma^{()}, \Gamma^{()})$ in Theorem 5.4 below requires a determination of the coproduct (rather than product) operation in $\mathrm{Ext}^*_\mathcal{P}(\Gamma^{()}, \Lambda^{()})$. Recall that comultiplication in the exterior algebra defines natural homomorphisms (for any $0 \leq j \leq r$ and any $d \geq 0$)

$$\Lambda^{dp^{r-j}(j)} \hookrightarrow \Lambda^{p^{r-j}(j)} \otimes \cdots \otimes \Lambda^{p^{r-j}(j)} ,$$
$$\Lambda^{d(r)} \hookrightarrow I^{(r)} \otimes \cdots \otimes I^{(r)}$$

which induce homomorphisms on Ext-groups

$$\mathrm{Ext}^*_\mathcal{P}(\Gamma^{d(r)}, \Lambda^{dp^{r-j}(j)}) \to \mathrm{Ext}^*_\mathcal{P}(\Gamma^{d(r)}, \Lambda^{p^{r-j}(j)} \otimes \cdots \otimes \Lambda^{p^{r-j}(j)}) = W_j^{\otimes d},$$
$$\mathrm{Ext}^*_\mathcal{P}(\Gamma^{dp^{r-j}(j)}, \Lambda^{d(r)}) \to \mathrm{Ext}^*_\mathcal{P}(\Gamma^{dp^{r-j}(j)}, I^{(r)} \otimes \cdots \otimes I^{(r)}) = U_j^{\otimes d}.$$

COROLLARY 5.3. *For all $0 \leq j \leq r$ and all $d \geq 0$ the above coproduct operations induce natural isomorphisms of graded vector spaces*:

$$\mathrm{Ext}^*_\mathcal{P}(\Gamma^{d(r)}, \Lambda^{dp^{r-j}(j)}) @>\sim>> \Lambda^d(W_j),$$
$$\mathrm{Ext}^*_\mathcal{P}(\Gamma^{dp^{r-j}(j)}, \Lambda^{d(r)}) @>\sim>> \Lambda^d(U_j).$$

*Proof.* A straightforward computation shows that the composite maps

$$\Lambda^d(W_j) @>\sim>> \mathrm{Ext}^*_\mathcal{P}(\Gamma^{d(r)}, \Lambda^{dp^{r-j}(j)}) \to W_j^{\otimes d} ,$$
$$\Lambda^d(U_j) @>\sim>> \mathrm{Ext}^*_\mathcal{P}(\Gamma^{dp^{r-j}(j)}, \Lambda^{d(r)}) \to U_j^{\otimes d}$$

coincide with the corresponding natural embeddings. $\square$

Now we turn to the computation of $\mathrm{Ext}^*_\mathcal{P}(\Gamma^{()}, \Gamma^{()})$. We introduce one more basic Ext-space (defined for all $0 \leq j \leq r$) $\widetilde{V}_j = \widetilde{V}_{r,j} = \mathrm{Ext}^*_\mathcal{P}(I^{(r)}, \Gamma^{p^{r-j}(j)})$. The differential $d_{p^{r-j}}$ in the first hypercohomology spectral sequence, corresponding to the dual Koszul complex

$$Kz^{\#(j)}_{p^{r-j}}: \quad 0 \to \Gamma^{p^{r-j}(j)} \to \Gamma^{p^{r-j}-1(j)} \otimes \Lambda^{1(j)} \to \cdots \to \Lambda^{p^{r-j}(j)} \to 0$$

defines an isomorphism

$$\eta: \widetilde{V}_j @>\sim>> W_j$$

of graded vector spaces (shifting degrees by $p^{r-j} - 1$).

In the same way as above, comultiplication in the divided power algebra defines natural homomorphisms

$$\operatorname{Ext}^*_{\mathcal{P}}(\Gamma^{d(r)}, \Gamma^{dp^{r-j}(j)}) \to \operatorname{Ext}^*_{\mathcal{P}}(\Gamma^{d(r)}, \Gamma^{p^{r-j}(j)} \otimes \cdots \otimes \Gamma^{p^{r-j}(j)})$$
$$= \operatorname{Ext}^*_{\mathcal{P}}(I^{(r)}, \Gamma^{p^{r-j}(j)})^{\otimes d} = \widetilde{V}_j^{\otimes d} \ ,$$
$$\operatorname{Ext}^*_{\mathcal{P}}(\Gamma^{dp^{r-j}(j)}, \Gamma^{d(r)}) \to \operatorname{Ext}^*_{\mathcal{P}}(\Gamma^{dp^{r-j}(j)}, I^{(r)} \otimes \cdots \otimes I^{(r)})$$
$$= \operatorname{Ext}^*_{\mathcal{P}}(\Gamma^{p^{r-j}(j)}, I^{(r)})^{\otimes d} = U_j^{\otimes d}.$$

Lemma 1.11 shows that these homomorphisms are $\Sigma_d$-invariant and hence their images are contained in the subspaces of $\Sigma_d$-invariant tensors, i.e. in $\Gamma^d(\widetilde{V}_j)$ and $\Gamma^d(U_j)$ respectively.

THEOREM 5.4. *For all $0 \le j \le r$ and all $d \ge 0$ the above constructed homomorphisms*

$$\operatorname{Ext}^*_{\mathcal{P}}(\Gamma^{d(r)}, \Gamma^{dp^{r-j}(j)}) \to \Gamma^d(\widetilde{V}_j) \ ,$$
$$\operatorname{Ext}^*_{\mathcal{P}}(\Gamma^{dp^{r-j}(j)}, \Gamma^{d(r)}) \to \Gamma^d(U_j)$$

*are isomorphisms.*

*Proof.* The proof is essentially the same in both cases. We consider only the first one (which is slightly more complicated). We proceed by induction on $d$. Consider the first hypercohomology spectral sequence corresponding to the dual Koszul complex

$$Kz^{\#(j)}_{dp^{r-j}}: \quad 0 \to \Gamma^{dp^{r-j}(j)} \to \Gamma^{dp^{r-j}-1(j)} \otimes \Lambda^{1(j)} \to \cdots \to \Lambda^{dp^{r-j}(j)} \to 0 \ ,$$

$$E_1^{m,n} = \operatorname{Ext}^n_{\mathcal{P}}(\Gamma^{d(r)}, \Gamma^{dp^{r-j}-m(j)} \otimes \Lambda^{m(j)}) \Longrightarrow 0.$$

The columns $E_1^{m,*}$ of this spectral sequence with $m \not\equiv 0 \mod p^{r-j}$ are trivial. Furthermore

$$E_1^{sp^{r-j},*} = \operatorname{Ext}^*_{\mathcal{P}}(\Gamma^{d(r)}, \Gamma^{(d-s)p^{r-j}(j)} \otimes \Lambda^{sp^{r-j}(j)})$$
$$= \operatorname{Ext}^*_{\mathcal{P}}(\Gamma^{d-s(r)}, \Gamma^{(d-s)p^{r-j}(j)}) \otimes \operatorname{Ext}^*_{\mathcal{P}}(\Gamma^{s(r)}, \Lambda^{sp^{r-j}(j)}).$$

Our induction hypothesis on $d$ and Corollary 5.3 show further that comultiplication defines (for $s > 0$) natural isomorphisms $E_1^{sp^{r-j},*} @>\sim>> \Gamma^{d-s}(\widetilde{V}_j) \otimes \Lambda^s(W_j)$. The first (and the only) nontrivial differential in this spectral sequence is
$d_{p^{r-j}} : E_{p^{r-j}}^{sp^{r-j},*} \to E_{p^{r-j}}^{(s+1)p^{r-j},*}$. Considering the natural homomorphism of complexes $Kz^{\#(j)}_{dp^{r-j}} \to (Kz^{\#(j)}_{p^{r-j}})^{\otimes d}$ (defined by comultiplication in the divided power and exterior algebras) and using Proposition 4.1, we finally get a commutative diagram of complexes

$$0 \to E^{0,*}_{p^{r-j}} \xrightarrow{d_{p^{r-j}}} E^{p^{r-j},*}_{p^{r-j}} \xrightarrow{d_{p^{r-j}}} \cdots \xrightarrow{d_{p^{r-j}}} E^{dp^{r-j},*}_{p^{r-j}} \to 0$$

$$\downarrow \quad\quad \downarrow \cong \quad\quad \downarrow \cong \quad\quad \downarrow$$

$$0 \to \Gamma^d(\widetilde{V}_j) \to \Gamma^{d-1}(\widetilde{V}_j) \otimes W_j \to \cdots \to \Lambda^d(W_j) \to 0.$$

Here the bottom row is the dual of the generalized Koszul complex for the isomorphism of graded vector spaces $\eta : \widetilde{V}_j \to W_j$. Since the bottom row is acyclic, we conclude that $E^{sp^{r-j},*}_{p^{r-j}+1} = 0$ for $s \geq 1$. This implies the vanishing of all higher differentials and hence gives the following relations : $0 = E^{0,*}_\infty = E^{0,*}_{p^{r-j}+1}$. Thus the top row in the above diagram is also acyclic and hence the left-hand side vertical arrow is also an isomorphism. □

The following corollary follows immediately from Theorem 5.4 and the duality isomorphism (1.12).

COROLLARY 5.5. *For any $0 \leq j \leq r$ and any $d \geq 0$, there are natural isomorphisms*

$$\operatorname{Ext}^*_{\mathcal{P}}(S^{dp^{r-j}(j)}, S^{d(r)}) \xrightarrow{\sim} \Gamma^d(\operatorname{Ext}^*_{\mathcal{P}}(S^{p^{r-j}(j)}, I^{(r)})) \, ,$$

$$\operatorname{Ext}^*_{\mathcal{P}}(S^{d(r)}, S^{dp^{r-j}(j)}) \xrightarrow{\sim} \Gamma^d(\operatorname{Ext}^*_{\mathcal{P}}(I^{(r)}, S^{p^{r-j}(j)})).$$

To compute $\operatorname{Ext}^*(\Lambda^{()}, \Lambda^{()})$ we proceed in the same way as above. For any $0 \leq j \leq r$ and any $d \geq 0$, the comultiplication in the exterior algebra defines natural homomorphisms

$$\operatorname{Ext}^*_{\mathcal{P}}(\Lambda^{d(r)}, \Lambda^{dp^{r-j}(j)}) \to \operatorname{Ext}^*_{\mathcal{P}}(\Lambda^{d(r)}, \Lambda^{p^{r-j}(j)} \otimes \cdots \otimes \Lambda^{p^{r-j}(j)})$$

$$= \operatorname{Ext}^*_{\mathcal{P}}(I^{(r)}, \Lambda^{p^{r-j}(j)})^{\otimes d} \, ,$$

$$\operatorname{Ext}^*_{\mathcal{P}}(\Lambda^{dp^{r-j}(j)}, \Lambda^{d(r)}) \to \operatorname{Ext}^*_{\mathcal{P}}(\Lambda^{dp^{r-j}(j)}, I^{(r)} \otimes \cdots \otimes I^{(r)})$$

$$= \operatorname{Ext}^*_{\mathcal{P}}(\Lambda^{p^{r-j}(j)}, I^{(r)})^{\otimes d}.$$

By Lemma 1.11, these homomorphisms are $\Sigma_d$-invariant and hence their images are contained in $\Gamma^d(\operatorname{Ext}^*_{\mathcal{P}}(I^{(r)}, \Lambda^{p^{r-j}(j)}))$ and $\Gamma^d(\operatorname{Ext}^*_{\mathcal{P}}(\Lambda^{p^{r-j}(j)}, I^{(r)})$ respectively.

We may now conclude the following theorem in a manner strictly parallel to that for Theorem 5.4. We leave details to the reader.

THEOREM 5.6. *For all $0 \leq j \leq r$ and all $d \geq 0$ the natural homomorphisms*

$$\mathrm{Ext}^*_{\mathcal{P}}(\Lambda^{d(r)}, \Lambda^{dp^{r-j}(j)}) \to \Gamma^d(\mathrm{Ext}^*_{\mathcal{P}}(I^{(r)}, \Lambda^{p^{r-j}(j)})) ,$$

$$\mathrm{Ext}^*_{\mathcal{P}}(\Lambda^{dp^{r-j}(j)}, \Lambda^{d(r)}) \to \Gamma^d(\mathrm{Ext}^*_{\mathcal{P}}(\Lambda^{p^{r-j}(j)}, I^{(r)}))$$

*are isomorphisms.*

We phrase our next result in terms of a chosen vector-space basis in our Ext-groups.

For an integer $j$, let $\phi_j$ denote the natural embedding $\phi_j : I^{(j)} \to S^{p^j}$ and let $\phi_j^{\#}$ denote the dual homomorphism $\phi_j^{\#} : \Gamma^{p^j} \to I^{(j)}$. We recall from [F-S, 4.5.6] that composition with $\phi_{r-j}^{(j)}$, $j \leq r$, induces an isomorphism

$$\mathrm{Ext}^s_{\mathcal{P}}(I^{(r)}, I^{(r)}) \to \mathrm{Ext}^s_{\mathcal{P}}(I^{(r)}, S^{p^{r-j}(j)})$$

for $s \equiv 0 \bmod 2p^{r-j}$. Let $kz_j$ denote the extension class in $\mathrm{Ext}^{p^j-1}_{\mathcal{P}}(S^{p^j}, \Lambda^{p^j})$ represented by the Koszul complex and let $kz_j^{\#}$ denote the extension class in $\mathrm{Ext}^{p^j-1}_{\mathcal{P}}(\Lambda^{p^j}, \Gamma^{p^j})$ represented by the dual Koszul complex. Define $e_r$ to be the class in $\mathrm{Ext}^{2p^{r-1}}_{\mathcal{P}}(I^{(r)}, I^{(r)})$ defined by the following equation in $\mathrm{Ext}^{2p^{r-1}}_{\mathcal{P}}(I^{(r)}, S^{p^{r-1}(1)})$ (cf. [F-S, §4]):

$$(5.7) \qquad \phi^{(1)}_{r-1} \cdot e_r = d_{p^{r-1}+1}(kz^{(1)}_{r-1} \cdot \phi^{(1)}_{r-1}) ,$$

where $d_{p^{r-1}+1}$ is the differential (applied to the $E^{p^{r-1}-1, p^{r-1}}_{p^{r-1}}$-term) of the second hypercohomology spectral sequence associated to applying $\mathrm{RHom}(I^{(r)}, -)$ to the de Rham complex $\Omega^{\bullet}_{p^r}$.

Recall that the graded vector space $\mathrm{Ext}^*_{\mathcal{P}}(I^{(r)}, I^{(r)})$ is one dimensional in even degrees $< 2p^r$ and is zero otherwise. We use the basis of $\mathrm{Ext}^{2m}_{\mathcal{P}}(I^{(r)}, I^{(r)})$ given by the Yoneda products

$$e_r(m) = e_1^{m_0(r-1)} \cdots e_r^{m_{r-1}} ,$$
$$m = m_0 + m_1 p + \cdots + m_{r-1} p^{r-1} \ (0 \leq m_i < p).$$

Thus, $e_r(m)^{(1)} = e_{r+1}(m)$ for all $r$ and $0 \leq m < p^r$. It is shown in [F-S, §4] that the elements $\phi^{(j)}_{r-j} e_r(mp^{r-j})$, $0 \leq m < p^j$ form a basis of the vector space $V_j = \mathrm{Ext}^*_{\mathcal{P}}(I^{(r)}, S^{p^{r-j}(j)})$. By duality, dual elements $e_r(mp^{r-j}) \phi^{\#(j)}_{r-j}$ form a basis of $U_j = \mathrm{Ext}^*_{\mathcal{P}}(\Gamma^{p^{r-j}(j)}, I^{(r)})$.

In the following theorem, we summarize the computations of Theorem 4.5 and the preceding results of this section, augmenting these computations by explicitly describing the Hopf algebra structure (cf. Lemma 1.11).

THEOREM 5.8. *Let $j$ and $r$ be integers, $0 \leq j \leq r$.*

(1) *The tri-graded Hopf algebra*
$$\operatorname{Ext}_{\mathcal{P}}^*(\Gamma^{*(j)}, S^{*(r)})$$
*is a primitively generated polynomial algebra on generators $e_r(mp^{r-j})\phi_{r-j}^{\#(j)}$ in $\operatorname{Ext}_{\mathcal{P}}^{2p^{r-j}m}(\Gamma^{p^{r-j}(j)}, I^{(r)})$, $0 \leq m < p^j$, of tri-degree $(2p^{r-j}m, p^{r-j}, 1)$.*

*Similarly, the tri-graded Hopf algebra*
$$\operatorname{Ext}_{\mathcal{P}}^*(\Gamma^{*(r)}, S^{*(j)})$$
*is a primitively generated polynomial algebra on generators $\phi_{r-j}^{(j)}e_r(mp^{r-j})$ in $\operatorname{Ext}_{\mathcal{P}}^{2p^{r-j}m}(I^{(r)}, S^{p^{r-j}(j)})$, $0 \leq m < p^j$, of tri-degree $(2p^{r-j}m, 1, p^{r-j})$.*

(2) *The tri-graded Hopf algebra*
$$\operatorname{Ext}_{\mathcal{P}}^*(\Gamma^{*(j)}, \Lambda^{*(r)})$$
*is a primitively generated exterior algebra on generators $e_r(mp^{r-j})\phi_{r-j}^{\#(j)}$ in $\operatorname{Ext}_{\mathcal{P}}^{2p^{r-j}m}(\Gamma^{p^{r-j}(j)}, I^{(r)})$, $0 \leq m < p^j$, of tri-degree $(2p^{r-j}m, p^{r-j}, 1)$.*

*A similar statement holds by duality for the tri-graded Hopf algebra*
$$\operatorname{Ext}_{\mathcal{P}}^*(\Lambda^{*(r)}, S^{*(j)}).$$

(3) *The tri-graded Hopf algebra*
$$\operatorname{Ext}_{\mathcal{P}}^*(\Gamma^{*(r)}, \Lambda^{*(j)})$$
*is a primitively generated exterior algebra on generators $kz_{r-j}^{(j)}\phi_{r-j}^{(j)}e_r(mp^{r-j})$ in $\operatorname{Ext}_{\mathcal{P}}^{2mp^{r-j}+p^{r-j}-1}(I^{(r)}, \Lambda^{p^{r-j}(j)})$, $0 \leq m < p^j$, of tri-degree $(2mp^{r-j} + p^{r-j} - 1, 1, p^{r-j})$.*

*A similar statement holds by duality for the tri-graded Hopf algebra*
$$\operatorname{Ext}_{\mathcal{P}}^*(\Lambda^{*(j)}, S^{*(r)}).$$

(4) *The tri-graded Hopf algebra*
$$\operatorname{Ext}_{\mathcal{P}}^*(\Gamma^{*(r)}, \Gamma^{*(j)})$$
*is a primitively generated divided power algebra on generators $(kz_{r-j}^{\#}kz_{r-j})^{(j)}\phi_{r-j}^{(j)}e_r(mp^{r-j})$ in $\operatorname{Ext}^{2mp^{r-j}+2p^{r-j}-2}(I^{(r)}, \Gamma^{p^{r-j}(j)})$, $0 \leq m < p^j$, of tri-degree $(2mp^{r-j} + 2p^{r-j} - 2, 1, p^{r-j})$.*

*A similar statement holds by duality for the tri-graded Hopf algebra*
$$\operatorname{Ext}_{\mathcal{P}}^*(S^{*(j)}, S^{*(r)}).$$

(5) *The tri-graded Hopf algebra*
$$\operatorname{Ext}_{\mathcal{P}}^*(\Gamma^{*(j)}, \Gamma^{*(r)})$$

*is a primitively generated divided power algebra on generators* $e_r(mp^{r-j})\phi_{r-j}^{\#(j)}$ *in* $\mathrm{Ext}_{\mathcal{P}}^{2mp^{r-j}}(\Gamma^{p^{r-j}(j)}, I^{(r)})$, $0 \leq m < p^j$, *of tri-degree* $(2mp^{r-j}, p^{r-j}, 1)$.

*A similar statement holds by duality for the tri-graded Hopf algebra*

$$\mathrm{Ext}_{\mathcal{P}}^*(S^{*(r)}, S^{*(j)}).$$

(6) *The tri-graded Hopf algebra*

$$\mathrm{Ext}_{\mathcal{P}}^*(\Lambda^{*(r)}, \Lambda^{*(j)})$$

*is a primitively generated divided power algebra on generators* $kz_{r-j}^{(j)}\phi_{r-j}^{(j)}e_r(mp^{r-j})$ *in* $\mathrm{Ext}_{\mathcal{P}}^{2mp^{r-j}+p^{r-j}-1}(I^{(r)}, \Lambda^{p^{r-j}(j)})$, $0 \leq m < p^j$, *of tri-degree* $(2mp^{r-j} + p^{r-j} - 1, 1, p^{r-j})$.

*A similar statement holds by duality for the tri-graded Hopf algebra*

$$\mathrm{Ext}_{\mathcal{P}}^*(\Lambda^{*(j)}, \Lambda^{*(r)}).$$

*Proof.* The construction of the map in Theorem 4.5 implies that the isomorphism

(5.8.1) $$S^*(V_j) @>\sim>> \mathrm{Ext}_{\mathcal{P}}^*(\Gamma^{*(r)}, S^{*p^{r-j}(j)})$$

is multiplicative. This identifies $\mathrm{Ext}_{\mathcal{P}}^*(\Gamma^{*(r)}, S^{*(j)})$ as a tri-graded algebra. Finally the elements of $V_j$ are primitive by dimension considerations. This proves (1).

The proofs of the remaining assertions follow in a similar, straightforward manner from the computations given earlier in this section. □

We derive from Theorem 5.8 the following corollary which gives a partial answer to a question raised in [F-S, 5.8].

COROLLARY 5.9. *The image of* $e_r \in \mathrm{Ext}_{\mathcal{P}}^{2p^{r-1}}(I^{(r)}, I^{(r)})$ *in* $\mathrm{Ext}_{\mathcal{P}}^{2p^{r-1}}(\Gamma^{p^{r-1}(1)}, S^{p^{r-1}(1)})$ *equals* $e_1^{p^{r-1}}$ *(up to a nonzero scalar factor). Moreover if one modifies the definition of* $e_r$ *the way it was done in* [S-F-B], *then the image of* $e_r$ *is exactly* $e_1^{p^{r-1}}$.

*Proof.* We start by computing the image of $e_r$ in $\mathrm{Ext}_{\mathcal{P}}^{2p^{r-1}}(\Gamma^{p(r-1)}, S^{p(r-1)})$. To do so, note that the homomorphism $\mathrm{Ext}_{\mathcal{P}}^{2p^{r-1}}(I^{(r)}, I^{(r)}) \to \mathrm{Ext}_{\mathcal{P}}^{2p^{r-1}}(\Gamma^{p(r-1)}, S^{p(r-1)})$ is a homogeneous component of a homomorphism of tri-graded Hopf algebras

$$\xi : \mathrm{Ext}_{\mathcal{P}}^*(\widetilde{\Gamma}^{*(r)}, \widetilde{S}^{*(r)}) \to \mathrm{Ext}_{\mathcal{P}}^*(\Gamma^{*(r-1)}, S^{*(r-1)})$$

induced by natural homomorphisms of exponential functors $\Gamma^{*(r-1)} \to \widetilde{\Gamma}^{*(r)}$, $\widetilde{S}^{*(r)} \to S^{*(r-1)}$. Here $\sim$ denotes the $p$-rarefaction of an exponential functor, i.e.
$$\widetilde{A}^n = \begin{cases} A^m & \text{if } n = pm \\ 0 & \text{if } n \not\equiv 0 \mod p \end{cases}.$$

Since the element $e_r \in \operatorname{Ext}_{\mathcal{P}}^{2p^{r-1}}(I^{(r)}, I^{(r)})$ is primitive in $\operatorname{Ext}_{\mathcal{P}}^*(\widetilde{\Gamma}^{*(r)}, \widetilde{S}^{*(r)})$ its image in $\operatorname{Ext}_{\mathcal{P}}^*(\Gamma^{*(r-1)}, S^{*(r-1)}) = S^*(\operatorname{Ext}_{\mathcal{P}}^*(I^{(r-1)}, I^{(r-1)}))$ also has to be primitive. The subspace of primitive elements in a (primitively generated) symmetric algebra $S^*(V)$ coincides with a linear span $<V, V^p, V^{p^2}, \ldots>$. From this we derive immediately that the only primitive element of tri-degree $(2p^{r-1}, p, p)$ in $S^*(\operatorname{Ext}_{\mathcal{P}}^*(I^{(r-1)}, I^{(r-1)}))$ is $e_{r-1}^p$. Thus $\xi(e_r)$ is a scalar multiple of $e_{r-1}^p$. Since $\xi$ is multiplicative we conclude further that $\xi(e_r^n)$ is a scalar multiple of $e_{r-1}^{pn}$ for each $n$. Repeating this construction $r-1$ times we finally conclude that the image of $e_r$ in $\operatorname{Ext}_{\mathcal{P}}^{2p^{r-1}}(\Gamma^{p^{r-1}}(1), S^{p^{r-1}}(1))$ is a scalar multiple of $e_1^{p^{r-1}}$.

We proceed to show now that with the modification made in [S-F-B, 3.4] the image of $e_r$ in $\operatorname{Ext}_{\mathcal{P}}^{2p^{r-1}}(\Gamma^{p(r-1)}, S^{p(r-1)})$ is exactly $e_{r-1}^p$ (i.e. the corresponding scalar factor is the identity). Recall that (with the modification made in [S-F-B]) for any $n \geq 2$ the image of $e_r$ under the natural homomorphism
$$\operatorname{Ext}_{\mathcal{P}}^{2p^{r-1}}(I^{(r)}, I^{(r)}) \to \operatorname{Ext}_{\operatorname{GL}_n}^{2p^{r-1}}(k^{n(r)}, k^{n(r)}) \to \operatorname{Ext}_{\operatorname{GL}_{n(1)}}^{2p^{r-1}}(k^{n(r)}, k^{n(r)})$$
$$= H^{2pr-1}(\operatorname{GL}_{n(1)}, gl_n^{(r)}) = \operatorname{Hom}_k(gl_n^{(r)\#}, H^{2p^{r-1}}(\operatorname{GL}_{n(1)}, k)).$$
is nontrivial [F-S, §6] and coincides with the composition
$$gl_n^{(r)\#} \to S^{p^{r-1}}(gl_n^{(1)\#}) \to H^{2p^{r-1}}(\operatorname{GL}_{n(1)}, k)$$
where the first arrow is the standard embedding and the second one is the edge homomorphism in the May spectral sequence. Consider further the natural homomorphism
$$\operatorname{Ext}_{\mathcal{P}}^{2p^{r-1}}(\Gamma^{p(r-1)}, S^{p(r-1)}) \to H^{2p^{r-1}}(\operatorname{GL}_{n(1)}, (\Gamma^{p(r-1)}(k^n))^\# \otimes S^{p(r-1)}(k^n))$$
$$= \operatorname{Hom}_k(\Gamma^{p(r-1)}(k^n) \otimes \Gamma^{p(r-1)}(k^{n\#}), H^{2p^{r-1}}(\operatorname{GL}_{n(1)}, k)).$$
The image of $\xi(e_r)$ under this map coincides with the composite
$$\Gamma^{p(r-1)}(k^n) \otimes \Gamma^{p(r-1)}(k^{n\#}) \twoheadrightarrow k^{n(r)} \otimes k^{n(r)\#}$$
$$= gl_n^{(r)\#} \to S^{p^{r-1}}(gl_n^{(1)\#}) \to H^{2p^{r-1}}(\operatorname{GL}_{n(1)}, k)$$

and, in particular, is nontrivial. Keeping in mind that the tensor product operation
$$\operatorname{Ext}_{\operatorname{GL}_{n(1)}}^*(k, k) \otimes \operatorname{Ext}_{\operatorname{GL}_{n(1)}}^*(k, k) \to \operatorname{Ext}_{\operatorname{GL}_{n(1)}}^*(k, k)$$

coincides with the usual product operation in the cohomology ring $\mathrm{Ext}^*_{\mathrm{GL}_{n(1)}}(k,k) = H^*(\mathrm{GL}_{n(1)}, k)$, using the description of the product operation in $\mathrm{Ext}^*_{\mathcal{P}}(\Gamma^{*(r-1)}, S^{*(r-1)})$ in terms of the tensor product operation and the multiplicative properties of the May homomorphism $S^*(gl_n^{(1)}) \to H^{2*}(\mathrm{GL}_{n(1)}, k)$, we easily conclude that the image of $e^p_{r-1} \in \mathrm{Ext}^{2p^{r-1}}_{\mathcal{P}}(\Gamma^{p(r-1)}, S^{p(r-1)})$ in

$$\mathrm{Hom}_k(\Gamma^{p(r-1)}(k^n) \otimes \Gamma^{p(r-1)}(k^{n\#}), H^{2p^{r-1}}(\mathrm{GL}_{n(1)}, k)$$

is given by the composition

$$\Gamma^{p(r-1)}(k^n) \otimes \Gamma^{p(r-1)}(k^{n\#}) \hookrightarrow (k^{n(r-1)})^{\otimes p} \otimes (k^{n(r-1)\#})^{\otimes p}$$
$$= (gl_n^{(r-1)})^{\otimes p} \to (S^{p^{r-2}}(gl_n^{(1)}))^{\otimes p} \twoheadrightarrow S^{p^{r-1}}(gl_n^{(1)}) \to H^{2p^{r-1}}(\mathrm{GL}_{n(1)}, k).$$

A straightforward verification shows that the resulting maps $\Gamma^{p(r-1)}(k^n) \otimes \Gamma^{p(r-1)}(k^{n\#}) \to H^{2p^{r-1}}(\mathrm{GL}_{n(1)}, k)$ coincide. Since this map is moreover nontrivial we conclude that the constant relating $\xi(e_r)$ and $e^p_{r-1}$ is equal to the identity. $\square$

To simplify the notation, we adopt the following convention: for homogeneous strict polynomial functors $A, B$ in $\mathcal{P}$ we define

$$\mathrm{Ext}^*_{\mathcal{P} \to \mathcal{F}}(A, B) = \varinjlim_m \mathrm{Ext}^*_{\mathcal{P}}(A^{(m)}, B^{(m)}).$$

Letting $r$ vary over positive integers and $j$ vary between 0 and $r$, and setting $h = r - j$, we have the following consequence of Theorem 5.8.

COROLLARY 5.10. *Let $h$ be a nonnegative integer.*

(1) *The tri-graded Hopf algebra*

$$\mathrm{Ext}^*_{\mathcal{P} \to \mathcal{F}}(\Gamma^*, S^{*(h)})$$

*is a primitively generated polynomial algebra on generators $e(mp^h)\phi_h^{\#}$ in $\mathrm{Ext}^{2p^h m}_{\mathcal{P} \to \mathcal{F}}(\Gamma^{p^h}, I^{(h)})$, $0 \leq m$, of tri-degree $(2p^h m, p^h, 1)$.*

*Similarly, the tri-graded Hopf algebra*

$$\mathrm{Ext}^*_{\mathcal{P} \to \mathcal{F}}(\Gamma^{*(h)}, S^*)$$

*is a primitively generated polynomial algebra on generators $\phi_h e(mp^h)$ in $\mathrm{Ext}^{2p^h m}_{\mathcal{P} \to \mathcal{F}}(I^{(h)}, S^{p^h})$, $0 \leq m$, of tri-degree $(2p^h m, 1, p^h)$.*

(2) *The tri-graded Hopf algebra*

$$\mathrm{Ext}^*_{\mathcal{P} \to \mathcal{F}}(\Gamma^*, \Lambda^{*(h)})$$

is a primitively generated exterior algebra on generators $e(mp^h)\phi_h^{\#}$ in $\text{Ext}_{\mathcal{P}\to\mathcal{F}}^{2p^h m}(\Gamma^{p^h}, I^{(h)})$, $0 \leq m$, of tri-degree $(2p^h m, p^h, 1)$.

A similar statement holds by duality for the tri-graded Hopf algebra

$$\text{Ext}_{\mathcal{P}\to\mathcal{F}}^*(\Lambda^{*(h)}, S^*).$$

(3) The tri-graded Hopf algebra

$$\text{Ext}_{\mathcal{P}\to\mathcal{F}}^*(\Gamma^{*(h)}, \Lambda^*)$$

is a primitively generated exterior algebra on generators $kz_h\phi_h e(mp^h)$ in $\text{Ext}_{\mathcal{P}\to\mathcal{F}}^{2mp^h + p^h - 1}(I^{(h)}, \Lambda^{p^h})$, $0 \leq m$, of tri-degree $(2mp^h + p^h - 1, 1, p^h)$.

A similar statement holds by duality for the tri-graded Hopf algebra

$$\text{Ext}_{\mathcal{P}\to\mathcal{F}}^*(\Lambda^*, S^{*(h)}).$$

(4) The tri-graded Hopf algebra

$$\text{Ext}_{\mathcal{P}\to\mathcal{F}}^*(\Gamma^{*(h)}, \Gamma^*)$$

is a primitively generated divided power algebra on generators $kz_h^{\#}kz_h\phi_h e_r(mp^h)$ in $\text{Ext}_{\mathcal{P}\to\mathcal{F}}^{2mp^h + 2p^h - 2}(I^{(h)}, \Gamma^{p^h})$, $0 \leq m$, of tri-degree $(2mp^h + 2p^h - 2, 1, p^h)$.

A similar statement holds by duality for the tri-graded Hopf algebra

$$\text{Ext}_{\mathcal{P}\to\mathcal{F}}^*(S^*, S^{*(h)}).$$

(5) The tri-graded Hopf algebra

$$\text{Ext}_{\mathcal{P}\to\mathcal{F}}^*(\Gamma^*, \Gamma^{*(h)})$$

is a primitively generated divided power algebra on generators $e(mp^h)\phi_h^{\#}$ in $\text{Ext}_{\mathcal{P}\to\mathcal{F}}^{2mp^h}(\Gamma^{p^h}, I^{(h)})$, $0 \leq m$, of tri-degree $(2mp^h, p^h, 1)$.

A similar statement holds by duality for the tri-graded Hopf algebra

$$\text{Ext}_{\mathcal{P}\to\mathcal{F}}^*(S^{*(h)}, S^*).$$

(6) The tri-graded Hopf algebra

$$\text{Ext}_{\mathcal{P}\to\mathcal{F}}^*(\Lambda^{*(h)}, \Lambda^*)$$

is a primitively generated divided power algebra on generators $kz_h\phi_h e(mp^h)$ in $\text{Ext}_{\mathcal{P}\to\mathcal{F}}^{2mp^h + p^h - 1}(I^{(h)}, \Lambda^{p^h})$, $0 \leq m$, of tri-degree $(2mp^h + p^h - 1, 1, p^h)$. A similar statement holds by duality for the tri-graded Hopf algebra

$$\text{Ext}_{\mathcal{P}\to\mathcal{F}}^*(\Lambda^*, \Lambda^{*(h)}).$$

## 6. Ext-groups between the classical functors in the category $\mathcal{F}$

As we see in this section, the extension-of-scalars formula (Theorem 3.4) and the comparison theorem (Theorem 3.10) allow us in many cases to reduce the computation of Ext-groups in the category $\mathcal{F}$ to the computation of appropriate Ext-groups in the category $\mathcal{P}$.

Throughout this section $k$ is a finite field with $q = p^N$ elements, $\mathcal{F} = \mathcal{F}(k)$ is the category of functors $\mathcal{V}_k^f \to \mathcal{V}_k$ and $\mathcal{P} = \mathcal{P}(k)$ is the category of strict polynomial functors. Finally, $A^*$ and $B^*$ are exponential (homogeneous) strict polynomial functors with $\deg(A^i) = \deg(B^i) = i$.

For nonnegative integers $h < N$, and $j$, $l$, $n$, let $Pr_n^h(j,l)$ be the set of all sequences $j_0, j_1, \ldots, j_{n-1}, l_1, \ldots, l_n$ of nonnegative integers with

$$j = j_0 + \cdots + j_{n-1} + l_1 p^{N-h} + \cdots + l_n p^{nN-h},$$
$$l = j_0 p^h + \cdots + j_{n-1} p^{(n-1)N+h} + l_1 + \cdots + l_n.$$

The obvious embedding $Pr_n^h(j,l) \hookrightarrow Pr_{n+1}^h(j,l)$ (taking $j_n = l_{n+1} = 0$) is an isomorphism, provided that $p^{(n+1)N-h} > j$ and $p^{nN+h} > l$. We denote by $Pr^h(j,l)$ the limit set $\varinjlim_n Pr_n^h(j,l)$.

For any sequence $j_0, \ldots, j_{n-1}, l_1, \ldots, l_n$ in $Pr_n^h(j,l)$ consider the following homomorphism

$$\bigotimes_{s=0}^{n-1} \mathrm{Ext}^*_{\mathcal{P} \to \mathcal{F}}(A^{j_s(sN+h)}, B^{j_s p^{sN+h}}) \otimes \bigotimes_{t=1}^{n} \mathrm{Ext}^*_{\mathcal{P} \to \mathcal{F}}(A^{l_t p^{tN-h}(h)}, B^{l_t(tN)})$$

$$\to \bigotimes_{s=0}^{n-1} \mathrm{Ext}^*_{\mathcal{F}}(A^{j_s(sN+h)}, B^{j_s p^{sN+h}}) \otimes \bigotimes_{t=1}^{n} \mathrm{Ext}^*_{\mathcal{F}}(A^{l_t p^{tN-h}(h)}, B^{l_t(tN)})$$

$$= \bigotimes_{s=0}^{n-1} \mathrm{Ext}^*_{\mathcal{F}}(A^{j_s(h)}, B^{j_s p^{sN+h}}) \otimes \bigotimes_{t=1}^{n} \mathrm{Ext}^*_{\mathcal{F}}(A^{l_t p^{tN-h}(h)}, B^{l_t})$$

$$\to \mathrm{Ext}^*_{\mathcal{F}}(A^{j(h)}, B^l).$$

Our next result determines $\mathrm{Ext}^*_{\mathcal{F}}(A^*, B^*)$ in terms of $\mathrm{Ext}^*_{\mathcal{P} \to \mathcal{F}}(A^*, B^*)$.

THEOREM 6.1. *Fix a nonnegative integer $h$ less than $N$. The resulting homomorphism*

$$\bigoplus_{Pr^h(j,l)} \bigotimes_{s \geq 0} \mathrm{Ext}^*_{\mathcal{P} \to \mathcal{F}}(A^{j_s(sN+h)}, B^{j_s p^{sN+h}}) \otimes \bigotimes_{t \geq 1} \mathrm{Ext}^*_{\mathcal{P} \to \mathcal{F}}(A^{l_t p^{tN-h}(h)}, B^{l_t(tN)})$$

$$\to \mathrm{Ext}^*_{\mathcal{F}}(A^{j(h)}, B^l)$$

*is an isomorphism.*

*Proof.* Choose an integer $n$ so that $p^{nN-h} > j$, $p^{nN} > l$ and consider a tower of finite field extensions $L \supset K \supset k$ in which each stage is of degree $n$. Theorem 3.4 gives a natural isomorphism

$$K \otimes_k \mathrm{Ext}^*_{\mathcal{F}(k)}(A^{j(h)}, B^l) = \mathrm{Ext}^*_{\mathcal{F}(K)}(A_K^{j(h)} \circ (t \circ \tau), B_K^l).$$

Moreover, according to Lemma 3.6, the functor $t \circ \tau$ may be identified with $I \oplus I^{(N)} \oplus \cdots \oplus I^{((n-1)N)}$. The exponential property of $A^*$ gives now a natural isomorphism

$$A_K^{j(h)} \circ (t \circ \tau) = \bigoplus_{j_0+\cdots+j_{n-1}=j} A_K^{j_0(h)} \otimes \cdots \otimes A_K^{j_{n-1}((n-1)N+h)}.$$

Extending scalars further from $K$ to $L$ and applying now the previous procedure to the second variable, we get the following formula:

$$L \otimes_k \mathrm{Ext}^*_{\mathcal{F}(k)}(A^{j(h)}, B^l) = \bigoplus_{\substack{j_0+\cdots+j_{n-1}=j \\ l_0+\cdots+l_{n-1}=l}} \mathrm{Ext}^*_{\mathcal{F}(L)}(\bigotimes_{s=0}^{n-1} A_L^{j_s(sN+h)}, \bigotimes_{t=0}^{n-1} B_L^{l_t(tnN)}).$$

In view of the exponential properties of $A^*$ and $B^*$, each summand in the above formula may be further identified (cf. Corollary 1.8) with

$$\bigoplus_{\substack{j_{0,0}+\cdots+j_{0,n-1}=j_0 \\ \cdots \\ j_{n-1,0}+\cdots+j_{n-1,n-1}=j_{n-1}}} \bigoplus_{\substack{l_{0,0}+\cdots+l_{0,n-1}=l_0 \\ \cdots \\ l_{n-1,0}+\cdots+l_{n-1,n-1}=l_{n-1}}} \bigotimes_{s,t=0}^{n-1} \mathrm{Ext}^*_{\mathcal{F}(L)}(A_L^{j_{s,t}(sN+h)}, B_L^{l_{t,s}(tnN)}).$$

Note further that

$$\deg(A_L^{j_{s,t}(sN+h)}) = j_{s,t} \cdot p^{sN+h} \leq j \cdot p^{(n-1)N+h} = \frac{j \cdot \mathrm{card}(L)}{p^{(n^2-n+1)N-h}} < \mathrm{card}(L)$$

and also

$$\deg(B_L^{l_{t,s}(tnN)}) = l_{t,s} p^{tnN} \leq l \cdot p^{(n-1)nN} = \frac{l \cdot \mathrm{card}(L)}{p^{nN}} < \mathrm{card}(L).$$

Theorem 3.10 and Lemma 3.11 now show that we may replace $\mathrm{Ext}^*_{\mathcal{F}}$ by $\mathrm{Ext}^*_{\mathcal{P} \to \mathcal{F}}$ everywhere in the above formula. This implies, in particular, that the summand corresponding to $\{j_{t,s}, l_{t,s}\}_{t,s}$ is trivial unless the corresponding degrees coincide, i.e. unless we have the following relations: $j_{s,t} p^{sN+h} = l_{t,s} p^{tnN}$ for all $s$, $t$. It is easy to see from the choice of $n$ that the above relation implies that $j_{s,t} = l_{t,s} = 0$ unless $t = 0$ or $t = 1$. The remaining $j_{t,s}$ and $l_{s,t}$ are related by equations

$$l_{0,s} = j_{s,0} p^{sN+h}, \quad j_{s,1} = l_{1,s} p^{(n-s)N-h}.$$

Since $\sum_{s,t} j_{s,t} = j$, $\sum_{s,t} l_{s,t} = l$, we conclude that $(j^0, \ldots, j^{n-1}, l^1, \ldots, l^n) = (j_{0,0}, \ldots, j_{n-1,0}, l_{1,n-1}, \ldots, l_{1,0})$ is an element of $Pr_n^h(j,l)$, i.e. we are left with the direct sum over $Pr_n^h(j,l)$:

$$L \otimes_k \mathrm{Ext}^*_{\mathcal{F}(k)}(A^{j(h)}, B^l) = \bigoplus_{Pr_n^h(j,l)} L \otimes_k \bigotimes_{s=0}^{n-1} \mathrm{Ext}^*_{\mathcal{P} \to \mathcal{F}}(A^{j^s(sN+h)}, B^{k^s p^{sN+h}})$$
$$\otimes \bigotimes_{t=1}^{n} \mathrm{Ext}^*_{\mathcal{P} \to \mathcal{F}}(A^{l^t p^{tN-h}(h)}, B^{l^t(tN)}).$$

Finally one checks without difficulty that the composite of the above isomorphism with the extension of scalars in the homomorphism we started with is the identity map. □

It is pleasing to reformulate the result of Theorem 6.1 in terms of the corresponding Hopf algebras. Observe that for any $s \geq 0$, $t \geq 1$ we have natural homomorphisms of tri-graded vector spaces

$$\mathrm{Ext}^*_{\mathcal{P} \to \mathcal{F}}(A^{*(sN+h)}, B^*) \to \mathrm{Ext}^*_{\mathcal{F}}(A^{*(sN+h)}, B^*) = \mathrm{Ext}^*_{\mathcal{F}}(A^{*(h)}, B^*)$$
$$\mathrm{Ext}^*_{\mathcal{P} \to \mathcal{F}}(A^*, B^{*(tN-h)}) = \mathrm{Ext}^*_{\mathcal{P} \to \mathcal{F}}(A^{*(h)}, B^{*(tN)}) \to \mathrm{Ext}^*_{\mathcal{F}}(A^{*(h)}, B^{*(tN)})$$
$$= \mathrm{Ext}^*_{\mathcal{F}}(A^{*(h)}, B^*).$$

Furthermore, it is easy to see that the tri-graded space $\mathrm{Ext}^*_{\mathcal{P} \to \mathcal{F}}(A^{*(sN+h)}, B^*)$ is $p^{sN+h}$-connected with respect to the total degree, whereas $\mathrm{Ext}^*_{\mathcal{P} \to \mathcal{F}}(A^*, B^{*(tN-h)})$ is $p^{tN-h}$-connected with respect to the total degree, which allows us to consider the infinite tensor product

$$\bigotimes_{s \geq 0} \mathrm{Ext}^*_{\mathcal{P} \to \mathcal{F}}(A^{*(sN+h)}, B^*) \otimes \bigotimes_{t \geq 1} \mathrm{Ext}^*_{\mathcal{P} \to \mathcal{F}}(A^*, B^{*(tN-h)}).$$

Suppose now that $A^*$ and $B^*$ are Hopf functors. The product operations in $\mathrm{Ext}^*_{\mathcal{F}}(A^{*(h)}, B^*)$ define a natural homomorphism from the above infinite tensor product to $\mathrm{Ext}^*_{\mathcal{F}}(A^{*(h)}, B^*)$. Moreover, assuming that $A^*$ (resp. $B^*$) is $\varepsilon(A)$-commutative (resp. $\varepsilon(B)$-commutative), the above map is even a homomorphism of Hopf algebras, provided that one endows the tensor product algebra with multiplication and comultiplication taking into account the sign convention given by Lemma 1.11.

COROLLARY 6.2. *Let $A^*$ and $B^*$ be commutative Hopf functors. The natural homomorphism*

$$\bigotimes_{s \geq 0} \mathrm{Ext}^*_{\mathcal{P} \to \mathcal{F}}(A^{*(sN+h)}, B^*) \otimes \bigotimes_{t \geq 1} \mathrm{Ext}^*_{\mathcal{P} \to \mathcal{F}}(A^*, B^{*(tN-h)}) \to \mathrm{Ext}^*_{\mathcal{F}}(A^{*(h)}, B^*)$$

*is an isomorphism of tri-graded Hopf algebras.*

Combining Corollary 6.2 with Corollary 5.9 we are able to compute various important Ext-algebras in the category $\mathcal{F}$.

THEOREM 6.3. *Let h be a nonnegative integer.*

(1) *The tri-graded Hopf algebra*

$$\operatorname{Ext}_{\mathcal{F}}^{*}(\Gamma^{*(h)}, S^{*})$$

*is a primitively generated polynomial algebra on generators $\phi_{h+sN}e(mp^{h+sN})$ in $\operatorname{Ext}_{\mathcal{F}}^{2p^{h+sN}m}(I^{(h)}, S^{p^{h+sN}})$, $0 \leq m$, $0 \leq s$, of tri-degree $(2p^{h+sN}m, 1, p^{h+sN})$ and generators $e(mp^{tN-h})\phi_{tN-h}^{\#}$ in $\operatorname{Ext}_{\mathcal{F}}^{2p^{tN-h}m}(\Gamma^{p^{tN-h}(h)}, I)$, $0 \leq m$, $1 \leq t$, of tri-degree $(2p^{tN-h}m, p^{tN-h}, 1)$.*

(2) *The tri-graded Hopf algebra*

$$\operatorname{Ext}_{\mathcal{F}}^{*}(\Gamma^{*(h)}, \Lambda^{*})$$

*is a primitively generated exterior algebra on generators $kz_{h+sN}\phi_{h+sN}e(mp^{h+sN})$ in $\operatorname{Ext}_{\mathcal{F}}^{2mp^{h+sN}+p^{h+sN}-1}(I^{(h)}, \Lambda^{p^{h+sr}})$, $0 \leq m$, $0 \leq s$, of tri-degree $(2mp^{h+sN} + p^{h+sN} - 1, 1, p^{h+sN})$ and generators $e(mp^{tN-h})\phi_{tN-h}^{\#(h)}$ in $\operatorname{Ext}_{\mathcal{F}}^{2p^{tN-h}m}(\Gamma^{p^{tN-h}(h)}, I)$, $0 \leq m$, $0 \leq s$, of tri-degree $(2p^{tN-h}m, p^{tN-h}, 1)$.*

(3) *The tri-graded Hopf algebra*

$$\operatorname{Ext}_{\mathcal{F}}^{*}(\Gamma^{*(h)}, \Gamma^{*})$$

*is a primitively generated divided power algebra on generators $kz_{h+sN}^{\#}kz_{h+sN}\phi_{h+sN}e(mp^{h+sN})$ in $\operatorname{Ext}_{\mathcal{F}}^{2mp^{sN+h}+2p^{sN+h}-2}(I^{(h)}, \Gamma^{p^{sN+h}})$, $0 \leq m$, $0 \leq s$, of tri-degree $(2mp^{h+sN} + 2p^{h+sN} - 2, 1, p^{h+sN})$ and generators $e(mp^{tN-h})\phi_{tN-h}^{\#(h)}$ in $\operatorname{Ext}_{\mathcal{F}}^{2mp^{tN-h}}(\Gamma^{p^{tN-h}(h)}, I)$, $0 \leq m$, $1 \leq t$, of tri-degree $(2mp^{tN-h}, p^{tN-h}, 1)$.*

*A similar statement holds by duality for the tri-graded Hopf algebra*

$$\operatorname{Ext}_{\mathcal{F}}^{*}(S^{*(N-h)}, S^{*}).$$

(4) *The tri-graded Hopf algebra*

$$\operatorname{Ext}_{\mathcal{F}}^{*}(\Lambda^{*(h)}, \Lambda^{*})$$

*is a primitively generated divided power algebra on generators $kz_{h+sN}\phi_{h+sN}e(mp^{sN+h})$ in $\operatorname{Ext}_{\mathcal{F}}^{2mp^{sN+h}+p^{sN+h}-1}(I^{(h)}, \Lambda^{p^{sN+h}})$, $0 \leq m$, $0 \leq s$, of tri-degree $(2mp^{sN+h} + p^{sN+h} - 1, 1, p^{sN+h})$ and generators $e(mp^{tN-h})\phi_{tN-h}^{\#(h)}kz_{tN-h}^{\#(h)}$ in $\operatorname{Ext}_{\mathcal{F}}^{2mp^{tN-h}+p^{tN-h}-1}(\Lambda^{p^{tN-h}(h)}, I)$, $0 \leq m$, $1 \leq t$, of tri-degree $(2mp^{tN-h} + p^{tN-h} - 1, p^{tN-h}, 1)$.*

## Appendix: GL-cohomology and $\mathrm{Ext}_\mathcal{F}$-groups
## By A. Suslin

The main subject discussed in this appendix is the map on Ext-groups

$$\mathrm{Ext}^*_\mathcal{F}(P,Q) \to \mathrm{Ext}^*_{k[\mathrm{GL}_n(k)]}(P(k^n), Q(k^n)) = H^*(\mathrm{GL}_n(k), P(k^n)^\# \otimes Q(k^n))$$

induced by the exact functor $\mathcal{F}(k) \to \{k[\mathrm{GL}_n(k)] - \text{modules}\}: \quad P \mapsto P(k^n)$ (here $k[\mathrm{GL}_n(k)]$ stands for the group ring of the discrete group $\mathrm{GL}_n(k)$ over the field $k$). All $\mathrm{GL}_n(k)$-modules considered below are always assumed to be $k$-vector spaces on which the group $\mathrm{GL}_n(k)$ acts by $k$-linear transformations (i.e. $k[\mathrm{GL}_n(k)]$-modules); we abbreviate the notation $\mathrm{Ext}^*_{k[\mathrm{GL}_n(k)]}$ to $\mathrm{Ext}^*_{\mathrm{GL}_n(k)}$. Note that the $\mathrm{Ext}^*_{\mathrm{GL}_n(k)}(P(k^n), Q(k^n))$-groups above stabilize with the growth of $n$, according to a well-known theorem of W. Dwyer [D], provided that the functors $P$ and $Q$ are finite. We use the notation $\mathrm{Ext}^*_{GL(k)}(P,Q)$ for the stable values of the corresponding Ext-groups. The main result proved below is the following theorem:

THEOREM A.1. *Let $k$ be a finite field. Then the natural map*

$$\mathrm{Ext}^*_\mathcal{F}(P,Q) \to \mathrm{Ext}^*_{\mathrm{GL}(k)}(P,Q)$$

*is an isomorphism for any finite functors $P, Q \in \mathcal{F} = \mathcal{F}(k)$.*

Over more general fields one has to compare $\mathrm{Ext}^*_\mathcal{F}$-groups with the so-called stable K-theory (see [B-P]). The above theorem shows that the conjecture stated in [B-P] (saying that for finite functors stable K-theory equals topological Hochschild homology) holds at least over finite fields. Apparently an appropriate modification of the argument given below should work over any field of positive characteristic. I did not try to work out this more general case, the more so that $\mathrm{Ext}_\mathcal{F}$-groups do not seem to be computable unless we are dealing with finite fields.

Theorem A.1 was known to be true (in a more general situation) in case $P = Q = I$ by a work of Dundas and McCarthy [D-M]. For finite fields a purely algebraic proof of this result was given in [F-S]. In a recent paper [B2], S. Betley gives a different proof of Theorem A.1, following the approach developed in [F-S].

Throughout this appendix $k = \mathbb{F}_q$ is a finite field of characteristic $p$ (with $q$ elements).

LEMMA A.2. *For any finite functors $A, B \in \mathcal{F}$ there is a natural duality isomorphism $\mathrm{Ext}^*_{GL(k)}(A,B) = \mathrm{Ext}^*_{\mathrm{GL}(k)}(B^\#, A^\#)$ which makes the following*

*diagram commute*

$$\begin{CD} \operatorname{Ext}^*_{\mathcal{F}}(A,B) @>\sim>> \operatorname{Ext}^*_{\mathcal{F}}(B^{\#},A^{\#}) \\ @VVV @VVV \\ \operatorname{Ext}^*_{\operatorname{GL}(k)}(A,B) @>\sim>> \operatorname{Ext}^*_{\operatorname{GL}(k)}(B^{\#},A^{\#}). \end{CD}$$

*Proof.* For any $s \geq 0$ we can find an integer $N(s)$ such that for all $n \geq N(s)$ we have canonical identifications

$$\operatorname{Ext}^s_{\operatorname{GL}(k)}(A,B) = \operatorname{Ext}^s_{\operatorname{GL}_n(k)}(A(k^n),B(k^n))$$
$$\operatorname{Ext}^s_{\operatorname{GL}(k)}(B^{\#},A^{\#}) = \operatorname{Ext}^s_{\operatorname{GL}_n(k)}(B^{\#}(k^n),A^{\#}(k^n))$$
$$= \operatorname{Ext}^s_{\operatorname{GL}_n(k)}(B(k^{n\#})^{\#},A(k^{n\#})^{\#}).$$

Furthermore for any finite dimensional $\operatorname{GL}_n(k)$-modules $M,N$ the exact contravariant functor $M \mapsto M^{\#}$ defines natural duality isomorphisms (cf. Lemma 1.12)

$$\operatorname{Ext}^*_{\operatorname{GL}_n(k)}(M,N) @>\sim>> \operatorname{Ext}^*_{\operatorname{GL}_n(k)}(N^{\#},M^{\#}).$$

Taking here $M = A(k^n), N = B(k^n)$, we identify $\operatorname{Ext}^s_{\operatorname{GL}_n(k)}(A(k^n),B(k^n))$ with $\operatorname{Ext}^s_{\operatorname{GL}_n(k)}(B(k^n)^{\#},A(k^n)^{\#}) = \operatorname{Ext}^s_{\operatorname{GL}_n(k)}(B^{\#}(k^{n\#}),A^{\#}(k^{n\#}))$. Finally for any $\operatorname{GL}_n(k)$-module $M$ consider a new $\operatorname{GL}_n(k)$-module $\widetilde{M}$, which coincides with $M$ as a $k$-vector space. The action of $\operatorname{GL}_n(k)$ on $\widetilde{M}$ is obtained from the original one by means of the automorphism

$$\operatorname{GL}_n(k) @>\sim>> \operatorname{GL}_n(k): \quad \alpha \mapsto (\alpha^{\mathrm{T}})^{-1}.$$

The functor $M \mapsto \widetilde{M}$ is obviously an equivalence of the category of $\operatorname{GL}_n(k)$-modules with itself and hence defines natural isomorphisms

$$\operatorname{Ext}^*_{\operatorname{GL}_n(k)}(M,N) @>\sim>> \operatorname{Ext}^*_{\operatorname{GL}_n(k)}(\widetilde{M},\widetilde{N}).$$

To finish the proof it suffices to note now that (for any functor $A \in \mathcal{F}$ and any $n$) the $\operatorname{GL}_n(k)$-modules $A(\widetilde{k^n})$ and $A(k^{n\#})$ are canonically isomorphic. $\square$

LEMMA A.3. *Let $V \in \mathcal{V}^f$ be a finite dimensional vector space over $k$. Assume that $n \geq \dim V$. Then for any $\operatorname{GL}_n(k)$-module $M$ we have isomorphisms*:

$$\operatorname{Ext}^*_{\operatorname{GL}_n(k)}(P_V(k^n),M) = \bigoplus_{W \subset V} H^*(\begin{pmatrix} 1_{n-\dim W} & * \\ 0 & \operatorname{GL}_{\dim W}(k) \end{pmatrix}, M).$$

*Here $P_V$ is the projective generator of $\mathcal{F}$, discussed in Section 1 and $W \subset V$ on the right runs through all subspaces of the vector space $V$.*

*Proof.* $P_V(k^n) = k[\mathrm{Hom}(V, k^n)]$ is a permutation $\mathrm{GL}_n(k)$-module. Two homomorphisms $f, g : V \to k^n$ are in the same $\mathrm{GL}_n(k)$-orbit if and only if $\mathrm{Ker}\, f = \mathrm{Ker}\, g$. Thus orbits of the action of $\mathrm{GL}_n(k)$ on $\mathrm{Hom}(V, k^n)$ are in one-to-one correspondence with subspaces $W \subset V$. For any $W \subset V$ choose a representative $f_W$ of the corresponding orbit so that $\mathrm{Im}\, f_W = k^{n-\dim W}$, the standard coordinate subspace in $k^n$. It is clear that the stabilizer of $f_W$ coincides with the affine group

$$\begin{pmatrix} 1_{n-\dim W} & * \\ 0 & \mathrm{GL}_{\dim W}(k) \end{pmatrix}.$$

Thus our statement follows from Shapiro's Lemma. □

PROPOSITION A.4. *Let $B \in \mathcal{F}$ be a finite functor. For any $s \geq 0$ there exists $N(s)$ such that for $n \geq N(s)$*

$$H^i(\mathrm{GL}_n(k), B(k^n)) = \begin{cases} B(0) & \text{for } i = 0 \\ 0 & \text{for } 0 < i \leq s. \end{cases}$$

*Proof.* This follows immediately from the theorem of Betley [B1], stating that over any field $k$ cohomology of $\mathrm{GL}_n(k)$ with coefficients in $B(k^n)$ coincides (for $n$ big enough with respect to the cohomology index) with cohomology of $\mathrm{GL}_n(k)$ with coefficients in the trivial submodule $B(0)$, and also from the theorem of Quillen [Q], giving the vanishing of higher dimensional cohomology in case of finite fields. □

Let $A \in \mathcal{F}$ be any functor and let $j$ be a nonnegative integer. For any $V \in \mathcal{V}$ set

$$H^i_{A,j}(V) = H^i(\mathrm{Hom}(V, k^j), A(k^j \oplus V)),$$

where the action of $\mathrm{Hom}(V, k^j)$ on $A(k^j \oplus V)$ is defined via the natural embedding

$$\mathrm{Hom}(V, k^j) = \begin{pmatrix} 1_j & \mathrm{Hom}(V, k^j) \\ 0 & 1_V \end{pmatrix} \hookrightarrow \mathrm{GL}(k^j \oplus V).$$

Each $k$-space homomorphism $V \to V'$ defines a homomorphism of (abelian) groups $\mathrm{Hom}(V', k^j) \to \mathrm{Hom}(V, k^j)$ and a homomorphism of $\mathrm{Hom}(V', k^j)$-modules $A(k^j \oplus V) \to A(k^j \oplus V')$, thus defining a homomorphism in cohomology $H^i_{A,j}(V) \to H^i_{A,j}(V')$. One checks immediately that in this way $H^i_{A,j}$ becomes a functor from $\mathcal{V}^f$ to $\mathcal{V}$, i.e. an element of $\mathcal{F}$.

LEMMA A.5. *Assume that the functor $A$ is finite. Then all functors $H^i_{A,j}$ are finite as well. Moreover, $H^i_{A,j}(0) = 0$ for $i > 0$ and $H^0_{A,j}(0) = A(k^j)$.*

*Proof.* The second part of the statement is obvious. To prove the first one we consider first a special case $A = S^{m_1} \otimes \cdots \otimes S^{m_l}$. To simplify notation we use the following multi-index notation: for every multi-index $i_\bullet = (i_1, \ldots, i_l)$ set

$$S^{i_\bullet} = S^{i_1} \otimes \cdots \otimes S^{i_l}.$$

We consider the usual partial ordering on the set of multi-indices, and for each multi-index $i_\bullet$ we set $|i_\bullet| = i_1 + \cdots + i_l$. Note that the functor $S^{m_\bullet}$ is "exponential"; more precisely we have the following obvious formula:

$$S^{m_\bullet}(k^j \oplus V) = \bigoplus_{0 \leq i_\bullet \leq m_\bullet} S^{m_\bullet - i_\bullet}(k^j) \otimes S^{i_\bullet}(V).$$

Consider an increasing filtration on $S^{m_\bullet}(k^j \oplus V)$ defined by the formula

$$\Phi_t = \bigoplus_{\substack{0 \leq i_\bullet \leq m_\bullet \\ |i_\bullet| \leq t}} S^{m_\bullet - i_\bullet}(k^j) \otimes S^{i_\bullet}(V).$$

One checks easily that the action of $\mathrm{Hom}(V, k^j)$ respects this filtration and moreover the action of $\mathrm{Hom}(V, k^j)$ on subsequent factors of this filtration $\Phi_t / \Phi_{t-1} = \bigoplus_{|i_\bullet| = t} S^{m_\bullet - i_\bullet}(k^j) \otimes S^{i_\bullet}(V)$ is trivial. Consider the spectral sequence defined by the above filtration:

$$E_1^{s,t} = H^s(\mathrm{Hom}(V, k^j), k) \otimes \bigoplus_{\substack{0 \leq i_\bullet \leq m_\bullet \\ |i_\bullet| = t}} S^{m_\bullet - i_\bullet}(k^j) \otimes S^{i_\bullet}(V) \Rightarrow H^{s+t}_{S^{m_\bullet},k}(V).$$

Observe finally that the above filtration and hence also the above spectral sequence depend functorially on $V$. In particular $E_r^{s,t} \in \mathcal{F}$ for all $s, t, r$. A well-known computation of cohomology of a finite dimensional vector $k$-space with trivial coefficients [C] implies that the functors $E_1^{s,t}$ are finite. Immediate induction on $r$ shows further that all functors $E_r^{s,t}$, $1 \leq r \leq \infty$, are finite. Finally, the functor $H^i_{S^{m_\bullet},j}$ has a finite filtration all subsequent factors of which are finite and hence is finite as well.

For a general finite $A$, one can find a resolution [K1, Th. 5.1]

$$0 \to A \to A^0 \to A^1 \to \ldots$$

in which every functor $A^i$ is a finite direct sum of functors of the form $S^{m_\bullet}$. To conclude the proof consider the spectral sequence defined by this resolution

$$E_1^{s,t} = H^s_{A^t,k}(V) \Rightarrow H^{s+t}_{A,k}(V)$$

(which depends functorially on $V$) and apply the same argument as above. □

COROLLARY A.6. *Let $B \in \mathcal{F}$ be a finite functor. Then for every $s, j \geq 0$ there exists $N(s, j)$ such that for $n \geq N(s, j)$*

$$H^i\left(\begin{pmatrix} 1_j & * \\ 0 & \mathrm{GL}_{n-j}(k) \end{pmatrix}, B(k^n)\right) = \begin{cases} B(k^j) & \text{if } i = 0 \\ 0 & \text{if } 0 < i \leq s. \end{cases}$$

*Proof.* Consider the Hochschild-Serre spectral sequence corresponding to the group extension

$$1 \to \operatorname{Hom}(k^{n-j}, k^j) \to \begin{pmatrix} 1_j & * \\ 0 & \operatorname{GL}_{n-j}(k) \end{pmatrix} \to \operatorname{GL}_{n-j}(k) \to 1 \;,$$

$$E_2^{a,b} = H^a(\operatorname{GL}_{n-j}(k), H^b_{B,j}(k^{n-j})) \Rightarrow H^{a+b}(\begin{pmatrix} 1_j & * \\ 0 & \operatorname{GL}_{n-j}(k) \end{pmatrix}, B(k^n)) \;,$$

and apply Lemma A.5 and Proposition A.4. □

Combining Lemma A.3 with Corollary A.5, we immediately conclude the following:

COROLLARY A.7. *Let $B$ be a finite functor. Then for any $V \in \mathcal{V}$ and any $s \geq 0$, there exists $N(V, s)$ such that for $n \geq N(V, s)$*

$$\operatorname{Ext}^i_{\operatorname{GL}_n(k)}(P_V(k^n), B(k^n)) = \begin{cases} \bigoplus_{W \subset V} B(V/W) & \text{if } i = 0 \\ 0 & \text{if } 0 < i \leq s. \end{cases}$$

Note that unlike the situation of Lemma A.3, the isomorphism

$$\operatorname{Hom}_{\operatorname{GL}_n(k)}(P_V(k^n), B(k^n)) = \bigoplus_{W \subset V} B(V/W)$$

is now canonical: it associates to a family $\{b_W \in B(V/W)\}_{W \subset V}$ a $\operatorname{GL}_n(k)$-equivariant homomorphism

$$P_V(k^n) = k[\operatorname{Hom}(V, k^n)] \to B(k^n)$$

which sends $f : V \to k^n$ to $B(\overline{f})(b_{\operatorname{Ker} f})$, where $\overline{f} : V/\operatorname{Ker} f \to k^n$ is the induced homomorphism.

If $A, B \in \mathcal{F}$ are functors for which the sequence $\operatorname{Ext}^s_{\operatorname{GL}_n(k)}(A(k^n), B(k^n))$ stabilizes with the growth of $n$, then we use the notation $\operatorname{Ext}^s_{\operatorname{GL}(k)}(A, B)$ for the stable value of this Ext-group. Thus Corollary A.7 tells us, in particular, that

$$\operatorname{Hom}_{\operatorname{GL}(k)}(P_V, B) = \bigoplus_{W \subset V} B(V/W) = \prod_{W \subset V} B(V/W)$$

for a finite $B$. Note that every $k$-linear homomorphism $V \to V'$ defines a homomorphism of functors $P_{V'} \to P_V$ and hence induces also a homomorphism

$$\operatorname{Hom}_{\operatorname{GL}(k)}(P_V, B) \to \operatorname{Hom}_{\operatorname{GL}(k)}(P_{V'}, B).$$

This makes $\operatorname{Hom}_{\operatorname{GL}(k)}(P_V, B) = \bigoplus_{W \subset V} B(V/W)$ into a functor from $\mathcal{V}^f$ to itself. One checks easily that the map $\bigoplus_{W \subset V} B(V/W) \to \bigoplus_{W' \subset V'} B(V'/W')$ corresponding to a $k$-linear homomorphism $f : V \to V'$ looks as follows: the $W'$-component of the image is obtained from the $f^{-1}(W')$-component of the pre-image by applying the homomorphism $B(\overline{f}) : B(V/f^{-1}(W')) \to B(V'/W')$ where $\overline{f} : V/f^{-1}(W') \to V'/W'$ is the homomorphism induced by $f$.

More generally, for any (not necessarily finite) functor $B \in \mathcal{F}$ we can define a new functor $aB \in \mathcal{F}$, by setting

$$(aB)(V) = \bigoplus_{W \subset V} B(V/W)$$

and defining homomorphisms $(aB)(f)$ (for $f : V \to V'$) in the same way as above. Obviously we get in this way an exact functor $a : \mathcal{F} \to \mathcal{F}$. By the universal property of $P_V \in \mathcal{F}$, we conclude that for any finite functor $B$ we have a natural isomorphism $\mathrm{Hom}_{\mathrm{GL}}(P_V, B) = \mathrm{Hom}_{\mathcal{F}}(P_V, aB)$.

Moreover, for any $B \in \mathcal{F}$ there is a natural monomorphism $B \hookrightarrow aB$, mapping $B(V)$ diagonally to $(aB)(V) = \bigoplus_{W \subset V} B(V/W)$. Define $\widetilde{a}B = aB/B$. This gives us another exact functor from $\mathcal{F}$ to itself. Except for the trivial case $B = 0$, the functor $\widetilde{a}B$ is never finite.

THEOREM A.8. *For any finite functors $A, B \in \mathcal{F}$ there is a natural identification*

$$\mathrm{Ext}^*_{\mathrm{GL}(k)}(A, B) = \mathrm{Ext}^*_{\mathcal{F}}(A, aB).$$

*After this identification, the natural homomorphism $\mathrm{Ext}^*_{\mathcal{F}}(A, B) \to \mathrm{Ext}^*_{\mathrm{GL}(k)}(A, B)$ corresponds to the homomorphism $\mathrm{Ext}^*_{\mathcal{F}}(A, B) \to \mathrm{Ext}^*_{\mathcal{F}}(A, aB)$ defined by the natural embedding $B \hookrightarrow aB$.*

*Proof.* Choose a resolution of $A$:

$$0 \leftarrow A \leftarrow P_0 \leftarrow P_1 \leftarrow \cdots$$

in which each $P_i$ is a finite direct sum of functors of the form $P_V$ and consider the spectral sequence

$$E_1^{i,j} = \mathrm{Ext}^i_{\mathrm{GL}_n(k)}(P_j(k^n), B(k^n)) \Rightarrow \mathrm{Ext}^{i+j}_{\mathrm{GL}_n(k)}(A(k^n), B(k^n)).$$

Corollary A.7 shows that for $n$ big enough with respect to $s$, all $E_1^{i,j}$-terms with $0 < i \leq s, 0 \leq j \leq s$ are zero. Thus for $n$ big enough, $\mathrm{Ext}^s_{\mathrm{GL}_n(k)}(A(k^n), B(k^n))$ coincides with the $s^{\mathrm{th}}$ homology group of the complex

$$0 \to \mathrm{Hom}_{\mathrm{GL}_n(k)}(P_0(k^n), B(k^n)) = \mathrm{Hom}_{\mathcal{F}}(P_0, aB)$$
$$\to \mathrm{Hom}_{\mathrm{GL}_k(F)}(P_1(k^n), B(k^n)) = \mathrm{Hom}_{\mathcal{F}}(P_1, aB) \to \ldots ,$$

i.e., coincides with $\mathrm{Ext}^s_{\mathcal{F}}(A, aB)$. $\square$

COROLLARY A.9. *For finite functors $A, B \in \mathcal{F}$ the following conditions are equivalent*:

(1) *The natural homomorphism* $\operatorname{Ext}^*_{\mathcal{F}}(A, B) \to \operatorname{Ext}^*_{\operatorname{GL}(k)}(A, B)$ *is an isomorphism.*
(2) *The natural homomorphism* $\operatorname{Ext}^*_{\mathcal{F}}(A, B) \to \operatorname{Ext}^*_{\mathcal{F}}(A, aB)$ *is an isomorphism.*
(3) $\operatorname{Ext}^*_{\mathcal{F}}(A, \widetilde{a}B) = 0.$

We need a natural generalization of the functor $a : \mathcal{F} \to \mathcal{F}$ to the case of bi-functors (and more generally poly-functors). Denote by $\mathcal{F}_n$ the category whose objects are covariant functors of $n$ variables $\mathcal{V}^f \times \cdots \times \mathcal{V}^f \to \mathcal{V}$ (and whose morphisms are natural transformations). For every functor $B \in \mathcal{F}_n$ define a new functor $a_n(B) \in \mathcal{F}_n$ using the formula

$$(a_n B)(V_1, \ldots, V_n) = \bigoplus_{W_1 \subset V_1 \cdots W_n \subset V_n} B(V_1/W_1, \ldots, V_n/W_n)$$

and defining $a_n B$ on morphisms in the same way as above. Denote further by

$$\mathcal{V}^f \xrightarrow{\Delta} \mathcal{V}^f \times \cdots \times \mathcal{V}^f$$

the diagonal and the direct sum functors.

PROPOSITION A.10. *For any functors $A \in \mathcal{F}, B \in \mathcal{F}_n$ there is a natural split monomorphism*

$$\operatorname{Ext}^*_{\mathcal{F}}(A, a(B \circ D)) \hookrightarrow \operatorname{Ext}^*_{\mathcal{F}_n}(A \circ \Pi, a_n B).$$

*Proof.* The standard adjunction formula (1.7.1) shows that

$$\operatorname{Ext}^*_{\mathcal{F}_n}(A \circ \Pi, a_n B) = \operatorname{Ext}^*_{\mathcal{F}}(A, (a_n B) \circ D).$$

Now it suffices to note that there are homomorphisms (natural in $V$)

$$(a(B \circ D))(V) = \bigoplus_{W \subset V} B(V/W, \ldots, V/W) \underset{p}{\overset{i}{\rightleftarrows}} ((a_n B) \circ D)(V)$$
$$= \bigoplus_{W_1 \subset V, \ldots, W_n \subset V} B(V/W_1, \ldots, V/W_n).$$

Here the composition of the $(W_1, \ldots, W_n)$-indexed projection on the right with $i$ equals the $W_1 \cap \cdots \cap W_n$-indexed projection on the left followed by the natural map

$$B(V/W_1 \cap \cdots \cap W_n, \ldots, V/W_1 \cap \cdots \cap W_n) \to B(V/W_1, \ldots, V/W_n) \ ;$$

and the composition of the $W$-indexed projection on the left with $p$ equals the $(W, \ldots, W)$-indexed projection on the right. Obviously $pi = 1$ so that $a(B \circ D)$ is canonically a direct summand in $(a_n B) \circ D$. $\square$

The following lemma is obvious from the construction

LEMMA A.11. *For any $B_1, \ldots, B_n \in \mathcal{F}$ there is a natural isomorphism of functors*

$$a_n(B_1 \boxtimes \cdots \boxtimes B_n) = aB_1 \boxtimes \cdots \boxtimes aB_n.$$

COROLLARY A.12. *Let $A^\bullet \in \mathcal{F}$ be a (graded) exponential functor. Then for any functors $B_1, \ldots, B_n \in \mathcal{F}$ and any $s \geq 0$ we have a natural split monomorphism*

$$\operatorname{Ext}_{\mathcal{F}}^*(A^s, a(B_1 \otimes \cdots \otimes B_n)) \hookrightarrow \bigoplus_{s_1+\cdots+s_n=s} \bigotimes_{i=1}^n \operatorname{Ext}^*(A^{s_i}, a(B_i)).$$

*Proof.* This follows immediately from Proposition A.10, Lemma A.11 and the Kunneth formula (1.7.2). □

COROLLARY A.13. *With conditions and notation of Corollary* A.12 *assume in addition that the functors $B_i$ are without constant term (i.e. $B_i(0) = 0$). Assume further that $n > 1$ and the natural homomorphisms $\operatorname{Ext}_{\mathcal{F}}^*(A^{s'}, B_i) \to \operatorname{Ext}_{\mathcal{F}}^*(A^{s'}, aB_i)$ are isomorphisms for all $i$ and all $s' < s$. Then the homomorphism*

$$\operatorname{Ext}_{\mathcal{F}}^*(A^s, B_1 \otimes \cdots \otimes B_n) \to \operatorname{Ext}_{\mathcal{F}}^*(A^s, a(B_1 \otimes \cdots \otimes B_n))$$

*is an isomorphism.*

*Proof.* This follows immediately from the commutative diagram

$$\operatorname{Ext}_{\mathcal{F}}^*(A^s, B_1 \otimes \cdots \otimes B_n) @>\sim>> \bigoplus_{s_1+\cdots+s_n=s} \bigotimes_{i=1}^n \operatorname{Ext}_{\mathcal{F}}^*(A^{s_i}, B_i)$$
$$@VVV @VV\cong V$$
$$\operatorname{Ext}_{\mathcal{F}}^*(A^s, a(B_1 \otimes \cdots \otimes B_n)) @>\hookrightarrow>> \bigoplus_{s_1+\cdots+s_n=s} \bigotimes_{i=1}^n \operatorname{Ext}^*(A^{s_i}, aB_i).$$

The right-hand vertical arrow in this diagram is an isomorphism since summands, corresponding to $n$-tuples $(s_1, \ldots, s_n)$ with $s_i = 0$ for some $i$, are trivial (the functors $B_i$ and $aB_i$ being without constant term). □

Now we are prepared to start the proof of Theorem A.1. We begin with a special case when both functors $A$ and $B$ are additive.

PROPOSITION A.14. *The natural homomorphism*

$$\operatorname{Ext}_{\mathcal{F}}^*(I^{(s)}, I) \to \operatorname{Ext}_{\operatorname{GL}(k)}^*(I^{(s)}, I) = \operatorname{Ext}_{\mathcal{F}}^*(I^{(s)}, aI)$$

*is an isomorphisms for all $s \geq 0$.*

*Proof.* This follows immediately from a theorem of Dundas and McCarthy [D-M]. Alternatively one can proceed as follows. We may clearly assume that $0 \leq s < r$ ( where $q = p^r$). Both sides are zero (by weight considerations) unless $s = 0$. So assume that $s = 0$ and consider the following commutative diagram

$$\varinjlim_{m \geq 0} \mathrm{Ext}^*_{\mathcal{P}}(I^{(m)}, I^{(m)}) @>\sim>> \mathrm{Ext}^*_{\mathrm{GL}(k)}(I, I)$$
$$@V \cong VV @V = VV$$
$$\mathrm{Ext}^*_{\mathcal{F}}(I, I) @>>> \mathrm{Ext}^*_{\mathrm{GL}(k)}(I, I).$$

The top horizontal arrow is an isomorphism; this is just a reformulation of Theorem 7.6 [F-S]. The left vertical arrow is an isomorphism; this follows from our strong comparison theorem (Theorem 3.10) or directly from computation of both groups presented in [F-L-S] and [F-S] respectively. Thus the bottom horizontal arrow is an isomorphism as well. □

PROPOSITION A.15. *The natural homomorphisms*

$$\mathrm{Ext}^*_{\mathcal{F}}(I^{(s)}, S^n) \to \mathrm{Ext}^*_{\mathcal{F}}(I^{(s)}, aS^n)$$

$$\mathrm{Ext}^*_{\mathcal{F}}(I^{(s)}, \Lambda^n) \to \mathrm{Ext}^*_{\mathcal{F}}(I^{(s)}, a\Lambda^n)$$

*are isomorphisms for all $n$ and $s$.*

*Proof.* We proceed by induction on $n$, using Proposition A.14 as the induction base. Applying the exact functor $\widetilde{a}$ to the Koszul and de Rham complexes we get the following two complexes:

$$0 \to \widetilde{a}\Lambda^n \to \widetilde{a}(\Lambda^{n-1} \otimes S^1) \to \cdots \to \widetilde{a}(\Lambda^1 \otimes S^{n-1}) \to \widetilde{a}S^n \to 0 ,$$

$$0 \to \widetilde{a}S^n \to \widetilde{a}(S^{n-1} \otimes \Lambda^1) \to \cdots \to \widetilde{a}(S^1 \otimes \Lambda^{n-1}) \to \widetilde{a}\Lambda^n \to 0.$$

The first of these two complexes is always acyclic, the second one is acyclic if $n \not\equiv 0 \mod (p)$ and has homology isomorphic to $\widetilde{a}((S^{m-i} \otimes \Lambda^i)^{(1)}) = (\widetilde{a}(S^{m-i} \otimes \Lambda^i))^{(1)}$ if $n = pm$. The induction hypothesis and Corollary A.13 show that $\mathrm{Ext}^*_{\mathcal{F}}(I^{(s)}, (\widetilde{a}(S^{m-i} \otimes \Lambda^i))^{(1)}) = 0$. Thus the second hypercohomology spectral sequences in both cases consist of zeroes only and hence the limit of the first spectral sequences is zero. On the other hand, Corollary A.13 implies that the $E_1$-term of both first hypercohomology spectral sequences could have only two nonzero columns and hence the differentials

$$d_n : \mathrm{Ext}^*_{\mathcal{F}}(I^{(s)}, \widetilde{a}(\Lambda^n)) \to \mathrm{Ext}^{*-(n-1)}_{\mathcal{F}}(I^{(s)}, \widetilde{a}(S^n)),$$

$$d'_n : \mathrm{Ext}^*_{\mathcal{F}}(I^{(s)}, \widetilde{a}(S^n)) \to \mathrm{Ext}^{*-(n-1)}_{\mathcal{F}}(I^{(s)}, \widetilde{a}(\Lambda^n))$$

are isomorphisms. Using an additional induction on the cohomology index we easily conclude the proof. □

Applying the duality isomorphisms of Lemma A.2 we derive from Proposition A.15 the following:

COROLLARY A.16. *The natural homomorphisms*

$$\operatorname{Ext}^*_{\mathcal{F}}(\Gamma^n, I^{(s)}) \to \operatorname{Ext}^*_{\mathcal{F}}(\Gamma^n, aI^{(s)}),$$
$$\operatorname{Ext}^*_{\mathcal{F}}(\Lambda^n, I^{(s)}) \to \operatorname{Ext}^*_{\mathcal{F}}(\Lambda^n, aI^{(s)})$$

*are isomorphisms for all $n, s$.*

PROPOSITION A.17. *The natural homomorphisms*

$$\operatorname{Ext}^*_{\mathcal{F}}(\Gamma^{m(s)}, S^n) \to \operatorname{Ext}^*_{\mathcal{F}}(\Gamma^{m(s)}, aS^n),$$
$$\operatorname{Ext}^*_{\mathcal{F}}(\Gamma^{m(s)}, \Lambda^n) \to \operatorname{Ext}^*_{\mathcal{F}}(\Gamma^{m(s)}, a\Lambda^n)$$

*are isomorphisms for all $m, n, s$.*

*Proof.* We use induction on $n$ and repeat word by word the argument used in the proof of the Proposition A.15, using now Corollary A.16 as the induction base. □

COROLLARY A.18. *The natural homomorphism*

$$\operatorname{Ext}^*_{\mathcal{F}}(\Gamma^{m_1} \otimes \cdots \otimes \Gamma^{m_s}, S^{n_1} \otimes \cdots \otimes S^{n_t}) \to \operatorname{Ext}^*_{\mathcal{F}}(\Gamma^{m_1} \otimes \cdots \otimes \Gamma^{m_s}, a(S^{n_1} \otimes \cdots \otimes S^{n_t}))$$

*is an isomorphism for all $m_1, \ldots, m_s, n_1, \ldots, n_t$.*

*Proof.* Assume first that $s = 1$. In this case our statement follows from Corollary A.13 and Proposition A.17. Lemma A.2 shows further that the statement holds also (for any $s$) in case $t = 1$. We may now obtain the general case in the same way as above, applying Corollary A.13 to the exponential functor $\underbrace{\Gamma^\bullet \otimes \cdots \otimes \Gamma^\bullet}_{s}$. □

*End of the Proof of Theorem* A.1. We can construct resolutions

$$0 \to B \to B^0 \to B^1 \to \ldots,$$
$$0 \leftarrow A \leftarrow A_0 \leftarrow A_1 \leftarrow \cdots$$

in which all $B^i$ are finite direct sums of functors of the form $S^{n_1} \otimes \cdots \otimes S^{n_t}$ and all $A_i$ are finite direct sums of functors of the form $\Gamma^{m_1} \otimes \cdots \otimes \Gamma^{m_s}$. Comparing the corresponding hypercohomology spectral sequences and using Corollary A.18 we get the desired result. □

The following is an immediate corollary of Theorem A.1 and the explicit stability bound of van der Kallen [vdK].

COROLLARY A.19. *Let $k$ be a finite field and let $P, Q \in \mathcal{F} = \mathcal{F}(k)$ be finite functors each of degree $\leq d$. Then the natural map*

$$\operatorname{Ext}^m_{\mathcal{F}}(F, Q) \to \operatorname{Ext}^m_{\operatorname{GL}_n(k)}(P(k^n), Q(k^n))$$

*is an isomorphism whenever $n \geq 2m + 2d$.*

UNIVERSITY OF NANTES, NANTES, FRANCE
*E-mail address*: franjou@math.univ-nantes.fr

NORTHWESTERN UNIVERSITY, EVANSTON, IL
*E-mail addresses*: eric@math.nwu.edu
scor@math.nwu.edu
suslin@math.nwu.edu